\documentclass[10pt]{article}
\date{29.5.2019}
\usepackage[english]{babel}
\usepackage{hyperref}
\usepackage{amssymb, amsfonts}
\usepackage{stmaryrd}
\usepackage{mathtools}
\usepackage{enumerate}
\usepackage{amsmath}
\usepackage{latexsym}
\usepackage{amssymb}

\usepackage{amsthm}

\font\tengoth=eufm10 at 10pt
\font\sevengoth=eufm7 at 6pt
\newfam\gothfam
\textfont\gothfam=\tengoth
\scriptfont\gothfam=\sevengoth

\newcommand{\mlabel}[1]{\marginpar{#1}\label{#1}}

\newcommand{\fS}{{\mathfrak S}}

\newcommand{\g}{{\mathfrak g}}
\newcommand{\s}{{\mathfrak s}}

\newcommand{\fh}{{\mathfrak h}}

\newcommand{\fu}{{\mathfrak u}}

\renewcommand\sp{\mathfrak {sp}}

\renewcommand{\:}{\colon}
\newcommand{\1}{\mathbf{1}}

\newcommand{\cA}{\mathcal{A}}
\newcommand{\cB}{\mathcal{B}}
\newcommand{\cC}{\mathcal{C}}
\newcommand{\cD}{\mathcal{D}}
\newcommand{\cE}{\mathcal{E}}
\newcommand{\cF}{\mathcal{F}}

\newcommand{\cH}{\mathcal{H}}
\newcommand{\cK}{\mathcal{K}}
\newcommand{\cL}{\mathcal{L}}
\newcommand{\cM}{\mathcal{M}}
\newcommand{\cN}{\mathcal{N}}

\newcommand{\cS}{\mathcal{S}}

\newcommand{\cU}{\mathcal{U}}

\newcommand{\cZ}{\mathcal{Z}}

\newcommand{\eset}{\emptyset}

\renewcommand{\phi}{\varphi}
\newcommand{\dd}{{\tt d}}

\newcommand{\subeq}{\subseteq}
\newcommand{\supeq}{\supseteq}

\newcommand{\into}{\hookrightarrow}
\newcommand{\eps}{\varepsilon}

\newcommand{\shalf}{{\textstyle{\frac{1}{2}}}}

\newcommand{\N}{{\mathbb N}}
\newcommand{\Z}{{\mathbb Z}}
\newcommand{\R}{{\mathbb R}}
\newcommand{\C}{{\mathbb C}}

\newcommand{\Q}{{\mathbb Q}}

\newcommand{\T}{{\mathbb T}}

\newcommand{\bH}{{\mathbb H}}
\newcommand{\bS}{{\mathbb S}}

\renewcommand{\hat}{\widehat}

\renewcommand{\tilde}{\widetilde}

\newcommand{\SU}{\mathop{{\rm SU}}\nolimits}

\newcommand{\U}{\mathop{\rm U{}}\nolimits}
\newcommand{\PU}{\mathop{\rm PU{}}\nolimits}

\newcommand{\Sp}{\mathop{{\rm Sp}}\nolimits}

\newcommand{\ad}{\mathop{{\rm ad}}\nolimits}
\newcommand{\Ad}{\mathop{{\rm Ad}}\nolimits}

\renewcommand{\Im}{\mathop{{\rm Im}}\nolimits}

\newcommand{\tr}{\mathop{{\rm tr}}\nolimits}

\newcommand{\Hom}{\mathop{{\rm Hom}}\nolimits}

\newcommand{\Ext}{\mathop{{\rm Ext}}\nolimits}

\newcommand{\Heis}{\mathop{{\rm Heis}}\nolimits}
\newcommand{\Homeo}{\mathop{{\rm Homeo}}\nolimits}

\newcommand{\Aut}{\mathop{{\rm Aut}}\nolimits}

\newcommand{\diag}{\mathop{{\rm diag}}\nolimits}

\newcommand{\id}{\mathop{{\rm id}}\nolimits}

\renewcommand{\dim}{\mathop{{\rm dim}}\nolimits}

\newcommand{\Ker}{\mathop{{\rm Ker}}\nolimits}

\newcommand{\supp}{\mathop{{\rm supp}}\nolimits}

\newcommand{\Inn}{\mathop{{\rm Inn}}\nolimits}

\newcommand{\conv}{\mathop{{\rm conv}}\nolimits}

\newcommand{\Spann}{\mathop{{\rm span}}\nolimits}

\newcommand{\Rarrow}{\Rightarrow}
\newcommand{\nin}{\noindent}
\newcommand{\oline}{\overline}

\newcommand{\la}{\langle}
\newcommand{\ra}{\rangle}

\newcommand{\res}{\vert}

\newcommand{\Spec}{{\rm Spec}}

\newcommand{\ssssarr}{\hbox to 15pt{\rightarrowfill}}
\newcommand{\sssarr}{\hbox to 20pt{\rightarrowfill}}
\newcommand{\ssarr}{\hbox to 30pt{\rightarrowfill}}
\newcommand{\sarr}{\hbox to 40pt{\rightarrowfill}}
\newcommand{\arr}{\hbox to 60pt{\rightarrowfill}}
\newcommand{\larr}{\hbox to 60pt{\leftarrowfill}}
\newcommand{\Arr}{\hbox to 80pt{\rightarrowfill}}

\def\theoremname{Theorem}
\def\propositionname{Proposition}
\def\corollaryname{Corollary}
\def\lemmaname{Lemma}
\def\remarkname{Remark}
\def\conjecturename{Conjecture}

\def\definitionname{Definition}
\def\exercisename{Exercise}
\def\examplename{Example}
\def\examplesname{Examples}
\def\problemname{Problem}
\def\problemsname{Problems}
\def\proofname{Proof}

\def\satzname{Satz}
\def\koroname{Korollar}
\def\folgname{Folgerung}
\def\bemerkname{Bemerkung}
\def\aufgname{Aufgabe}

\def\beisname{Beispiel}
\def\beissname{Beispiele}
\def\bewname{Beweis}

\def\@thmcounter#1{\noexpand\arabic{#1}}
\def\@thmcountersep{}
\def\@begintheorem#1#2{\it \trivlist \item[\hskip
\labelsep{\bf #1\ #2.\quad}]}
\def\@opargbegintheorem#1#2#3{\it \trivlist
      \item[\hskip \labelsep{\bf #1\ #2.\quad{\rm #3}}]}
\makeatother
\newtheorem{theor}{\theoremname}[section]
\newtheorem{propo}[theor]{\propositionname}
\newtheorem{coro}[theor]{\corollaryname}
\newtheorem{lemm}[theor]{\lemmaname}

\newenvironment{thm}{\begin{theor}\it}{\end{theor}}

\newenvironment{prop}{\begin{propo}\it}{\end{propo}}

\newenvironment{cor}{\begin{coro}\it}{\end{coro}}

\newenvironment{lem}{\begin{lemm}\it}{\end{lemm}}

\newenvironment{Lemma}{\begin{lemm}\it}{\end{lemm}}

\newtheorem{rema}[theor]{\remarkname}

\newenvironment{rem}{\begin{rema}\rm}{\end{rema}}

\newtheorem{stepnow}[theor]{}

\newtheorem{defin}[theor]{\definitionname} 

\newenvironment{defn}{\begin{defin}\rm}{\end{defin}}

\newtheorem{exerc}[theor]{\exercisename}

\newtheorem{exa}[theor]{\examplename}

\newenvironment{example}{\begin{exa}\rm}{\end{exa}}
\newenvironment{ex}{\begin{exa}\rm}{\end{exa}}

\newtheorem{exas}[theor]{\examplesname}

\newenvironment{exs}{\begin{exas}\rm}{\end{exas}}

\newtheorem{conj}[theor]{\conjecturename}

\newtheorem{pro}[theor]{\problemname}

\newtheorem{prs}[theor]{\problemsname}

\newenvironment{Proof*}{\begin{trivlist}\item[\hskip%
\labelsep{\bf\proofname.\quad}]}%
{\end{trivlist}}

\newenvironment{prf}{\begin{proof}}{\end{proof}}

{\hfill\qed\end{trivlist}}

\newenvironment{beweis*}{\begin{trivlist}\item[\hskip%
\labelsep{\bf\bewname.\quad}]}%
{\end{trivlist}}

\newtheorem{satzn}[theor]{\satzname}

\newtheorem{koro}[theor]{\koroname}

\newtheorem{folg}[theor]{\folgname}

\newtheorem{bem}[theor]{\bemerkname}

\newtheorem{aufg}[theor]{\aufgname}

\newtheorem{aufgn}[theor]{\aufgname}

\newtheorem{beis}[theor]{\beisname}

\newtheorem{beiss}[theor]{\beissname}

\DeclareMathAlphabet{\Ma}{U}{msa}{m}{n}
\DeclareMathAlphabet{\Mb}{U}{msb}{m}{n}
\DeclareMathAlphabet{\Meuf}{U}{euf}{m}{n}

\def\got#1{\Meuf{#1}}
\def\br#1.{\llbracket #1 \rrbracket}
\def\lbr{\llbracket}
\def\rbr{\rrbracket}
\def\strong{point-norm }
\def\strongH{strong operator }
\def\usual{regular }
\def\minimal{inner minimal positive }

\newcommand{\cJ}{\mathcal{J}}
\newcommand{\Fol}{\mathop{{\rm Fol}}\nolimits}
\newcommand{\Proj}{\mathop{{\rm Proj}}\nolimits}
\renewcommand{\phi}{\varphi}
\newcommand{\fA}{{\mathfrak A}}

\renewcommand{\baselinestretch}{1.2}
\normalsize
 \def\ot #1.{{\got{#1}}}

\def\cR{\mathcal{R}}
\def\al#1.{{\cal #1}}
\def\wt{\widetilde}

\def\s #1.{_{\smash{\lower2pt\hbox{\mathsurround=0pt $\scriptstyle #1$}}\mathsurround=5pt}}
\def\XP#1!{\renewcommand{\baselinestretch}{.7}\marginpar{{\footnotesize
$\leftarrow$#1}\hfil}\renewcommand{\baselinestretch}{1.2}}
\def\ccr #1,#2.{\overline{\Delta(#1,\,#2)}}
\def\rsl{\mathord{\al R.(\cH,\sigma)}}
\def\b#1.{{\bf #1}}

\newcommand{\be}{{\bf{e}}}

\renewcommand\mlabel{\label}

\begin{document}


\title{Covariant representations for possibly singular actions on $C^*$-algebras}
\author{
Daniel Belti\c t\u a\footnote{Institute of Mathematics ``Simion Stoilow'' of the Romanian Academy,
		P.O. Box 1-764, Bucharest, Romania.
		\texttt{Email:} beltita@gmail.com, Daniel.Beltita@imar.ro},
Hendrik Grundling\footnote{Department of Mathematics,
	University of New South Wales,
	Sydney, NSW 2052, Australia.
	\texttt{Email:} h.grundling@unsw.edu.au},
Karl-Hermann Neeb\footnote{Department of Mathematics,
	Friedrich-Alexander-Universit\"at Erlangen-N\"urnberg,
	Cauerstr. 11, 91058 Erlangen, Germany.
	\texttt{Email:} neeb@mi.uni-erlangen.de}}

\maketitle

\begin{abstract}
Singular actions on $C^*$-algebras are automorphic group actions on $C^*$-algebras,
where the group is not locally compact, or the action is not  strongly continuous.
We study the covariant representation theory of actions which may be singular.
In the usual case of strongly continuous  actions of locally compact groups on $C^*$-algebras,
this is done via crossed products, but this approach is not available for  singular $C^*$-actions.
We explored extension of  crossed products to singular actions in a previous paper.
The literature regarding covariant representations
for possibly singular actions is already large and scattered, and in need of some consolidation.
We collect in this survey  a range of results in this field, mostly known. We  improve some proofs and
elucidate some interconnections. These include existence theorems by Borchers and Halpern,
Arveson spectra, the Borchers--Arveson theorem, standard representations and Stinespring dilations
as well as ground states, KMS states and ergodic states and the spatial structure of their GNS representations.
\end{abstract}

\tableofcontents

\section{Introduction}

Covariant representations of $C^*$- and $W^*$-dynamical systems $(\cA, G, \alpha)$
are fundamental
objects in both $C^*$-algebra theory, as well as in mathematical quantum physics.
Our interest here is in covariant representations for possibly singular $C^*$-actions,
i.e.\ automorphic group actions on $C^*$-algebras, where the group need not
be locally compact, or the action need not be strongly continuous (i.e.
continuous w.r.t. the pointwise convergence topology).
Such actions are abundant in physics and arise naturally in
mathematics.
For example, for bosonic field theories, the field $C^*$-algebra is usually
chosen to be either the Weyl algebra, or the resolvent algebra, and then
nonconstant one parameter symplectic groups
produce one parameter automorphism groups on these algebras which are not
strongly continuous (cf.\ Examples~\ref{ccrexmp1} and \ref{ccrexmp2} below).
On the other hand, for a gauge theory, the gauge group
has to act on the field algebra,  and this is infinite dimensional,
hence not locally compact. Any unitary representation of the gauge
group can lead to a singular action, either on the CAR algebra,
or on the Weyl algebra associated with the Hilbert space.
Other naturally occurring actions of infinite dimensional Lie groups
in physics, are loop groups, restricted orthogonal and symplectic
groups on a Hilbert space, the diffeomorphism group of the circle,
and other diffeomorphism groups in gravity models.

Many of the usual mathematical tools break down
for singular actions, e.g.\ $C^*$-crossed products,
which means that there is not a good structure theory
for their covariant representation theory,
but a great deal of analysis has  been done for special
subsets of it. Though much of the theory is collected
in monographs such as \cite{BR02}, \cite{BR96}, \cite{Sa91},
\cite{Pe89} and \cite{Ta03}, unfortunately many important
results are still widely scattered in the literature. We feel it necessary
to  collect here some of these scattered results, improve proofs where we can,
and add some new examples and results which seem interesting.
Our intention is to augment the material in the monographs, not to
replace any of these sources.
Whilst the usefulness of this is primarily for ourselves, we hope that
this review will also be of use to practitioners in the area.

In a previous work, we studied crossed product constructions which are possible
for a subclass of possibly singular actions (cf.~\cite{GrN14}). This does provide a good
structure theory for (a subclass) of their covariant representations.
However in the interest of brevity, in this review we will not include these.
We will mainly concentrate on singular actions and their associated $W^*$-dynamical
systems in specific covariant representations. Whilst for a $W^*$-dynamical
system of a locally compact group we can construct a $W^*$-crossed product,
the structure theory this gives is limited to the structure theory possible
on the predual of  a von Neumann algebra. In particular, it can only produce
covariant representations which are normal w.r.t. the defining representation
of the $W^*$-dynamical system.
Whilst the strong operator dense (strongly continuous)
$C^*$-dynamical subsystem of the $W^*$-dynamical system
does have different covariant representations, these need
not extend to the original $C^*$-algebra on which the singular action is defined.
Because of these considerations, we will not study  $W^*$-crossed products here.

For singular actions, we  focus on structural issues for covariant representations,
leaving applications aside. Some of these issues include existence, spectrum
conditions (cf.~Borchers--Arveson Theorem), innerness, standard representation structures and
Stinespring dilations.
The most important types of states associated with a singular action are ground states,
KMS states, and ergodic states, and we will briefly review these, as well
as the properties of their GNS representations.

In more detail, what we will cover are the following.
We start
 with the natural topologies of the automorphism groups
of $C^*$-algebras and von Neumann algebras, and discuss the
Borchers--Halpern Theorem characterizing existence
of covariant representations in terms of their folia
of normal states.
We refine these conditions and
consider covariance for cyclic representations where the generating vector
is not necessarily $G$-invariant.
We  also consider conditions for covariant representations
to be inner.
The universal covariant representation is a
useful tool for analyzing a singular action in terms of a $W^*$-dynamical system.

Next, we  consider the standard form representations  of a
$W^*$-dynamical system, which is a special and heavily used
covariant representation (Section~\ref{SectStandardForm}).
For a projection $P$ in a von Neumann algebra $\cM$, we consider the reduced
von Neumann algebra $P\cM P$, and composing the reduction map with the
standard representation of the image, we obtain a completely positive
map $\varphi_P$ for which we can construct a Stinespring dilation representation
 $\pi_{\varphi_P}$ of $\cM$.
In particular, given a  $W^*$-dynamical
system and an invariant projection $P$, then  $\pi_{\varphi_P}$ is covariant,
and this generalizes the analogous theorem for the GNS representation
of an invariant state for a  $C^*$-dynamical system
(Subsection~\ref{ImpW-dynsys}).

We then consider covariant representations satisfying a spectral condition,
 study issues around the Borchers--Arveson Theorem  (Section~\ref{SpecCovRep})
and characterize the ground
states whose GNS representations give rise to such covariant
representations  (Section~\ref{sec:5}).
This is motivated by the fact that such states are
of central importance in physics, in fact the existence of such an invariant state is an
axiom for algebraic quantum field theory (cf. \cite[Axiom 4, p.104]{Ar99},
\cite{HK64}).
We also consider the structure of
these representations and clarify the role of ground states.
We study the case where zero is isolated in the Arveson spectrum
in detail.

We continue in Section~\ref{KMSsect} by recalling the basic structural
facts of the GNS representation of a KMS state,
since thermal quantum physics is based on such a setting.
This is followed by a very short section on
ergodic states.

\subsection{Notation and terminology}
\label{NotTerm}

For a $C^*$-algebra $\cA$, we write $\cA^*$ for the space of continuous linear functionals on
$\cA$ and $\fS(\cA)\subeq \cA^*$ for the set of states.
For a $W^*$-algebra $\cM$, we write $\cM_* \subeq \cM^*$ for the predual of $\cM$, i.e.~the
subspace of normal functionals and $\fS_n(\cM) \subeq \cM_*$ for the set of normal states.
For a topological group $G$, we write $G_d$ for the underlying discrete group.

If $X$ and $Y$ are Banach spaces, we write $\cB(X,Y)$ for the space of bounded
operators from $X$ to $Y$. Then $\cB(X,Y)$ has two topologies,
the {\it norm topology} (w.r.t. the supremum norm over unit balls), and the {\it strong topology}
which is the topology of pointwise convergence  for maps from $X$ to $Y$.
In the strong topology, an open neighborhood base of an $A\in \cB(X,Y)$ is given by the sets
\[
N_\varepsilon(A;\,x_1,\ldots,x_n):=\{B\in\cB(X,Y)\,\mid\,\|B(x_i)-A(x_i)\|<\varepsilon,\;i=1,\ldots,n\}
\]
for $\varepsilon>0$, $x_i\in X$ and $n\in\N$. The strong topology is also
referred to as the {\it point-norm topology} (\cite[II.5.5.3]{Bla06}),
and we will use this terminology for $\cB(\cA)$ where $\cA$ is a $C^*$-algebra.
The  norm topology (supremum norm over unit balls),
will be the {\it uniform} topology for $\cB(\cA)$.

If $\cH$ is a Hilbert space, then the strong topology of $\cB(\cH)$ coincides with the
strong operator topology. All unitary representations $U:G\to\cU(\cH)$ will be assumed
to be strong operator continuous, including unitary one parameter groups $U:\R\to\cU(\cH)$.
 The cases where continuity is not required, will be
covered by taking the underlying discrete group, i.e. by considering
unitary representations $U:G_d\to\cU(\cH).$

We include an index of terms and notation at the end of this paper.

\section{Covariant representations}
\mlabel{sec:2}

\subsection{$C^*$-and $W^*$-dynamical systems}
\mlabel{subsec:2.1}

For a $C^*$-algebra $\cA,$   as $\Aut(\cA)\subset\cB(\cA),$  there are
two natural topologies for its automorphism group $\Aut(\cA)$
 with respect to  which it  is a topological group.
The  norm topology of $\cB(\cA)\supset\Aut(\cA)$, and the \strong topology.
Therefore if we want a topological group $G$ to act on $\cA$,
it is natural to look for homomorphisms $\alpha \: G \to \Aut(\cA)$, $g\mapsto\alpha_g$,
which are continuous
with respect to  one of these two topologies.
The norm topology is too restrictive for most applications, hence one normally requires
continuity with respect to  the \strong topology. We fix some terminology:
\begin{defn}\mlabel{def:1.1}
\begin{itemize}
\item[(i)]  A (discrete group) {\it $C^*$-action} is a triple $(\cA, G, \alpha)$, where
$\cA$ is a $C^*$-algebra, $G$ is a topological group and
$\alpha \: G \to \Aut(\cA)$ is a homomorphism, which is not assumed to have any continuity
property.  We will usually omit the ``discrete group''.
\item[(ii)]
If $\alpha \: G \to \Aut(\cA)$ is {\it \strong continuous}, i.e.
 for every point $A \in \cA$, the
orbit map $\alpha^A \: G \to \cA,$ $g \mapsto \alpha_g(A)$ is continuous,
we call $(\cA, G, \alpha)$ a {\it $C^*$-dynamical system} (cf.~\cite{Pe89},
\cite[Def.~2.7.1]{BR02}).
The {\it \usual case} will mean that the action is \strong continuous
and the group $G$ is locally compact. A {\it singular} action is one which is not the
\usual case.
\item[(iii)]
A  $C^*$-action $(\cA, G, \alpha)$ has a dual action
$\alpha^*:G\to\cB(\cA^*)$ by isometries on the topological dual $\cA^*$ given by
\begin{equation}
  \label{eq:1a}
(\alpha_g^*\omega)(A) := \omega(\alpha_g^{-1}(A)) \quad \mbox{ for } \quad g \in G,
A \in \cA, \omega \in \cA^*\,.
\end{equation}
 The space of norm continuous elements of $\alpha^*$ is denoted by
\begin{equation}
  \label{eq:a*c}
(\cA^*)_c := \{ \omega \in  \cA^*| {\lim\limits_{g\to \be}\|\alpha_g^*\omega-\omega\|}=0\}.
\end{equation}
Since $G$ acts on $\cA^*$ by isometries, this subspace is norm closed and
maximal with respect to the property that the $G$-action on $(\cA^*)_c$ is continuous
with respect to the norm topology on $(\cA^*)_c$ (see \cite[Thm.~II.2.2]{Bo96} for further
properties). We write
\[ \fS(\cA)_c := \fS(\cA) \cap (\cA^*)_c \]
for the set of states with continuous orbit maps. If $\alpha \: G \to \Aut(\cA)$ is
continuous with respect to the operator norm on $\cB(\cA)$, then $(\cA^*)_c=\cA^*$.
\end{itemize}
\end{defn}

Examples of singular actions were mentioned in the introduction, and below we will
give some typical examples (cf. Examples~\ref{ccrexmp1} and \ref{ccrexmp2}).

The  \usual case, is  what is normally assumed in the literature.
For the \usual case, $C^*$-actions have been extensively analyzed, and there are many tools available,
such as crossed products. However, this is frequently too restrictive, e.g.\ if we have a
\strong continuous one-parameter
automorphism group $\alpha \: \R \to \Aut(\cA)$ where $\cA$ is a $W^*$-algebra, then the action
must be inner (cf. \cite[Exercise~XI.3.6]{Ta03}).
In physics and some natural examples in mathematics, we have
singular actions, and then the available theory is more limited. To analyze a singular action,
one is often forced to choose
some representation $\pi$ with respect to ~which the  $\alpha_g$ are normal maps (i.e.\
each $\alpha_g^*$ preserves the set of normal states of the von Neumann algebra  $\pi(\cA)''$),
and the orbit maps $g \mapsto \pi(\alpha_g(A))$ are strong operator continuous
and then analyze the action on the von Neumann algebra $ \pi(\cA)''$.
The cost of this strategy is that
the analysis is subject to the chosen representation $\pi$.
Not every automorphism of $\pi(\cA)$ will extend to $ \pi(\cA)''$, only those
which are normal maps with respect to ~$\pi$. On the other hand, every automorphism of $ \pi(\cA)''$
is automatically normal,
but not all will preserve $\pi(\cA)$.
We fix terminology for this context.

Let $\cM$ be a $W^*$-algebra, then every automorphism $\rho$ of the
$W^*$-algebra $\cM$ is already a normal map, i.e.\  a $W^*$-automorphism
(cf.~\cite[Thm.~2.5.2]{Pe89} or \cite[Cor.~4.1.23]{Sa71}), hence there is no need to restrict
$\Aut(\cM)$.
As any $\rho\in\Aut(\cM)$ is a normal map, the isometry
$\rho^*:\cM^*\to\cM^*$ (given by $\rho^*(\omega)=\omega\circ\rho$)
preserves the predual $\cM_*$, hence by $\cM=(\cM_*)^*$
the map $\rho\to\rho^*\restriction\cM_*$ embeds $\Aut(\cM)$ as a group of isometries of the Banach space $\cM_*$.

The natural topology one would like to give $\Aut(\cM)$, is the coarsest topology which makes the orbit maps
$ \Aut(\cM) \to \cM,\;\rho\to\rho(A)$ continuous with respect to  any
of the strong operator, weak operator, ultraweak or ultrastrong topologies.
Unfortunately $\Aut(\cM)$ is not a topological group with respect to ~such a topology, which leads us to the following.
As $\Aut(\cM)$ is identified with a group of isometries of $\cM_*$, there are two natural group topologies
on it (cf.~\cite{Haa75}):

\begin{defn}
\label{PUtops}
Let $\cM$ be a $W^*$-algebra. Then the {\it $u$-topology} of $ \Aut(\cM)$ is defined to be
 the coarsest topology which makes the
orbit maps $ \Aut(\cM) \to \cM_*,\;\rho\to \rho^*(\omega)\in\cM_*$ norm continuous for each $\omega\in\cM_*$.
This topology is also called the
$\sigma\hbox{--weak}$ topology (cf.~\cite[p.~12]{Sa91}), and
 $\Aut(\cM)$ is a topological group with respect to ~this topology.\\[2mm]
The {\it $p$-topology}  of $ \Aut(\cM)$
 is the coarsest topology for which all maps $ \Aut(\cM) \to \C,\;
\rho \mapsto \omega(\rho(M))$
for $\omega\in\cM_*$ and $M\in\cM$
are continuous, and this also makes  $\Aut(\cM)$ into a topological group.
\end{defn}
Clearly, the $u$-topology is finer than the $p$-topology, and we will derive the corresponding
inequality in Example~\ref{ex:linft} below.
However, the two topologies coincide for factors of type I and II$_1$ (\cite[Cor.~3.8]{Haa75}).
We define:

\begin{defn}\mlabel{def:1.1.1}
Let $G$ be a topological group and $\cM$ be a  $W^*$-algebra, and assume we have
a homomorphism $\alpha \: G \to \Aut(\cM)$.
We call $(\cM, G, \alpha)$ a {\it $W^*$-dynamical system} if
$\alpha$ is continuous with respect to ~the $u$-topology, i.e.
$\cM_* \subeq (\cM^*)_c$, i.e.\  the action of $G$ on the Banach space
$\cM_*$ is continuous.
\end{defn}
For locally compact groups, this coincides with the naive notion by
the following (\cite[Cor.~2.4]{Hal72}, \cite{Arv74}, \cite[\S 13.5]{Str81},
\cite[Thm.~III.3.2.2]{Bla06}):

\begin{thm}
Let $G$ be a locally compact group, $\cM$ be a  von Neumann algebra, and
\break  $\alpha \: G \to \Aut(\cM)$ a homomorphism. Then the following are equivalent:
 \begin{itemize}
 \item[\rm(i)] For each $M\in\cM$, the map $\alpha^M \: G\to\cM,\; g\mapsto \alpha_g(M)$ is continuous with respect to ~the strong (or weak) operator topology.
 \item[\rm(ii)] For each $\omega\in\cM_*$, the orbit
map  $\alpha^\omega \: G\to\cM_*,\; g\mapsto \alpha_g^*(\omega)$ is norm continuous.
 \item[\rm(iii)] For each $\omega\in\cM_*$  and $M\in\cM$,
the map $\alpha^{\omega,M} \: G\to\C,\;g \mapsto \omega(\alpha_g(M))$ is continuous.
  \end{itemize}
\end{thm}

\begin{rem} \mlabel{rem:2.1}  That (ii) and (iii) need not  be equivalent  for a general topological group
follows from the fact that the  $u$-topology is strictly finer than the $p$-topology for some
von Neumann algebras (cf.~Example~\ref{ex:linft}).
For general topological groups it follows from properties of the
standard representation that this extension of the definition of
a $W^*$-dynamical system is the most useful one (cf.\ equation \eqref{eq:standard-auto} below).
\end{rem}
\begin{example}
  \mlabel{ex:linft} (see~\cite[Cor.~3.15]{Haa73} for a similar discussion of
$\Aut(L^\infty([0,1]))$).

We consider
$\cM = L^\infty([0,1])$, $\cH = L^2([0,1])$ and note
that $\cM_* \cong L^1([0,1])$.
Let $G := \Homeo([0,1])_\mu \subeq \Homeo([0,1])$
be the subgroup consisting of all homeomorphisms mapping Lebesgue
zero sets to Lebesgue zero sets, i.e.\ $g$ and $g^{-1}$ are absolutely continuous.
We topologize $G$ as a subgroup of $\Homeo([0,1])$ which carries the group
topology defined by the embedding
\[ \Homeo([0,1]) \to C([0,1])^2, \quad g \mapsto (g,g^{-1})\]
(\cite[Cor.~9.15]{Stp06}).
Then $G$ acts by automorphisms on the von Neumann algebra $\cM$
by $\alpha_g(f) := g_*f := f \circ g^{-1}$. We show that this action is continuous with respect to the
$p$-topology but not with respect to the $u$-topology (Remark~\ref{rem:2.1}(b)). This implies in particular
that on the group $\Aut(L^\infty([0,1]))$, these two topologies do not coincide.

{\bf Continuity in $p$-topology:} We consider the continuous bilinear map
\[ \beta \: L^\infty([0,1]) \times L^1([0,1]) \to \ell^\infty(G), \quad
\beta(f,h)(g) := \int_0^1 (g_* f)(x)h(x)\, dx.\]
We have to show that all functions $\beta(f,h)$ are continuous on $G$.
In view of $\beta(f,h)(g_1 g_2) = \beta((g_2)_*f,h)(g_1)$,
it suffices to verify continuity in $\be = \id_{[0,1]} \in G$.

Since $\beta$ is continuous and bilinear and the subspace $C(G)\cap \ell^\infty(G)$
is closed in $\ell^\infty(G)$, it suffices to do that for the case
where $h$ is bounded and $f = \chi_{[a,b]}$ is a characteristic function of an
interval $[a,b] \subeq [0,1]$. For
$\|g - \be\|_\infty < \eps$, we observe that
\[ E := g^{-1}([a,b]) \Delta [a,b] \subeq [a-\eps,a+\eps] \cup [b-\eps,b+\eps],\]
which leads to
\begin{align*}
|\beta(\chi_{[a,b]},h)(g)- \beta(\chi_{[a,b]},h)(\be)|
&= \Big|\int_0^1 (g_*\chi_{[a,b]}-\chi_{[a,b]})(x) h(x)\, dx\Big|
\leq \int_0^1 \chi_E(x) |h(x)|\, dx
\leq 4\eps \|h\|_\infty.
\end{align*}
This proves that the function $\beta(\chi_{[a,b]},h)$ is continuous at $\be$,
and hence that the homomorphism $\alpha \: G \to \Aut(\cM)$ is continuous with respect to the
$p$-topology.

{\bf Discontinuity in the $u$-topology:} Since the $u$-topology on
$\Aut(\cM)$ corresponds to the strong operator topology
for the action on $L^2([0,1])$
(see Example~\ref{ex:3.4}(a) and Remark~\ref{rem:aut-w*}),
we have to show that
the representation $U \: G \to \U(L^2([0,1]))$ defined by
$U_{g^{-1}} f := \sqrt{|g'|} \cdot (f \circ g)$ is not continuous. This will be achieved
by showing that the orbit map $G \to L^2([0,1]),
g \mapsto \sqrt{|g'|}$ for the constant function $1$
is not continuous at $\be$.

For every $n \in \N$, we consider the piecewise linear continuous function
$h_n \: [0,1] \to \R$, determined by its values at the joining points to be:
\[ h_n(x) :=\begin{cases}
0 & \mbox{ for } x = \frac{k}{2^n}, k = 0,\ldots, 2^n, \\
\big(1 - \frac{1}{2^n}\big) \frac{1}{2^{n+1}}  & \mbox{ for }
x = \frac{2k+1}{2^{n+1}}, k = 0,\ldots, 2^n-1.
\end{cases}
\]
Then
\[ g_n \: [0,1] \to [0,1], \qquad g_n(x) :=x + h_n(x) \]
defines a sequence in $G$. Note that these homeomorphisms are piecewise linear
with
\[ g_n'(x) :=\begin{cases}
2 - \frac{1}{2^n} & \mbox{ for } \frac{k}{2^n} < x < \frac{2k+1}{2^{n+1}},\\
\frac{1}{2^n} & \mbox{ for } \frac{2k+1}{2^{n+1}} < x < \frac{k+1}{2^n}.
\end{cases}
\]
As $g_n(x) = x$ for $x = \frac{k}{2^n}$, $k = 0,\ldots, 2^n$, and
$g_n$ is strictly increasing, we have
\[  \|g_n - \id\|_\infty \leq \frac{1}{2^n} \quad \mbox{ and } \quad
\|g_n^{-1} - \id\|_\infty \leq \frac{1}{2^n}.\]
This implies that $\lim_{n \to \infty} g_n = \be$ in $G$.
Next we observe that
\[ \| \sqrt{g_n'} - \sqrt{g_{n+1}'}\|^2_2
= 2\Big(1 - \int_0^1 \sqrt{g_n'} \sqrt{g_{n+1}'}\Big).\]
From
\begin{align*}
&\int_0^1 \sqrt{g_n'}\sqrt{g_{n+1}'} \\
&= \frac{1}{4}\sqrt{2 - \frac{1}{2^n}}\sqrt{2 - \frac{1}{2^{n+1}}}  +
\frac{1}{4}\sqrt{2 - \frac{1}{2^n}}\sqrt{\frac{1}{2^{n+1}}}+
\frac{1}{4}\sqrt{\frac{1}{2^n}}\sqrt{2 - \frac{1}{2^{n+1}}}+
\frac{1}{4}\sqrt{\frac{1}{2^n}}\sqrt{\frac{1}{2^{n+1}}}\to \frac{1}{2}
\end{align*}
it follows that
$\| \sqrt{g_n'} - \sqrt{g_{n+1}'}\|_2 \to 1$.
This shows that the sequence $U_{g_n^{-1}} 1 = \sqrt{g_n'}$ does not converge to $1$
in $L^2([0,1])$.
\end{example}

\begin{rem}\label{cont}
Given any action $(\cA, G, \alpha)$, we can always define the \strong continuous part of it by
\[
{\cal A}_c:=\{A\in{\cal A}\,\mid\, \alpha^A \: G \to \cA, g\mapsto\alpha_g(A)\quad\hbox{is norm continuous}\}
\quad \mbox{ and } \quad
\alpha^c_g:=\alpha_g\restriction{\cal A}_c.
\]
Unfortunately, as we will see in Example~\ref{Ac0} below, it is possible that $\cA_c = \C\1$.

If we start from  a $W^*$-dynamical system $(\cM,G, \beta)$ with $G$ locally compact,
 then ${\cal M}_c$ is weakly dense in $\cal M$, and
\[\cM_c=C^*\big\{\beta_f(A)\mid\, f\in L^1(G),\, A\in\cM\big\},\]
where the integrals $\beta_f(A) := \int_G f(g)\beta_g(A)\, dg$
exist in the weak topology  (\cite[Lemma~7.5.1]{Pe89}).

In the case that $\cM=\cA''$ for some concrete $C^*$-algebra $\cA$  invariant with respect to ~$G$,
 it is unfortunately
possible that $\cA\cap\cM_c=\C\1$. Moreover, in general only the representations of $\cM_c$
which are the restrictions of normal representations of $\cM$ will extend from $\cM_c$ to $\cM$ to produce
representations on $\cA$. Thus the $C^*$-dynamical system  $(\cM_c,G, \beta)$ is not a good vehicle to study
the general covariant representations of $(\cA,G, \beta)$.
\end{rem}

\begin{example}
\label{HeisAction}
Let $(\cdot\mid\cdot)$ be
the usual scalar product on $\R^n$.
The Heisenberg algebra is
$\fh_{2n+1}={\mathbb R}^n\times{\mathbb R}^n\times{\mathbb R}$
with the Lie bracket
$[(p,q,t),(p',q',t')]=[(0,0,(q\mid p')-(q'\mid p))]$,
and its corresponding simply connected Heisenberg group
is $\bH_{2n+1}=(\fh_{2n+1},\cdot)$
with
$x\cdot y=x+y+\frac{1}{2}[x,y]$.
The Schr\"o\-dinger representation is the irreducible unitary representation
$U\colon\bH_{2n+1}\to \U(L^2(\R^n))$
defined by
$(U(p,q,t)f)(x)=e^{i((q\mid x)+\frac{1}{2}(q\mid p)+t)}f(p+x)$
 for a.e.
 $x\in\R^n$,
for arbitrary $f\in L^2(\R^n)$ and $(p,q,t)\in\bH_{2n+1}$.
Let $S^0(\R^{2n})$ be the space of symbols of order zero, that is, the functions  $a\in C^\infty(\R^{2n})$
for which the partial derivatives of $a$ of arbitrary order are bounded.
Recall that the pseudo-differential Weyl calculus is a linear mapping
$\text{Op}\colon \cS'(\R^{2n})\to \cL(\cS(\R^n),\cS'(\R^n))$
satisfying
$$(\text{Op}(a)f)(x)=\int_{\R^{2n}}e^{i(x-y\mid z)}a((x+y)/2,z)f(z) dy\,dz $$
if $a\in\cS(\R^{2n})$ and $f\in\cS(\R^n)$,
where we denote by $\cS(\cdot)$ the Schwartz space, and by $\cS'(\cdot)$ its topological dual, that is, the space of tempered distributions.
By the classical Calder\'on-Vaillancourt theorem,
for every $a\in S^0(\R^{2n})$ the operator $\text{Op}(a)\colon \cS(\R^n)\to\cS'(\R^n)$ extends to a bounded linear operator $\text{Op}(a)\in \cB(L^2(\R^n))$.

Now denote $G=\bH_{2n+1}$, $\cH=L^2(\R^n)$, $\cM= \cB(\cH)$, and define $\beta\colon G\to\Aut(\cM)$ by $\beta_g(A)=U(g)AU(g)^{-1}$ for all $A\in\cM$ and $g\in G$.
We thus obtain the $W^*$-dynamical system $(\cM,G, \beta)$ and, with the notation of Remark~\ref{cont}, one can prove that $\cM_c$ is the norm-closure of the space of pseudo-differential operators of order zero $\text{Op}(S^0(\R^{2n}))\subseteq \cB(\cH)$.
(See \cite[Th. 1.1]{No12} and \cite[Ex. 1]{BB13}.)
Moreover, $\cK(\cH)$ is the norm-closure of $\text{Op}(\cS(\R^{2n}))$ and $\cK(\cH)\subsetneqq\cM_c$.
\end{example}

Typical examples of singular actions occur for bosonic systems:

\begin{example}
\label{ccrexmp1}
 Let $\cH$ be a nonzero complex Hilbert space and define a
symplectic form
$\sigma:\cH\times\cH\to\R$ by $\sigma(x,y):={\rm Im}{\langle x,y\rangle}$ where
${\langle \cdot,\cdot\rangle}$ denotes the inner product.
Then $(\cH,\sigma)$ is a symplectic space over $\R$,
and we let $\Sp(\cH,\sigma)$ denote the group of linear symplectic transformations of it.
Note that
unitaries on $\cH$ define symplectic transformations.
For the quantum system based on this we choose for its field algebra $\cA$  the Weyl algebra $\ccr \cH,\sigma.$
(cf.~\cite{Ma68}).
It is defined through the generators
$\{\delta_f\mid f\in \cH\}$ and the Weyl relations
\[ \delta_f^*=\delta_{-f} \quad \mbox{ and } \quad
\delta_f\delta_g=e^{-i\sigma(f,g)/2}\delta_{f+g}
\quad \mbox{ for } \quad f,g \in \cH.\]
Define a $C^*$-action $\alpha:\Sp(\cH,\sigma)\to\Aut(\cA)$ by $\alpha\s T.(\delta_x):=\delta_{T(x)}.$
Thus, if $U:G\to \U(\cH)$ is a nontrivial unitary representation of a connected topological group $G$ on $\cH,$ then
by $\U(\cH)\subset \Sp(\cH,\sigma),$ it defines a $C^*$-action on  $\ccr \cH,\sigma.$ by $\beta:=\alpha\circ U:G
\to\Aut(\cA)$. This action is not \strong continuous, hence singular. To see this, note that
${\|\delta_x-\delta_y\|}=2$ if $x\not=y$, so that the restriction of the $C^*$-topology to the subset
 $S:=\{ \delta_x \,\mid\, x \in \cH\}$ is discrete, the map $x\mapsto \delta_x$ is a bijection between $\cH$ and $S$
 and $\beta$ preserves $S$.
 If for an $x\in\cH$ the map $g\mapsto\beta_g(\delta_x)$ were continuous, then its image $\beta_G(\delta_x)=\delta_{U_Gx}$
 must be connected, and as the only connected open neighbourhood of $\delta_x$ in  $S$ is $\{\delta_x\}$ itself, it follows that
$\delta_{U_Gx}=\delta_x,$ i.e. $U_gx=x$ for all $g\in G$. If this is true for all $x,$ then
$U$ must be trivial, contrary to assumption. Thus there is an $x$ for which $g\mapsto\beta_g(\delta_x)$ is
not continuous, so $\beta$ is not \strong continuous, hence singular.
%
%
 \end{example}

 A second approach to bosonic quantum systems, is to replace the Weyl algebra above with the
 resolvent algebra (cf.~\cite{BG08}).
 \begin{example}
\label{ccrexmp2}
Let  $(\cH,\sigma)$ be the symplectic space above, and choose for its field algebra the resolvent
algebra $\rsl=:\cA$ (cf.~\cite{BG08}). It has a definition by a set of generators
\[ \{R(\lambda,x)\mid x\in \cH,\, \lambda\in\R^\times\} \]
satisfying a list of relations, but a much quicker way is to use the fact that in any
regular representation $\pi$ of  $\ccr \cH,\sigma.,$ we have up to *-isomorphism that
$\rsl$ is generated as a $C^*$-algebra by the resolvents of the fields $\varphi(x),$
where each $\varphi(x)$ is the self-adjoint
generator of the one-parameter unitary group $t\mapsto\pi(\delta_{tx}).$
Then $\pi$ can be defined on $\rsl$ by setting $\pi(R(\lambda,x))=\big(i\lambda\1-\varphi(x)\big)^{-1}.$
Similar to above, we define a $C^*$-action $\alpha:\Sp(\cH,\sigma)\to\Aut(\cA)$ by $\alpha\s T.(R(\lambda,x)):=R(\lambda,T(x)).$

Next, assume as above that
  we are given a non-trivial unitary representation $U:G\to \U(\cH)$
of a connected topological group $G$ on $\cH$.
Then, by $\U(\cH)\subset \Sp(\cH,\sigma),$ it defines a $C^*$-action on  $\rsl$ by $\beta:=\alpha\circ U:G\to\Aut(\cA)$.
Now we have by \cite[Thm.~5.3(ii)]{BG08} that
${\|R(1,x)-R(1,y)\|}=1$ if $x\not\in\R y,$ i.e.\ $x$ and $y$ are not on the same real ray.
Note that $\U(\cH)$ preserves the unit sphere $\bS\subset\cH,$ and each nontrivial unitary acts nontrivially on $\bS$,
and $\bS$ intersects each real ray in $\cH$ in exactly two points. Moreover for these points
we also have ${\|R(1,x)-R(1,-x)\|}=1$ by spectral theory.
Thus the $C^*$-topology restricts on the set $S:=\{R(1,x)\,\mid\,x\in\bS\}$
to the discrete topology, and the map $x\mapsto R(1,x)$ is a bijection between $\bS$ and $S$.
Moreover $\beta_G$ preserves the set $S$ and acts nontrivially on it.
 If for an $x\in\bS$ the map $g\mapsto\beta_g(R(1,x))$ were continuous, then its
image $\beta_G(R(1,x))=R(1,U_Gx)\subset S$
 must be connected, and as the only connected open neighbourhood of
 $R(1,x)$ in  $S$ is $\{R(1,x)\}$ itself, it follows that
$R(1,U_Gx)=R(1,x),$  i.e. $U_gx=x$ for all $g\in G$. By nontriviality of $U$ this cannot
be true for all $x\in \bS,$ hence there is some $x\in \bS$ for which  $g\mapsto\beta_g(R(1,x))$ is not continuous,
i.e.  $\beta$ is not \strong continuous, hence singular.
%
\end{example}
Note that if the connected topological group $G$ is not locally compact, then by definition even if a $C^*$-action
$G\to\Aut(\cA)$ is \strong continuous, then it is also singular.

\subsection{Covariant representations for a $C^*$-action.}
\mlabel{subsec:2.2}

\begin{defn} \mlabel{def:2.2}
(a) A {\it covariant representation} for a $C^*$-action  $(\cA, G, \alpha)$
is a pair $(\pi,U)$, where \break
$\pi \: \cA \to \cB(\cH)$ is a non-degenerate representation of $\cA$
on the Hilbert space $\cH$ and $U \: G \to \U(\cH)$ is a
\strongH continuous unitary representation satisfying
\begin{equation}
  \label{eq:covar}
U(g)\pi(A)U(g)^* = \pi(\alpha_g(A)) \quad
\mbox{ for } \quad g \in G, a \in \cA.
\end{equation}
For a fixed Hilbert space $\cH$, we write
${\rm Rep}(\alpha,\cH)$ for the set of covariant representations
$(\pi,U)$ of $(\cA, G, \alpha)$ on $\cH$.

(b) A non-degenerate representation $(\pi, \cH)$ of $\cA$ is called {\it covariant}
if there exists a \strongH continuous representation $U$ of $G$ such that
$(\pi, U)$ is a covariant representation of $(\cA,G,\alpha)$. It is called
{\it quasi-covariant} if $(\pi,\cH)$ is quasi-equivalent to a covariant representation
(cf.\ \cite{Bo69}). See Remark~\ref{rem:2.15}(c) below for more on quasi-equivalent
representations.

(c)
We write $\fS_\alpha(\cA)$ for the set of those states of
$\cA$ arising as vector states in covariant representations
of $(\cA,G, \alpha)$. By the Lemma~\ref{lem:2.3} below,
we have in fact $\fS_\alpha(\cA) \subeq (\cA^*)_c$.
A state $\omega \in \fS_\alpha(\cA)$ is called {\it covariant}
(resp.~{\it quasi-invariant}) if the corresponding cyclic
representation $\pi_\omega$ obtained by the GNS construction is covariant (resp.~quasi-covariant) (cf.\ \cite[Def.~6]{GK70};
see Theorem~\ref{thm:2.7} for more on quasi-invariant states).
Below we will characterize the covariant states.
\end{defn}

\begin{lem} \mlabel{lem:2.3} If $(\pi, U)$ is a covariant representation for $(\cA,G,\alpha)$
and $S \in \cB_1(\cH)$ a trace class operator, then
the continuous linear functional $\omega_S(A) :=\tr(\pi(A)S)$ on $\cA$ is contained in $(\cA^*)_c$.
\end{lem}

\begin{prf} (\cite[Prop.~3]{GK70})
   For $S \in \cB_1(\cH)$, we have
\begin{align*}
 (\omega_S \circ \alpha_g^{-1} - \omega_S)(A)
= \tr\big(U_g^* \pi(A) U_g S - \pi(A)S\big) = \tr\big( \pi(A) (U_g S U_g^*-S)\big),
\end{align*}
and since the conjugation action of $G$ on $\cB_1(\cH)$ is point-norm continuous,
$\omega_S \in (\cA^*)_c$.
\end{prf}

\begin{rem} (1) For a covariant representation  $(\pi,U)\in{\rm Rep}(\alpha,\cH)$ of the
$C^*$-action $(\cA,G, \alpha)$, the map $U \: G \to \U(\cH)$
is strong operator continuous by definition. Therefore
$\beta_g(A) := U_g A U_g^*$ defines a homomorphism
$\beta \: G \to \Aut(\cB(\cH))$ and Lemma~\ref{lem:2.3} shows that
$(\cB(\cH), G,\beta)$ is a $W^*$-dynamical system.
As the von Neumann algebra $\pi(\cA)'' \subeq \cB(\cH)$ is $\beta_G$-invariant, we also obtain a
 $W^*$-dynamical system $(\pi(\cA)'', G, \beta\res_{\pi(\cA)''})$
(cf.~\cite[7.4.2]{Pe89}). Conversely, given a  $W^*$-dynamical system $(\cM, G,\beta)$,
 it always has a faithful normal representation which is covariant
(cf.~equation ~\eqref{eq:standard-auto} below).
Note that the strong operator continuity in $g$ of $\pi\circ\alpha_g(A)={\rm Ad}(U_g)\circ\pi(A)$
on $\pi(\cA)$ does {\it not} make $g\to\pi\circ\alpha_g$ \strong continuous,
unless $\pi(\cA)\subseteq\cK(\cH),$ i.e. it consists of compact operators.
(Of course the continuous part of the $W^*$-dynamical system $(\cB(\cH), G,\beta)$ can be
larger than $\cK(\cH)$, cf. Example~\ref{HeisAction}).

(2) In  the \usual case,  for vector states, Lemma~\ref{lem:2.3}
 is \cite[Lemma~II.2]{Bo69}.

(3) Note that if  $(\pi, \cH)$ is covariant, then $\ker(\pi)$ is preserved by $\alpha$,
hence one can easily find non-covariant representations if $\cA$ has ideals which are not
preserved by $\alpha$.
\end{rem}

Let $(\cA,G, \alpha)$ be a $C^*$-action.
In the \usual case, where $G$ is locally compact and $\alpha$ is \strong continuous,
the covariant representations
are in bijective correspondence with the non-degenerate representations of the crossed
product $C^*$-algebra $\cA \rtimes_\alpha G$.
For a singular action, it is not obvious in general that covariant representations exist.
There  always exist covariant representations of $(\cA,G_d, \alpha)$, which
is an instance of the \usual case,
and if covariant representations of $(\cA,G, \alpha)$ exist, they will be amongst these.
Here is an example of a singular action with no covariant representations.

\begin{ex}
\mlabel{ExoticEx}
(A $C^*$-action $(\cA,G,\alpha)$ with no non-zero covariant representations)\\
A topological
group $G$ is called {\it exotic} if  all its (\strongH continuous) unitary representations
are trivial. In \cite[Ch.~2]{Ba91} one finds various
constructions of such a group of the type
$G = E/\Gamma$, where $E$ is a Banach space regarded as an additive group, and
$\Gamma \subeq E$ is a discrete subgroup.

Let $G$ be an exotic topological group.
Take the left regular representation on $\ell^2(G)$, i.e.\ $(V_g\psi)(h):= \psi(g^{-1}h)$ for $\psi\in \ell^2(G)$,
$g,\,h\in G$. Then $V$  is a non-trivial unitary representation of $G_d,$ but not continuous for $G.$
Let $\cA=\cK(\ell^2(G))$ which is a simple ideal of  $\cB(\ell^2(G))$
(but not the only closed ideal as $\ell^2(G)$ is nonseparable).
Define $\alpha:G\to \Aut(\cA)$ by $\alpha_g(A):=V_gAV_g^*$ which is a non-trivial action.
If $(\pi,U)$
were a covariant representation, then $U$ must be trivial as $G$ is exotic, hence
$\pi\circ\alpha_g=\pi$ for all $g\in G$.
However $\cA$ is simple and $\pi$ is non-trivial hence
$\pi$ is injective, and then $\alpha_g(A)=A$ for all $g\in G$ and $A\in\cA$, which is a contradiction.
Thus there are no
non-trivial covariant representations, i.e.\  $\fS_\alpha(\cA)=\emptyset.$
\end{ex}

In the  subsections below, we will consider the problem of the existence of covariant representations
in some detail. In Corollary~\ref{corCovRepXst} we will obtain conditions
characterizing the existence of a covariant representation. Regarding explicit constructions, it is well-known that one can obtain a covariant
representation for singular actions either from
\begin{itemize}
\item{}
standard form representations of $W^*$-dynamical systems,
\item{}
 from the representations of $W^*$-crossed products of  $W^*$-dynamical systems, or
\item{} from invariant states with appropriate continuity conditions (cf.~\cite{DJP03}).
\end{itemize}
These will be considered below in Section~\ref{SectStandardForm}
     and  Proposition~\ref{fixstatecov} respectively.
     Below in Theorem~\ref{thm:Pinv} we will obtain a covariant representation via
     the Stinespring Dilation Theorem.

There are also natural uniqueness and structure questions, e.g.\ given a covariant representation
 $(\pi,\cH)$ for a $C^*$-action $(\cA,G, \alpha)$, find and analyze all unitary representations
 $U:G\to\cU(\cH)$ for which   $(\pi,U)\in{\rm Rep}(\alpha,\cH)$ is a covariant representation.
 Below we will see that if a spectral condition is added, then we can find a natural ``minimal''
such $U:\R\to\cU(\cH)$ which is unique.

\subsection{Folia in $\cA^*$ and the Borchers--Halpern Theorem}
\mlabel{subsec:2.3}

In this subsection we will characterize when a representation $\pi$
is covariant in terms of properties of its set of normal states,
i.e., the corresponding folium $F(\pi)$.

\begin{defn}
\label{FolDef}
(a) For a $C^*$-algebra $\cA$, we call a subset $F \subeq \cA^*$ a {\it folium} if there exists a
representation $(\pi, \cH)$ of $\cA$ with
\begin{equation}
 \label{eq:folpi}
F = F(\pi) := \{ \omega_S \in\fS(\cA)\mid 0 \leq S \in \cB_1(\cH),\ \tr S = 1\}
\end{equation}
where $\omega_S(A):= \tr(\pi(A)S)$   as in Lemma~\ref{lem:2.3}.

(b) We likewise define the folium $F(\pi) \subeq \cM_*$ of a normal representation $(\pi,\cH)$
of a $W^*$-algebra~$\cM$.
\end{defn}

As  the normal states of $\cB(\cH)$ are identified with trace class operators
by $\omega_S(A) = \tr(SA)$, we have
\begin{equation}
  \label{eq:4}
F(\pi) = \{ \omega\circ\pi \mid \omega\hbox{ is a normal state of }\pi(\cA)''\} \cong \fS_n(\pi(\cA)'')
\end{equation}
because all normal states of $\pi(\cA)''$ extend to normal states of $\cB(\cH)$
(cf.~\cite[Cor.~III.2.1.10]{Bla06}).
Clearly $F(\pi)$ inherits from $\pi(\cA)''_*$ the convexity and invariance under conjugations.
We verify that it also inherits norm closedness.

\begin{lem} \mlabel{lem:2.2} Let $(\pi, \cH)$ be a representation of $\cA$.
\begin{itemize}
\item[\rm(a)] The folium
$F(\pi)\subset\cA^*$ and its linear span are both norm
closed.
\item[\rm(b)] Moreover, $ F(\pi)$  coincides with the set of vector states of the representation
$(\rho, \cB_2(\cH))$ of $\cA$ given by $\rho(A)B := \pi(A)B$, $A\in\cA$, $B\in\cB_2(\cH)$.
\end{itemize}
\end{lem}

\begin{prf} (a) The restriction map
$(\pi(\cA)'')_* \to \pi(\cA)^*$ is isometric (\cite[Prop.~2.12]{BN12}),
and the subset $\fS_n(\pi(\cA)'') \subeq (\pi(\cA)'')_*$ of normal states is norm closed.
This implies that
\[ \Spann F(\pi) = \{ \omega_S | S \in \cB_1(\cH)\} \]  is norm-closed in $\cA^*$,
and this shows the norm closedness of $F(\pi)$.

(b) The vector states of $\rho$ are of the form
\[ \omega_B(A) = \la B, \pi(A)B\ra = \tr(B^* \pi(A)B) = \tr(BB^*\pi(A)),\]
 where $BB^*$ is a positive trace class operators with $\tr(B^*B) = \|B\|_2^2 = 1$.
Hence these are precisely the states of the form
$\omega_S$, $0 \leq S \in \cB_1(\cH)$ with $\tr S = 1$.
Therefore $F(\pi)$ coincides with the set of vector states of $\rho$.
\end{prf}

\begin{rem} \mlabel{rem:2.15}
(a) In \cite{Ka62} it is shown that the set of vector states $V(\pi)
\subeq F(\pi) \subeq \cA^*$ of a representation $(\pi, \cH)$
of a $C^*$-algebra $\cA$ is a norm closed
subset. This implies the closedness of $F(\pi)$ since $F(\pi) = V(\rho)$ by
Lemma~\ref{lem:2.2}(b).
However, the closedness of the larger set $F(\pi)$
is much easier to get.

(b) A folium $F \subeq \fS(\cA)$ can be abstractly characterized
as a convex set of states which is norm closed, and contains with any state $\omega\in F,$ all states of the form
 \begin{equation}
   \label{eq:twodideact}
 (B * \omega)(A) := \frac{\omega(B^*AB)}{\omega(B^*B)}, \quad \omega(B^*B) > 0
 \end{equation}
(\cite[p.~84]{HKK70}). This is a better intrinsic definition of a folium
as it does not rely on the existence of a representation $\pi$.

(c) For a state $\omega\in\fS(\cA)$,
the folium $F(\pi_\omega)$ is the norm-closed convex subset generated by the set
$\{B * \omega | \omega(B^*B) > 0\}$ (cf.~\eqref{eq:twodideact}).
By polarization, $F(\pi_\omega)$ generates the same norm closed subspace of
$\cA^*$ as $\cA\omega\cA$, where we define
\begin{equation}
  \label{eq:leftrightshift}
(A\omega)(B) := \omega(AB) \quad \mbox{ and } \quad
(\omega A)(B) := \omega(BA) \quad \mbox{ for } \quad A,B \in \cA.
\end{equation}
As $\Spann F(\pi_\omega)$ is norm closed by
Lemma~\ref{lem:2.2},
 we see that
\begin{equation}
  \label{eq:fol-rel}
\Spann F(\pi_\omega) = \lbr\cA \omega \cA \rbr.
\end{equation}
where $\lbr\cdot\rbr$ denotes the closed span of its argument.

(d) For two representations
$\pi_1$ and $\pi_2$ of  $\cA$,
their folia are equal $F(\pi_1) = F(\pi_2)$ if and only if
they are quasi-equivalent, i.e.\ there is an isomorphism $\beta:\pi_1(\cA)''\to\pi_2(\cA)''$ of
$W^*$-algebras such that $\beta(\pi_1(A))=\pi_2(A)$ for all $A\in\cA$
(cf.~\cite[Prop.~10.3.13]{KR86}). This
means that each subrepresentation of $\pi_1$ has a subrepresentation
which is unitarily equivalent to a subrepresentation of $\pi_2$, and vice versa
(\cite[Cor.~10.3.4(ii)]{KR86}). This statement is also contained in
 \cite[Cor.~5.11]{AS01} as the corresponding ``split faces''  are the corresponding
folia in our terminology.

(e) A subset
$E \subeq \fS(\cA)$ is contained in the folium $F(\pi)$ of a representation
$(\pi,\cH)$ if and only if the cyclic
representations $(\pi_\omega, \cH_\omega)$, $\omega \in E$, are contained in the
corresponding left multiplication representation $(\rho,\cB_2(\cH))$ with
$\rho(A)(B) = \pi(A)B$, which satisfies $F(\rho) = F(\pi)$.
Therefore every subset $E \subeq \fS(\cA)$ is contained in a {\it minimal
folium} $\Fol(E)$ which can be obtained as $F\big(\bigoplus_{\omega \in E} \pi_\omega\big)$.
This further implies that
\[ \Fol(E) = \Big\{ \sum_{n = 1}^\infty c_n \nu_n \;\Big\vert\; 0 \leq c_n \leq 1,
\sum_{n=1}^\infty c_n = 1,
\nu_n \in F(\pi_{\omega_n}), \omega_n \in E\Big\}.\]
\end{rem}

\begin{ex}
For $\cA = C_0(X)$, $X$ locally compact, and a
state $\omega \in \fS(\cA)$ obtained from a probability measure by
$\omega(A)=\int A\,d\mu$, the corresponding folium can be determined rather easily
from (\ref{eq:twodideact}).
For $f \in \cA$ with $\int_X |f|^2\, d\mu = 1$, we have
$f * \omega = |f|^2 \omega$.
Since the embedding $L^1(X,\mu)\into \cA^*, h \mapsto h \cdot \omega$ is isometric,
it follows  that
\[ \Fol(\omega) = \Big\{ F \omega \,\mid\, F \in L^1(X,\mu), 0 \leq F,
\int_X F\, d\mu = 1\Big\}\]
corresponds to
 the set of probability measures which are absolutely continuous with respect to~$\mu$.
\end{ex}

\begin{thm} {\rm(Borchers--Halpern Theorem)} \mlabel{thm:bor}
Let $(\cA,G,\alpha)$ be a $C^*$-action and $F \subeq \fS(\cA)$ be a folium.
Then there exists a covariant representation $(\pi,U,\cH)$ of $(\cA,G,\alpha)$
with $F = F(\pi)$ if and only if $F$ is $\alpha^*_G$-invariant and contained in $(\cA^*)_c$.
\end{thm}

\begin{prf} (cf.~\cite[p.~258]{Hal72}, \cite[Thm.~III.2]{Bo83})
Below we will also obtain a proof of this from standard forms in
Remark~\ref{SFPf}.
\end{prf}

Kadison's old paper \cite{Ka65} already contains an
interesting precursor of this theorem.

\begin{cor}\mlabel{cor:bor1}
Let $(\cA,G,\alpha)$ be a $C^*$-action and $(\pi, \cH)$ be a non-degenerate
representation of $\cA$. Then the following are equivalent:
\begin{itemize}
\item[\rm(i)] $\pi$ is quasi-covariant.
\item[\rm(ii)] $F(\pi)$ is $\alpha_G^*$-invariant and contained in $(\cA^*)_c$.
\item[\rm(iii)] We have that $\ker \pi$ is $\alpha_G$-invariant, hence
the induced action of $G$ on $\pi(\cA)$ is defined.
Moreover, this induced action of $G$ on $\pi(\cA)$ extends to an
action  $\beta:G\to\Aut(\pi(\cA)'')$, defining a $W^*$-dynamical system.
\end{itemize}
\end{cor}

\begin{prf} The equivalence between (i) and (ii) follows by applying the
Borchers--Halpern Theorem to the folium $F = F(\pi)$.

Next we show the equivalence of (ii) and (iii).
Note that if $\pi$ is quasi-covariant, then its kernel must coincide
with the kernel of a covariant representation, and this is always invariant
with respect to ~$\alpha_G$. Thus the induced action of $G$ on $\pi(\cA)$
is defined.
As  quasi-covariance of $\pi$ implies (iii), we only need to prove the converse.
That $(\pi(\cA)'',G,\beta)$
is a $W^*$-dynamical system implies that $(\pi(\cA)'')_*$ is
 $\beta_G^*$-invariant, hence  $\alpha_G^*$-invariant as a $\varphi\in(\pi(\cA)'')_*$
 is uniquely specified by its restriction to $\cA$. Thus
 $F(\pi)=(\pi(\cA)'')_*\cap\fS(\cA)$ is $\alpha_G^*$-invariant.
 By definition of a $W^*$-dynamical system, $G$ acts continuously on
$(\pi(\cA)'')_*$, hence $F(\pi) \subeq (\cA^*)_c$.
Thus we have obtained equivalence with (ii).
\end{prf}

\begin{rem}(a) The existence of non-zero covariant representations is equivalent to the existence
of non-zero quasi-covariant representations of $\cA$. Thus
Corollary~\ref{cor:bor1}(iii) is a criterion for the existence of covariant representations.

(b) The question of when a quasi-covariant representation is actually covariant, was
analyzed by Bulinskii in \cite{Bu73a, Bu73b}, but below in Subsection~\ref{sec:3} we
will see better conditions.
\end{rem}
Corollary~\ref{cor:bor1} has a specialization which can answer the following question.
Given a $C^*$-algebra $\cA\subset\cB(\cH)$ and an automorphism
$\gamma\in \Aut\cA$, when does  $\gamma$ extend to an automorphism of $\cA''$?
\begin{cor}\mlabel{cor:bor2}
Let $\cA$ be a $C^*$-algebra, let $(\pi, \cH)$ be a non-degenerate
representation of $\cA$ and let $\gamma\in \Aut\cA$ be an automorphism
such that $\ker \pi$ is $\gamma$-invariant.
  Then the following are equivalent:
\begin{itemize}
\item[\rm(i)] $F(\pi)$ is $\gamma$-invariant.
\item[\rm(ii)] The induced automorphism of $\gamma$ on $\pi(\cA)$ extends to an
automorphism on  $\pi(\cA)''$.
\end{itemize}
Moreover, if this is the case, and if $\pi$ is irreducible, then
$\gamma$ is unitarily implementable on $\pi(\cA)$.
\end{cor}

\begin{prf}
Let $G\subset\Aut\cA$ be the discrete group generated by $\gamma$. This defines a
$C^*$-action $(\cA,G,\alpha)$ for which   $F(\pi)\subset(\cA^*)_c$. Moreover, if
$\gamma$  extends to an
automorphism on  $\pi(\cA)''$, it automatically defines a  $W^*$-dynamical system
with respect to~$G$.
Thus by Corollary~\ref{cor:bor1}, it follows that (i) and (ii) are equivalent.

If $\pi$ is irreducible, then $\pi(\cA)''=\cB(\cH)$, so as all automorphisms of $\cB(\cH)$ are inner,
the last statement follows.
\end{prf}

We now consider covariance conditions for states.

\begin{thm} \mlabel{thm:2.7}
For  a $C^*$-action $(\cA,G,\alpha)$ and a state $\omega$ of $\cA$, the following are equivalent:
\begin{itemize}
\item[\rm(i)] $\omega \in \fS_\alpha(\cA)$, i.e.\ $\omega$ is a vector state of some
covariant representation $(\pi, U)$ of $(\cA,G, \alpha)$.
\item[\rm(ii)] $\cA \omega \cA \subeq (\cA^*)_c$.
\item[\rm(iii)] $\Fol(\omega) := F(\pi_\omega) \subeq (\cA^*)_c$.
\item[\rm(iv)] $\Fol_G(\omega) := \Fol(\alpha_G^*\omega)\subeq (\cA^*)_c$.
\end{itemize}
Furthermore, the following are equivalent for $\omega \in \fS_\alpha(\cA)$:
\begin{itemize}
\item[\rm(a)] $\omega$ is quasi-invariant (cf. Def.~\ref{def:2.2}(c)).
\item[\rm(b)]  $\pi_\omega$ is quasi-covariant.
\item[\rm(c)] $F(\pi_\omega) = \Fol(\omega)$ is $\alpha_G^*$-invariant.
\item[\rm(d)] $\pi_\omega$ is equivalent to a subrepresentation of a covariant
representation $\pi$ with $F(\pi) = F(\pi_\omega)$.
\end{itemize}
\end{thm}

\begin{prf}
Observe first, that for a subset $E \subeq \fS(\cA)$, the folium
$\Fol(E)$ generated by $E$ is equal to the norm closed convex hull
of the union of the folia $\Fol(\nu)=F(\pi_\nu) $, $\nu \in E$, and the span of each of these is equal to $\lbr \cA \nu \cA\rbr$ by (\ref{eq:fol-rel}). Hence $\cA E \cA \subeq (\cA^*)_c$
is equivalent to $\Fol(E) \subeq (\cA^*)_c$ as $(\cA^*)_c$ is a norm-closed subspace of
$\cA^*$.

(ii)$\Leftrightarrow$ (iii) follows directly from \eqref{eq:fol-rel} and the norm
closedness of $(\cA^*)_c$.

(i) $\Rarrow$ (ii): If $(\pi, U)$ is a covariant
representation with vector state $\omega$, then $\omega \in F(\pi)$. This implies that
$\cA \omega \cA \subeq \Spann F(\pi) \subeq (\cA^*)_c.$

(ii) $\Rarrow$ (i): Since $(\cA^*)_c$ is $G$-invariant, condition (ii) implies that
\[ A \alpha_g^*(\omega) B = \alpha_g^*\big(\alpha_g^{-1}(A) \cdot\omega
\cdot \alpha_g^{-1}(B)\big) \in(\cA^*)_c ,\]
and hence that $\cA (\alpha_G^*\omega) \cA \subeq (\cA^*)_c$.
We thus obtain by the first part of the proof that
$\Fol(\alpha_G^*\omega) \subeq (\cA^*)_c$.
The $G$-invariance of $\Fol(\alpha_G^*\omega)$ follows from the fact that it is
generated by a $G$-invariant subset of $(\cA^*)_c$. Therefore the
Borchers--Halpern Theorem~\ref{thm:bor} implies the existence of a covariant representation
$(\pi,U)$ of $(\cA,G,\alpha)$ with $F(\pi) = \Fol(\alpha_G^*\omega) \ni \omega$.

(i)$\Leftrightarrow$ (iv):
The $G$-orbit $\alpha_G^*\omega =\{ \omega \circ \alpha_g \mid  g \in G\}$
generates a folium $\Fol_G(\omega) = \Fol(\alpha_G^*\omega)$, and since
\[ \alpha_g^* \Fol(E) = \Fol(\alpha_g^*E) \quad \mbox{ for } \quad E \subeq \fS(\cA), g \in G, \]
the folium $\Fol_G(\omega)$ is $G$-invariant.
It is the minimal $G$-invariant folium containing
$\omega$. Hence the equivalence of (i) and (iv) follows from Theorem~\ref{thm:bor}.

Now we prove the second part of the theorem.
Assertions (a) and (b) are equivalent by definition. That (b) is equivalent to (c) follows
from Theorem~\ref{thm:bor} and Corollary~\ref{cor:bor1} since $\omega \in \fS_\alpha(\cA)$ implies
$F(\pi_\omega) \subeq (\cA^*)_c$ by (iii).

If (d) holds, i.e.\ $F(\pi) = F(\pi_\omega)$ for a covariant representation $\pi$, then (c) follows from
Corollary~\ref{cor:bor1}. Suppose, conversely, that (c) holds.
Then $F(\tau) = F(\pi_\omega)$ for a covariant representation $(\tau,U,\cH)$ by the
Borchers--Halpern Theorem. Consider the representation $(\pi, V, \cB_2(\cH))$
with $\pi(A)B := \tau(A)B$ and $V_g B := U_g B$. This covariant
representation satisfies $F(\pi) =F(\tau)$, but
$\omega \in F(\pi_\omega) = F(\tau) = F(\pi)$ is a vector state of
$\pi$ by Lemma~\ref{lem:2.2}. If  $B \in \cB_2(\cH)$ is such that
$\omega(A) = \tr(B^* \tau(A) B)$ for $A \in \cA$, then
the cyclic subrepresentation of $\cA$ on $\lbr \pi(A)B\rbr$ is
equivalent to~$\pi_\omega$. Then (d) follows.
\end{prf}

\begin{rem}
(a) Theorem~\ref{thm:2.7} improves \cite[Thm.~III.2]{Bo83}, in
that we already obtain the existence of a covariant representation
from the condition $\Fol(\omega) = F(\pi_\omega) \subeq (\cA^*)_c$, the
$\alpha_G^*$-invariance of $\Fol(\omega)$ is not required.

(b)  Note that if $\alpha$ is uniformly continuous (but $G$ need not be locally compact)
then $(\cA^*)_c=\cA^*$, so the properties (ii) and (iii) are trivially satisfied, hence by (i)
we have $\fS_\alpha(\cA)=\fS(\cA)$. So covariant representations always exist for this case.
\end{rem}

For the \usual case, the following corollary can already be found in \cite[Thm.~III.1]{Bo69}.
\begin{cor} \mlabel{cor:2.21} Let $(\cA,G,\alpha)$ be a $C^*$-dynamical system (hence $\alpha$
is point-norm continuous) and
$\omega \in \fS(\cA)$. Then $\omega \in \fS_\alpha(\cA)$ if and only if
$\omega \in (\cA^*)_c$, i.e.,
\[ \fS_\alpha(\cA) = \fS(\cA) \cap (\cA^*)_c.\]
\end{cor}

\begin{prf} In view of Theorem~\ref{thm:2.7},
we have to show that $\omega \in (\cA^*)_c$ implies
$\cA \omega \cA \subeq (\cA^*)_c.$
The trilinear map
\[ \cA \times \cA \times \cA^* \mapsto
(A,B,\omega) \mapsto A\omega B \]
is continuous because $\|A\omega B\| \leq \|A\|\|\omega\| \|B\|$
(cf.~\eqref{eq:leftrightshift}).
This map is $G$-equivariant, and this implies that
$\cA (\cA^*)_c \cA \subeq (\cA^*)_c$, using the \strong
continuity of $g\mapsto\alpha_g$.
\end{prf}
Thus covariant representations always exist if $\alpha$ is \strong
continuous and $\fS(\cA) \cap (\cA^*)_c\not=\emptyset$.

Theorem~\ref{thm:2.7} has the following corollary:
\begin{cor}\mlabel{corCovRepXst}
Given a $C^*$-action $(\cA,G,\alpha)$, the following are equivalent:
\begin{itemize}
\item[\rm(i)] There exists a non-zero covariant representation
$(\pi, U)$ of $(\cA,G, \alpha)$.
\item[\rm(ii)] There is a state $\omega\in(\cA^*)_c$ such that
$\cA \omega \cA \subeq (\cA^*)_c$.
\item[\rm(iii)] There exists a $\alpha_G^*$-invariant folium $F \subeq (\cA^*)_c$.
\item[\rm(iv)] There is a state $\omega\in(\cA^*)_c$ such that
$B * \omega \in (\cA^*)_c$ whenever $\omega(B^*B) > 0$.
\end{itemize}
\end{cor}

In the \usual case the GNS representations of
invariant states always produce covariant representations.
In the next example we see that for singular actions invariant states need not even be
in  $ \fS_\alpha(\cA)$. By Corollary~\ref{cor:2.21}, this requires $\alpha$
to be discontinuous in the \strong topology.

\begin{ex}
\mlabel{AcNotFol}
We construct an example of an invariant state
$\omega\not\in \fS_\alpha(\cA)$. Then $\Fol(\omega)$ is
$\alpha_G^*$-invariant but not contained in $(\cA^*)_c$, which implies that
$\fS(\cA)_c$ is not a folium.

Let $(X,\sigma)$ be the non-degenerate symplectic space over $\R$  given by $X=\C$,
$\sigma(z,w):=\Im(z\overline{w})$
and  $\cA=\ccr X,\sigma.$ is the associated Weyl $C^*$-algebra
with generating unitaries  $(\delta_z)_{z \in X}$ satisfying
\[ \delta_z^*=\delta_{-z} \quad \mbox{ and } \quad
\delta_z\delta_w=e^{-i\sigma(z,w)/2}\delta_{z+w}
\quad \mbox{ for } \quad z,w \in X.\]
The tracial state $\omega_0$ defined by $\omega_0(\delta_z)=\delta_{z,0}$
is invariant with respect to  the action of $G = \R$ on $\cA$ by
$\alpha_\theta (\delta_z) = \delta_{e^{i\theta}z}$.
For this action we clearly have  that $\omega_0\in(\cA^*)_c$ by its invariance.
Now $\delta_z\omega_0\in\cA \omega_0 \cA$
for $z\not=0$. Thus
\[
\alpha^*_\theta(\delta_z\omega_0)(\delta_{-z}) =
\omega_0\big(\delta_z\alpha_\theta(\delta_{-z})\big)=
\omega_0\big(\delta_z\delta_{-e^{i\theta}z}\big)
=e^{i\sigma(z,e^{i\theta}z)/2}\omega_0(\delta\s{(1-e^{i\theta})z}.)
\]
and this expression is nonzero only when $e^{i\theta}=1$ (when it has modulus $1$) hence it is discontinuous
with respect to $\theta$. Thus $\delta_z\omega_0\not\in(\cA^*)_c$ and Theorem~\ref{thm:2.7} implies that
$\omega_0\not\in\fS_\alpha(\cA)$.
\end{ex}

\begin{rem}
\mlabel{remFolia}
(a) Example~\ref{AcNotFol} shows that the inclusion $\fS_\alpha(\cA) \into
\fS(\cA)_c := \fS(\cA) \cap (\cA^*)_c$ may be proper.

(b) If $\fS(\cA)_c$ is a folium, by its  $G$-invariance, the Borchers-Halpern Theorem
implies that $\fS_\alpha(\cA) = \fS(\cA)_c$, which is not always the case by (a).
Therefore $\fS(\cA)_c$ is not always a folium.

(c) A similar situation arises
for a $W^*$-dynamical system $(\cM,G,\alpha)$ because the
{weak-$*$} continuity of orbit maps in $\cM_*$ does not imply that
$\cM_* = (\cM_*)_c$ (cf.~Example~\ref{ex:linft}).
Accordingly, Theorem~\ref{thm:bor} implies that the vector states of normal covariant
representations can be characterized by
\[ \fS_{n,\alpha}(\cM) = \{ \omega \in \fS_n(\cM) \mid \Fol_G(\omega) \subeq (\cM_*)_c\}.\]

(d) Suppose now that
$(\pi, \cM)$ is a normal representation of $\cM$ whose folium
$F(\pi)\cong \fS_n(\pi(\cM))$ is $G$-invariant and contained in $(\cM_*)_c$. Then
by Corollary~\ref{cor:bor1},
$\ker \pi$ is $G$-invariant, so that we obtain a natural $G$-action
$\alpha^\cN$ on $\cN := \pi(\cM)$ for which the action on $\pi(\cM)_*$ is continuous.
We thus obtain a $W^*$-dynamical system $(\cN,G,\alpha^\cN)$.
By the Borchers-Halpern Theorem we know that this has a covariant representation with folium
$F(\pi)$. Below, we will obtain such a  representation from
the standard form realization of $\cN$.
\end{rem}

As invariant states are important to construct covariant representations
(cf.~ground states, Definition~\ref{defgroundst0} below, as well as KMS states),
we need to characterize when they do produce
covariant representations for singular actions.
 Observe first that given a $C^*$-action $(\cA,G,\alpha)$ and an invariant state $\omega
\in\fS(\cA)$, then its GNS representation ${(\pi_\omega,\Omega_\omega,\cH_\omega)}$ always gives a
covariant representation  $(\pi_\omega,U^\omega)$ of $(\cA,G_d, \alpha)$,
where $U^\omega_g$ is uniquely determined by
\[
U^\omega_g\Omega_\omega=\Omega_\omega\quad\hbox{and}\quad
U^\omega_g\pi_\omega(A)\Omega_\omega:=\pi_\omega(\alpha_g(A))\Omega_\omega
\quad \mbox{ for all } \quad g\in G, A\in\cA \]
(cf.~\cite[Lemma~IV.4.4]{Bo96}). We then have:

\begin{prop} \mlabel{fixstatecov}
For a $C^*$-action $(\cA,G,\alpha)$ and an invariant state $\omega
\in\fS(\cA)$, the following are equivalent:
\begin{itemize}
\item[\rm(i)] $U^\omega:G\to\U(\cH_\omega)$ is continuous, i.e.\
$(\pi_\omega,U^\omega)\in{\rm Rep}(\alpha,\cH)$.
\item[\rm(ii)]  $\pi_\omega$ is covariant.
\item[\rm(iii)] $\cA \omega \cA \subeq (\cA^*)_c$.
\item[\rm(iv)] $ \cA\omega  \subeq (\cA^*)_c$.
\end{itemize}
If $(\cA,G,\alpha)$ is a $C^*$-dynamical system, then {\rm(i)}-{\rm(iv)} are satisfied.
\end{prop}

\begin{prf} (i) $\Rarrow$ (ii) is trivial.

(ii) $\Rarrow$ (iii): If  $(\pi_\omega,U)\in{\rm Rep}(\alpha,\cH_\omega)$ then $\omega \in \fS_\alpha(\cA)$,
hence by Theorem~\ref{thm:2.7} we obtain that $\cA \omega \cA \subeq (\cA^*)_c$.

(iii) $\Rarrow$ (iv): Assume that $\cA \omega \cA \subeq (\cA^*)_c$.
Observe that as $\pi_\omega(\cA)\Omega_\omega$ is dense in $\cH_\omega$, this condition just states
that the bounded
maps $g\mapsto\pi_\omega(\alpha_g(A))$ are continuous in the weak operator topology for  all
$A\in\cA$. As $\Omega_\omega\in\cH_\omega$, this implies that  $\cA \omega \cup \omega\cA\cup\{\omega\} \subeq (\cA^*)_c$.

(iv) $\Rarrow$ (i):
Condition (iv) implies that for $A, B \in \cA$, the function
\[ g \mapsto  \big(\pi_\omega(A)\Omega_\omega,U_g^\omega\pi_\omega(B)\Omega_\omega\big)
=\omega(A^*\alpha_g(B))
=(A^*\omega)(\alpha_g(B))\]
is continuous, hence that $g\mapsto U_g^\omega$ is weak operator continuous,
by density of $\pi_\omega(\cA)\Omega_\omega$. As the weak operator topology coincides with the strong operator topology
on the unitary group, we conclude that $U^\omega$ is continuous.

Finally, we assume that the $G$-action on $\cA$ is continuous, i.e.\  that
$(\cA,G,\alpha)$ is a $C^*$-dynamical system. Then
(i) follows from the fact that the subspace
$\cH_{\omega,c}$ of $U^\omega$-continuous vectors in $\cH_\omega$ is
$\pi_\omega(\cA)$-invariant and contains $\Omega_\omega$. Hence it coincides with
$\cH_\omega$.
\end{prf}

If $\omega \in \fS_\alpha(\cA)$ is not $G$-invariant,
then the remaining condition $\cA \omega \cA \subeq (\cA^*)_c$
is not enough to ensure that $\pi_\omega$ is covariant,
as Example~\ref{ex:3.13} shows. It does imply that
 $\pi_\omega$ is equivalent to a {\it sub}representation of a covariant representation by
Theorem~\ref{thm:2.7}. Below in Theorem~\ref{thm:Pinv} we will obtain a generalization
of this theorem to invariant projections, where the GNS representation $\pi_\omega$ has to be
replaced with a Stinespring dilation.

 \begin{ex}\mlabel{ex:3.13} {\rm(A non-quasi-covariant representation)}
Let $\alpha:\R\to\Aut C_0(\R)$ be the action of  translation on $\cA=C_0(\R)$. Consider the covariant representation
$(\pi,U)$, where
$C_0(\R)$ acts by multiplication on $\cH=L^2(\R)$ and the implementing unitaries $U_t$ act by right translation on $L^2(\R)$.
Let $\xi=\chi_{[0,1]}\in L^2(\R)$ and let $\omega_\xi$ be the associated vector state. Then
for the positive vector functionals $\omega_{A\xi}=|A|^2\omega_\xi,$ $A\in\cA$ we have
${\lim\limits_{t\to 0}\|\alpha_t^* \omega_{A\xi}-\omega_{A\xi}\|}=0$
(cf.~\cite[Lemma~II.2]{Bo69}). By polarization we thus
get $\cA\omega_\xi\cA\subset(\cA^*)_c$. Now $\pi_{\omega_\xi}$ is unitarily equivalent
to the restriction of $\pi(\cA)$ to $L^2[0,1]\subset\cH$. As the kernel of $\pi_{\omega_\xi}$ is $\{f\in C_0(\R)\mid\,
f\restriction[0,1]=0\}$ which is not translation invariant, the representation
$\pi_{\omega_\xi}$ is not quasi-covariant.
However, by construction there exists a covariant representation $(\pi, U)$ such that
$\pi_{\omega_\xi}$ is a subrepresentation of~$\pi$.
\end{ex}

\subsection{Covariance of cyclic representations}
\mlabel{cyclic-cov}

In Theorem~\ref{thm:2.7} we saw that a state $\omega$ is a vector state
of a covariant representation of $(\cA,G,\alpha)$ if and only if
 $\cA \omega \cA \subeq (\cA^*)_c$.  In the case that
$\omega$ is invariant, by Proposition~\ref{fixstatecov},
this condition is even enough to ensure that $\pi_\omega$
is covariant. This raises the question of how one can characterize
for the general case when  a GNS representation $\pi_\omega$ is covariant.
First, following the path of the Wigner Theorem, we have to characterize
for a single automorphism, whether it is implementable in $\pi_\omega$.

\begin{defn}
  \mlabel{def:carrier}
Let $\cM$ be a $W^*$-algebra. For $\omega \in \fS_n(\cM)$
we write $s(\omega) \in \cM$ for the corresponding
{\it carrier projection}, also called the {\it support of $\omega$}. It is the maximal
projection $p \in \cM$ with $\omega(p) = 1$ and
\begin{equation}
\label{carrierW}
  \{ M \in \cM | \omega(M^*M) = 0\} = \cM (\1 - p)\end{equation}
(cf.\ \cite[Def.~2.133]{AS01} or \cite[8.15.4]{Pe89}).
The  {\it central support of $\omega$}, denoted $z(\omega)$ is the infimum of all
central projections $q\in Z(\cM)$ such that $s(\omega)\leq q$.
\end{defn}
\begin{rem} \mlabel{rem:2.x}
(1) Let $\cA$ be a $C^*$-algebra and $\cA^{**}$ be its enveloping $W^*$-algebra.
Realizing any state $\omega$ of $\cA$ as a normal
state of the $W^*$-algebra $\cA^{**}$, we define
$s(\omega),\,z(\omega) \in \cA^{**}$ as above.
For a non-degenerate representation $(\pi,\cH)$ of $\cA$, we write
$\pi^{**} \: \cA^{**} \to \cB(\cH)$ for the corresponding weakly continuous
representation of $\cA^{**}$ extending~$\pi$.
 Then
$\ker \pi_\omega^{**} = \cA^{**}(\1 - z(\omega))$ and
$\pi_\omega^{**}(\cA^{**}) \cong z(\omega)\cA^{**}$ contains $s(\omega)$.\\[2mm]
(2) If $\cA$ is already a $W^*$-algebra, then for a normal state
 $\omega$ of $\cA$ we have that $s(\omega)\in\cA\subset \cA^{**}$
 by \cite[Lemma~2.132]{AS01}, and thus the two definitions for
  $s(\omega)$ coincide.
Note that if  $\cM$ is a $W^*$-algebra, and $\omega \in \fS_n(\cM)$,
then $\omega$ is faithful on $\pi_\omega(\cM)\cong z(\omega)\cM$ if and only if
$s(\omega)=z(\omega)$.\\[2mm]
(3)The relation of the central support $z(\omega)$ to the GNS-representation
$\pi_\omega$ can be generalized to any (non-degenerate) representation
$(\pi,\cH)$ of $\cA$ (cf. \cite[3.8.1]{Pe89}) as follows.
Define the   {\it central support (or cover)} of $(\pi,\cH)$, denoted $z(\pi),$ as the unique
central projection such that $z(\pi)\cA^{**}\cong \pi^{**}(\cA)'' (=\pi^{**}(\cA^{**})).$
By \cite[Theorem~3.8.2]{Pe89}, $z(\pi)$ determines $\pi$ up to quasi-equivalence.
As every folium $F \subeq \fS(\cA)$ determines a representation $(\pi,\cH)$ up to
quasi-equivalence, for which it is the set of normal states, we can define the
{\it central support (or cover)} of $F$, denoted $z(F),$ as $z(F):=z(\pi)$
where $(\pi,\cH)$ is the representation determined by $F.$
\end{rem}

For the next theorem we recall the Murray--von Neumann equivalence
relation $\sim$ on the set $\Proj(\cM)$ of projections in a $W^*$-algebra~$\cM$.
It is defined by $P \sim Q$ if and only if there exists a $V \in \cM$ with
$V^*V = P$ and $VV^* = Q$. We write $[\Proj(\cM)]$ for the set of equivalence classes
of projections.

\begin{thm}\mlabel{thm:cycl} {\rm(Equivalence Theorem for cyclic representations)}
Let $\cA$ be a $C^*$-algebra.
For states $\phi,\psi \in \fS(\cA)$, the corresponding
cyclic representations $(\pi_\phi, \cH_\phi)$ and $(\pi_\psi, \cH_\psi)$
are unitarily equivalent if and only if $s(\phi) \sim s(\psi)$,
i.e.\ their support projections are equivalent in $\cA^{**}$ in the sense of Murray--von Neumann.
\end{thm}

\begin{prf}
Given $\phi,\psi \in \fS(\cA)$,
recall from \cite[Cor.~V.1.11]{Ta02} that
$\pi_\phi$ is unitarily equivalent to a subrepresentation
of $\pi_\psi$, denoted $\pi_\phi \preccurlyeq \pi_\psi$,
 if and only if $s(\phi) \lesssim s(\psi),$ i.e.\
$s(\phi)$ is equivalent to a subprojection of $s(\psi)$.
Thus
\[\pi_\phi \preccurlyeq \pi_\psi \quad \Longleftrightarrow \quad
s(\phi) \lesssim s(\psi)\]
and thus
\[\pi_\phi \cong \pi_\psi \quad \Longleftrightarrow \quad
s(\phi) \sim s(\psi).\]
Here we use that $s(\phi) \lesssim s(\psi) \lesssim s(\phi)$ is equivalent to
$s(\phi) \sim s(\psi)$ (\cite[Prop.~V.1.3]{Ta02}) and
$\pi_\phi \preccurlyeq  \pi_\psi \preccurlyeq \pi_\phi$ is equivalent to
$\pi_\phi \cong \pi_\psi$ (\cite[Cor.~5.1.5]{Dix77}).
\end{prf}

\begin{rem} (a) An analogous statement holds for the central support
projections;-
for states $\phi,\psi \in \fS(\cA)$, we have that
$z(\phi) = z(\psi)$ if and only if $\pi_\phi$ and $\pi_\psi$ are
quasi-equivalent (cf. \cite[Prop.~5.10, Eq.~(5.6)]{AS01},
or \cite[Thm.~3.8.2]{Pe89}).
Thus, by Remark~\ref{rem:2.15}(d), their folia are equal
$F(\pi_\phi) = F(\pi_\psi)$.

(b)
Below we will present an alternative proof
of Theorem~\ref{thm:cycl}
based on standard representations in Theorem~\ref{thm:cycl2}.

(c) That $s(\omega) = \1$ means that $\omega$ is a faithful state,
i.e.\  $\Omega_\omega \in \cH_\omega$ is separating. So one particular
case of the preceding theorem is the fact that if
$\phi$ and $\psi$ are faithful, then $\pi_\phi \cong \pi_\psi$
(cf.~\cite[Thm.~III.2.6.7]{Bla06}).

(d) If $\cA$ is unital, then for
two pure states $\phi$ and $\psi$,
their GNS representations are equivalent if and only if $\phi(A) = \psi(UAU^{-1})$ for some
$U \in \cU(\cA)$ and all $A \in \cA$ (cf.\ \cite[Thm.~5.19]{AS01}).

(e) A set of states of which the support projections are equivalent, has a differential geometric structure.
This is studied in \cite{AV05}, \cite{ACS00}, and \cite{ACS01}.
\end{rem}

\begin{cor}
\label{autimplement}
Let $\cA$ be a $C^*$-algebra.
 For $\omega \in \fS(\cA)$,
an automorphism $\alpha\in \Aut(\cA)$ can be implemented in
$\cH_\omega$, i.e.\ there exists $U \in \cU(\cH_\omega)$ with
\[ U \pi_\omega(A) U^* = \pi_\omega(\alpha(A)) \quad \mbox{ for } \quad A \in \cA,\]
if and only if  $\alpha(s(\omega)) \sim s(\omega)$.
\end{cor}

\begin{prf} The implementability of $\alpha$ is equivalent to
$\pi_\omega \cong \pi_\omega \circ \alpha
\cong \pi_{\omega \circ \alpha}$ and hence to
$ s(\omega) \sim s(\omega \circ \alpha)$ by Theorem~\ref{thm:cycl}.
As $s(\omega \circ \alpha) = \alpha^{-1}(s(\omega))$, the claim follows.
\end{prf}

The equivalence  $\alpha(s(\omega)) \sim s(\omega)$ is in $\cA^{**}$ when $\cA$ is
a $C^*$-algebra, but if $\cA$ is a $W^*$-algebra, then by Remark~\ref{rem:2.x}
$s(\omega)\in\cA\subset \cA^{**}$ and hence the equivalence
 $\alpha(s(\omega)) \sim s(\omega)$ is in $\cA$.
 The next example applies these concepts concretely.

\begin{ex} Let $G = \R$ and
$\cM := L^\infty(\R,M_2(\C))=L^\infty(\R)\,\overline{\otimes}\, M_2(\C)$ with the natural
$\R$-action $\alpha$ by translation. It has a representation
$\rho:\cM\to\cB\big(L^2(\R,\C^2)\big)$ by pointwise matrix multiplication,
and $Z(\cM)=L^\infty(\R)\otimes \1$.
A projection $P \in \cM$ can be represented
by a measurable function $P \: \R \to M_2(\C)$ whose range consists of
projections in $M_2(\C)$. For projections in $\cM$, the relation
 $P \sim Q$ is equivalent to $\tr P = \tr Q$ in
$L^\infty(\R)$. Let $(E_{ij})_{1 \leq i,j \leq 2}$ in $M_2(\C)$ denote the standard matrix basis.

(a) For $f \in L^1(\R,\R)$ with $0 < f(x)$ for  all $x \in \R$ and
$\int_\R f(x)\, dx = 1$, we consider the state
\[ \omega(B) := \int_\R f(x) B_{11}(x)\, dx,
\qquad\hbox{where}\qquad B(x)=\left(\begin{matrix}B_{11}(x) & B_{12}(x)\\
B_{21}(x) & B_{22}(x)\end{matrix}\right).
\]
As the $\rho\hbox{-cyclic}$ vector $v(x):=\sqrt{f(x)}\left({1\atop 0}\right)$
in $L^2(\R,\C^2)$
produces the state $\omega(B)={\la v,\,\rho(B)v\ra}$, there is a unitary
$W:\cH_\omega \to L^2(\R,\C^2)$ which intertwines $\rho$ and $\pi_\omega$.
The support projection of $\omega$ is
$s(\omega)(x) = E_{11}$ and its central support $z(\omega)$ is $\1$.
Both are translation invariant, hence so are their equivalence classes.
Thus all $\alpha_t$ are implementable in $\pi_\omega$, and,
as $\rho$ is a product representation, it is easily seen to be
covariant, using the implementers $U_t\otimes\1$ on $L^2(\R)\otimes\C^2$,
where $U_t$ is translation.

(b) Now we consider a state of the form
\begin{align*}
\omega(B) &:=
\int_{-\infty}^0 f(x) B_{11}(x)\, dx +
\int_0^{\infty} g(x) \tr(B(x))\, dx \\
&=\int_\R \big(f(x)+g(x)\big)B_{11}(x)\, dx +
\int_0^{\infty} g(x) B_{22}(x)\, dx\\
& \mbox{ for } \quad
0 < f \in L^1((-\infty,0)),\;
0 < g \in L^1((0,\infty)) \quad \mbox{ with } \quad
\int_\R f = \frac{1}{2}, \quad
\int_\R g = \frac{1}{4}.
\end{align*}
Then $s(\omega) = E_{11} + \chi_{\R_+} E_{22}$ is not translation invariant
but $z(\omega) = \1$ is.
Then $\pi_\omega$ is not covariant because
$\tr s(\omega) = 1 + \chi_{\R_+}$ is not translation invariant.
If $\gamma$ is the representation of $\cM$ on
${L^2((-\infty, 0),\C^2) \oplus L^2([0,\infty),M_2(\C))}$ by matrix multiplication
(equipping $M_2(\C)$ with the inner product $\la C,D\ra:=\tr(CD^*)$), then
$\omega$ is the vector state obtained from the $\gamma\hbox{-cyclic}$ vector
\[
w(x):=\sqrt{f(x)}\left({1\atop 0}\right)\oplus \sqrt{g(x)}I.
\]
Thus  there is a unitary
$V:\cH_\omega \to {L^2((-\infty, 0),\C^2) \oplus L^2([0,\infty),M_2(\C))}$
which intertwines $\gamma$ and~$\pi_\omega$.
\end{ex}

\begin{rem}
 Let $(\cA,G, \alpha)$ be a $C^*$-action
 and $\omega \in \fS_\alpha(\cA)$. We would
 like to characterize situations when
$\pi_\omega$ is actually covariant, i.e.\  the $G$-action can be implemented on
$\cH_\omega$ by a continuous unitary representation.
By Theorem~\ref{thm:2.7} we need to assume at least that
$\pi_\omega$ is quasi-covariant, i.e.\ that
 $\Fol(\omega) = F(\pi_\omega)$ is $\alpha_G^*$-invariant.
Then we obtain a $W^*$-dynamical system
$(\pi_\omega(\cA)'',G,\beta)$ and
$\omega$ extends naturally to a state on $\cM := \pi_\omega(\cA)''$.
The implementability problem for $\cA$ is equivalent to the corresponding
problem for the von Neumann algebra $\cM$, so that it suffices to deal with it on the
$W^*$-level.

As a next condition, one should require implementability of $\alpha_G$
in $\pi_\omega$. By Corollary~\ref{autimplement}
 it is necessary that
the equivalence class
\[ [s(\omega)] = \{ P \in \Proj(\cM) \,\mid\; P \sim s(\omega)\} \]
is invariant under $\beta_G$. Suppose that this is the case.
Then each $\beta_g \in \Aut(\cM)$ can be implemented in $\cH_\omega$.
To characterize whether there are implementers which combine to give a
group representation, hence a covariant representation, is a well--known
problem in group cohomology. One chooses a set of unitary implementers, e.g.\ let $U_g$ implement $\beta_g$. Then the discrepancy $\sigma$
with the group law, i.e.\ $U_gU_h=\sigma(g,h)U_{gh}$,
 produces a (non-commutative) $2$-cocycle with coefficients in the unitary
group $\U(\cM')$. If $\cM'$ is commutative (the representation of $\cM$ is
multiplicity free), then
one needs to characterize when the cocycle $\sigma$  is a coboundary
within a suitably continuous class of cochains. If the appropriate
second cohomology group is trivial, this  would give a sufficient condition for
obtaining a covariant representation.
In the case that $G$ is locally compact, this leads to the study of Moore cohomology
for the group (cf. \cite{Ro86}, \cite{MOW16}).

More specifically, we consider the group
\[  \hat G_\omega := \{ (g,U) \in G \times \U(\cH_\omega) |\
(\forall M \in \cM)\ U \pi_\omega(M) U^{-1} = \pi_\omega(\beta_g(M))\}.\]
Then $\hat G_\omega$ is a closed subgroup of $G \times \U(\cH_\omega)$
and the projection onto the second factor provides a continuous unitary
representation of $\hat G_\omega$ on $\cH_\omega$.
Since every $\beta_g$ is implementable on $\cH_\omega$, the map
\[ q \: \hat G_\omega \to G, \quad (g,U) \mapsto g \]
is surjective and its kernel is isomorphic to the unitary group
$\U(\pi_\omega(\cM)') \cong \U(\cM_{s(\omega)}')$.
We thus have a short exact sequence
\[ \be \to \U(\cM_{s(\omega)}') \to \hat G_\omega \to G \to \be.\]
The covariance of the representation $\pi_\omega$ is equivalent to
the splitting of this extension of topological groups.
\end{rem}
The question of covariance of $\pi_\omega$
for a  a $C^*$-action $(\cA,G, \alpha)$ with $\omega \in \fS_\alpha(\cA)$, given
unitary implementability, can be answered in a more restricted context
(cf.~\cite{Ka71}):

 \begin{thm}\mlabel{thm:Kallman}
{\rm(Kallman's Theorem)}
Let $\cA\subseteq\cB(\cH)$ be a unital $C^*$-algebra where $\cH$ is separable.
Let  $\alpha:\R\to\Aut \cA$ be a $C^*$-action such that
\begin{itemize}
\item[\rm(i)] $t\mapsto\alpha_t(A)$ is
weak operator continuous for each $A\in\cA$, and
\item[\rm(ii)] for each $t\in\R$ there is a unitary $U_t\in \cA''$
such that $\alpha_t={\rm Ad}\,U_t$ on $\cA$.
\end{itemize}
Then there is a strong operator continuous one parameter unitary group
$W:\R\to\cA''$ such that $\alpha_t={\rm Ad}\,W_t$ on $\cA$.
\end{thm}

 As an application of this, consider a  $C^*$-action $(\cA,\R, \alpha)$ where $\cA$ is unital and separable,
and let $\omega \in \fS(\cA)$. Let  $F(\omega)$ be $\alpha_{\R}^*$-invariant and contained in $(\cA^*)_c$
(cf. Corollary~\ref{cor:bor1}), so (i) is satisfied and $\cH_\omega$ is separable.
We can obtain (ii) by e.g. assuming $\omega$ is pure as all automorphisms of $\cB(\cH_\omega)=\pi_\omega(\cA)''$
are inner. Thus a pure state is covariant if and only if   $F(\omega)$ is $\alpha_{\R}^*$-invariant and contained in $(\cA^*)_c$,
as the converse follows from  Corollary~\ref{cor:bor1}.

As a second application of Kallman's Theorem, consider a  $C^*$-action $(\cA,\R, \alpha)$, and define the discrete
crossed product $\cA \rtimes_\alpha \R_d=:\cB$.
Observe that $\alpha$ extends to an action
on $\cB$ by $\alpha_t(B):= ({\rm Ad}\,\delta_t)(B)$ for $B\in\cB$
where $\delta_t\in\ell^1(\R,\cA)\subset\cB$ is the function with value $\1$ at $t$ and zero elsewhere
(note that $({\rm Ad}\,\delta_t)(\delta_s)=\delta_s$).
Let $(\pi,\cH)$ be a representation of $\cB$ for which $\cH$ is
separable. This corresponds to a covariant representation  of $(\cA,\R_d, \alpha)$.
Then, using the unitaries $\pi(\delta_t)\in\pi(\cB)$, we have satisfied (ii)
of Kallman's Theorem for $(\cB,\R, \alpha)$.
To satisfy (i), we need to also assume that
on $\cB$,  $F(\pi)$ is $\alpha_{\R}^*$-invariant and contained in $(\cB^*)_c$. It then follows that
  $(\pi,\cH)$ is a covariant representation of  $(\cB,\R, \alpha)$, and restriction to $\cA$ produces a
covariant representation of  $(\cA,\R, \alpha)$.

\subsection{Continuity properties of covariant representations}
\label{StrucCovRep}

Henceforth we assume that non-zero covariant representations exist for
a $C^*$-action $(\cA,G,\alpha)$.
In the \usual case for  $(\cA, G, \alpha)$ (Subsection~\ref{subsec:2.1}),
the entire covariant representation theory is carried by
the crossed product $\cA \rtimes_\alpha G$.
When we do not have the \usual case, it may still be possible to find
 a $C^*$-algebra which
can fulfill the role of the crossed product. This has already been analyzed in~\cite{GrN14},
and in a subsequent paper (\cite{GrN18}) we have continued this analysis in the presence of spectral conditions.
First, we consider natural structures associated with covariant representations.

There is a  universal covariant representation, obtained as follows.
\begin{defn}\mlabel{def:1.1d}
Given a $C^*$-action $(\cA, G, \alpha)$,
 cyclic representations of $\cA \rtimes_\alpha G_d$ are obtained from
states through the GNS construction. Let $\fS_{co}$
denote the set of those states
$\omega$ on $\cA\rtimes_{\alpha} G_d=C^*(\cA\cup U_G)$
whose GNS-representations 
restrict on $\cA$
and on the implementing unitaries $U_G$ to a covariant pair $(\pi_\omega,U^\omega)$
for $(\cA, G, \alpha)$, i.e.
$(\pi_\omega,U^\omega)\in{\rm Rep}(\alpha,\cH_\omega)$.
Note that if $G$ is nondiscrete, then some
states on $\cA\rtimes_{\alpha} G_d$ need not be in $\fS_{co}$, due to the continuity requirement for $U^\omega$.
This allows us to define the {\it universal  covariant
representation}
$(\pi_{co},U_{co})\in{\rm Rep}(\alpha,\cH_{co})$ by
\[
\pi_{co}:= \bigoplus_{\omega\in\fS_{co}}
\pi_\omega,\quad U_{co}:= \bigoplus_{\omega\in\fS_{co}} U^\omega
\quad\hbox{on}\quad  \cH_{co}=\bigoplus_{\omega\in\fS_{co}}\cH_\omega.
\]
This is non-trivial as long as $\fS_{co}\not=\emptyset.$
We obtain a canonical $W^*$-dynamical system \break
$\alpha^{co}:G\to \Aut(\cM_{co})$, where $\cM_{co}:=\pi_{co}(\cA)''$ and
$\alpha^{co}(g)=\Ad U_{co}(g)$.
\end{defn}
\begin{prop} \mlabel{prop:posrep}
Assume that the $C^*$-action $(\cA, G, \alpha)$ has non-zero covariant representations. Then
 $(\pi_{co},U_{co})\in{\rm Rep}(\alpha,\cH_{co})$
is non-zero, and the folium $F(\pi_{co})$ is the unique
folium in $(\cA^*)_c$
which is maximal in the sense that it contains all other folia in $(\cA^*)_c$.
Moreover $F(\pi_{co})=\fS_\alpha(\cA)$ and this folium is $G$-invariant.
 \end{prop}

\begin{prf} Any covariant representation corresponds to a representation of $\cA \rtimes_\alpha G_d$, and the cyclic subrepresentations for this $C^*$-algebra
are still covariant, hence $\fS_{co}\not=\emptyset$ and  $(\pi_{co},U_{co})\in{\rm Rep}(\alpha,\cH_{co})$
is non-zero. Moreover every covariant representation
$(\pi,U)$ of $(\cA, G, \alpha)$ is a direct sum of subrepresentations of
$(\pi_{co},U_{co})$, hence $F(\pi)\subseteq F(\pi_{co})$. Since every
folium $F \subeq (\cA^*)_c$ is contained in a $G$-invariant folium
$F_G := \Fol(\alpha_G^*F) \subeq (\cA^*)_c$
and $F_G = F(\pi)$ for some covariant representation $(\pi,U)$
(Borchers--Halpern Theorem~\ref{thm:bor}),
it follows that $F(\pi_{co})$ contains every folium in $(\cA^*)_c$.
Clearly, there is only one folium in $(\cA^*)_c$ with this property.

Further, Theorem~\ref{thm:2.7} implies that every
$\omega \in \fS_\alpha(\cA)$ is contained in the folium
$\Fol_G(\omega) \subeq (\cA^*)_c$, so that we also obtain the inclusion
$\fS_\alpha(\cA) \subeq F(\pi_{co})$.
Conversely, let $\omega\in F(\pi_{co})$, then by  $G$-invariance of $F(\pi_{co})$,
we have that $\alpha_G^*\omega\subset F(\pi_{co})$ and hence
$\Fol_G(\omega) \subeq F(\pi_{co})\subset (\cA^*)_c$.
Thus by the above characterization of $\fS_\alpha(\cA)$ (Theorem~\ref{thm:2.7}),
it follows that $\omega\in\fS_\alpha(\cA)$ (Theorem~\ref{thm:2.7}).
This proves the reverse inclusion, hence the equality
$\fS_\alpha(\cA)= F(\pi_{co})$.
\end{prf}

Note that in general it is not true that $(\cA^*)_c
= \pi_{co}(\cA)''_*$ by Example~\ref{AcNotFol}.

\begin{prop}
Given a $C^*$-action $(\cA, G, \alpha),$ let $\tau_1 \supeq \tau_2$ be group topologies on $G$.
If $(\pi,U)$ is a
covariant representation with respect to ~$\tau_1$, then it contains a $\tau_2$-covariant subrepresentation
$(\pi_{\tau_2},U_{\tau_2})$  which is maximal, in the sense that it contains
all other $\tau_2\hbox{-covariant}$ subrepresentations
of $(\pi,U)$.
 \end{prop}

\begin{prf} Given  a covariant representation $(\pi,U)$
with respect to  $\tau_1$ on $\cH$, we consider the closed subspace
$\cH_c$ of continuous vectors for the representation of the
topological group $(G,\tau_2)$. Then
$\cH_c$ is $G$-invariant and maximal with respect to the property that
the action of $(G,\tau_2)$ on this subspace is continuous.

Now let $\cH_2 := \{ \xi \in \cH_c \,\mid\, (\forall A \in \cA) \ \pi(A)\xi \in \cH_c\}$
be the maximal $\cA$-invariant subspace of $\cH_c$. Then $\cH_2$ is also
$G$-invariant because, for $g \in G$,  $\xi \in \cH_2$ and $A \in \cA$, we have
$\pi(A)U_g \xi = U_g \pi(\alpha_g^{-1}(A)) \xi \in U_g \cH_c = \cH_c$.
It is also clear that $\cH_2$ is maximal with respect to the
property that it carries a covariant representation of
$(\cA,(G,\tau_2), \alpha)$.
\end{prf}

If $\tau_1$ is the discrete topology and $\tau_2$ the given topology on $G$,
then the preceding proposition
implies that every covariant representation $(\pi, U)$ of $(\cA,G_d, \alpha)$
contains a maximal covariant subrepresentation for $(\cA,G,\alpha)$.
If the covariant subrepresentation is zero, we will call $(\pi, U)$ a {\it purely discontinuous}
covariant representation.
An  irreducible covariant representation of $(\cA,G_d,\alpha)$
is either covariant or purely discontinuous.

Given any  $(\cA, G, \alpha)$, we can always define the \strong continuous part of it by
\[
{\cal A}_c:=\{A\in{\cal A}\,\mid\, \alpha^A \: G \to \cA, g\mapsto\alpha_g(A)\quad\hbox{is norm continuous}\}
\quad \mbox{ and } \quad
\alpha^c_g:=\alpha_g\restriction{\cal A}_c.
\]
Unfortunately, as we will see in Example~\ref{Ac0} below, it is possible that $\cA_c = \C\1$.
As we have seen in the Borchers--Halpern Theorem, it is much more the continuous portion
$(\cA^*)_c$ of the $G$-action on $\cA^*$ than the continuous portion $\cA_c$ of $\cA$ that is
responsible for the covariant representations.

\begin{rem}
If we start from  a $W^*$-dynamical system $(\cM,G, \beta)$ with $G$ locally compact,
 then ${\cal M}_c$ is weakly dense in $\cal M$, and
\[\cM_c=C^*\big\{\beta_f(A)\mid\, f\in L^1(G),\, A\in\cM\big\},\]
where the integrals $\beta_f(A) := \int_G f(g)\beta_g(M)\, dg$
exist in the weak topology  (\cite[Lemma~7.5.1]{Pe89}).
Thus, associated with any $C^*$-action $(\cA,G, \alpha)$, there is a $C^*$-dynamical system
$(\cM_{co,c},G, \beta^{co,c})$, which in the locally compact case encodes the covariant representation theory of
$(\cA,G, \alpha)$. As remarked, it is possible that $\cM_{co,c}$
intersects $\pi_{co}(\cA)$ only in $\C\1$, though in the \usual case $\cM_{co,c}\supseteq\pi_{co}(\cA)$
where the inclusion may be proper.
\end{rem}

\begin{rem}
In general, $\cM_{co} := \pi_{co}(\cA)''$ produces the $W^*$-dynamical system
$(\cM_{co}, G, \alpha^{co})$ whose covariant normal representations are in one-to-one correspondence
with the covariant representations of $(\cA,G,\alpha)$ because
the $\fS_n(\cM_{co})$ can be identified with $\fS_\alpha(\cA)$
(Proposition~\ref{prop:posrep}). This $W^*$-dynamical system is therefore a suitable tool
to analyze covariant representations of a given (possibly singular)
$C^*$-action $(\cA,G, \alpha)$.
\end{rem}

We list a few examples which will be useful for subsequent discussion. The reader in a hurry can proceed to the next subsection.

\begin{ex}
\label{Ac0} {\rm(A case of $\cA_c = \C \1$ and $\omega \in \fS(\cA)^G \setminus \fS_\alpha(\cA)$)}
 We consider the rotation action of $\T$ on the abelian group
$(\C,+)$ by multiplication. This produces an action of $G = \R$ on the Weyl algebra
$\cA := \ccr \C ,\sigma.$, where
$\sigma(z,w) = \Im(\oline z w)$, by
\[ \alpha_t(\delta_z) = \delta_{e^{it}z}, \quad t \in \R, z \in \C.\]

We claim that $\cA_c = \C \1$. To verify this claim, we consider the
covariant representation $(\pi_\omega, U_\omega, \cH_\omega)$ of $\cA$ on
$\ell^2(\C) \cong \cH_\omega \subeq \cA^*$ corresponding to the $\alpha$-invariant
tracial state $\omega$ defined by  $\omega(\delta_z) = \delta_{0,z}$, for which
$U_\omega$ fixes the cyclic vector $\delta_0$.

Since $\cA$ is simple by \cite[Thm.~5.2.8]{BR96},
the state $\omega$ is faithful by Lemma~\ref{lem:3.2b} below.
Therefore
\[ \eta \: \cA \to \ell^2(\C), \quad \eta(A) = A \delta_0 \]
is a faithful continuous injection mapping the generator
$\delta_z$ of $\cA$ to the basis element $\delta_z = (\delta_{z,w})_{w \in \C} \in \ell^2(\C)$.
Note that $\eta$ is equivariant with respect to the action
$\alpha$ of $\T$ on $\cA$ and the representation $U$ of $\T$ on $\ell^2(\C)$
defined by the permutation of the generators
\[ U_t \delta_z = \delta_{e^{it}z}, \quad t \in \R,z \in \C.\]
Lemma~\ref{lem:3.1b} implies that $\ell^2(\C)_c = \C \delta_0$
for the unitary one-parameter group $U$ which in particular entails that
$\omega \in \fS(\cA)^\T \setminus \fS_\alpha(\cA)$.
The continuity of the inclusion $\eta \: \cA \to \ell^2(\C)$
now yields $\cA_c = \C \1$ for~$\alpha$.
Nevertheless, the Schr\"odinger representation is an
example of a faithful covariant representation for $\alpha$.
\end{ex}

With a similar argument as in the previous example, we even obtain
an example where $\cA$ is commutative.

\begin{ex}  {\rm(A case of $\cA_c = \C \1$, $\omega \in \fS(\cA)^G \setminus \fS_\alpha(\cA)$ and
$\cA$ commutative)}
We consider the rotation action of $\T$ on the abelian group
$(\C,+)$ by multiplication and the $C^*$-algebra
$\cA := C^*(\C_d)$, where $\C_d$ is the discrete additive group of complex numbers.
We thus obtain an action of $\T$ on $\cA$ by
\[ \alpha_t(\delta_z) = \delta_{tz}, \quad t \in \T, z \in \C.\]

(a) We claim that $\cA_c = \C \1$. To verify this claim, we consider the
faithful covariant representation $(\pi_\omega, U_\omega, \cH_\omega)$ of $\cA$ on
$\ell^2(\C) \cong \cH_\omega \subeq \cA^*$ corresponding to the $\alpha$-invariant
state $\omega$ defined by $\omega(\delta_z) = \delta_{0,z}$, for which
$U_\omega$ fixes the cyclic vector $\delta_0\in \ell^2(\C)$.

Since $\cA$ is commutative, the annihilator of the state $\omega$ coincides with
$\ker \pi_\omega$ (Lemma~\ref{lem:3.2b}).
Further, the amenability of the discrete
abelian group $\C_d$ shows that the representation
of $\cA = C^*(\C_d)$ on $\cH_\omega \cong \ell^2(\C_d)$ is faithful (\cite[Thm.~7.3.9]{Pe89}).
Now $\eta \: \cA \to \ell^2(\C), \eta(A) = A \delta_0$
is a continuous linear injection mapping the generator
$\delta_z$ of $\cA$ to the basis element $\delta_z \in \ell^2(\C)$.
Therefore $\eta$ is equivariant with respect to the action
$\alpha$ of $\T$ on $\cA$ and the representation $U$ of $\T$ on $\ell^2(\C)$
defined by the permutation of the generators
\[ U_t \delta_z = \delta_{tz}, \quad t \in \T,z \in \C.\]
From Lemma~\ref{lem:3.1b} we know that $\ell^2(\C)_c = \C \delta_0$
for the unitary one-parameter group $U$, so that
the continuity of the inclusion $\eta \: \cA \to \ell^2(\C)$
implies that $\cA_c = \C \1$.

(b) Note that $C^*(\C_d) \cong C_0(\hat{\C_d})$ and that the compact
group $b\C := \hat{\C_d} = \Hom(\C_d,\T)$
is the {\it Bohr compactification} of the locally compact
abelian group $\C$.
Since the canonical image of $\C$ in $b\C$ is dense, the $C^*$-algebra $\cA$ embeds naturally
into $C^b(\C)$, which in turn injects into the algebra
 $\prod_{r > 0} C(r\T)$ by restricting
to circles of radius $r > 0$.
By considering the $L^2\hbox{--space}$ of a measure $\mu$
on $\C$ concentrated on $r\T$ where it is the invariant measure,
we  obtain for any $r > 0$ a covariant
representation $(\pi_r, U_r)$ of $(\cA,\R,\alpha)$. These representations
separate the points of $\cA$. By taking their direct sum we obtain
an example of a $C^*$-action $(\cA,\R,\alpha)$ with
$\cA_c  = \C \1$ and a faithful covariant representation~$(\pi,U)$.

Note that for the  associated $W^*$-dynamical system
$\beta:G\to \Aut(\cM)$, where $\cM:=\pi(\cA)''$ and
$\beta(g)=\Ad U(g)$, we do obtain a subalgebra
$\cM_c$ which is strong operator dense in $\cM$ (\cite[Lemma~7.5.1]{Pe89}),
but by the preceding, $\cM_c$ intersects $\pi(\cA)$ only in $\C\1$.
\end{ex}

\begin{ex} \mlabel{ex:3.6} Let $\cA := \ell^\infty(\T_d)$ with pointwise operations
and $\|\cdot\|_\infty$-topology. Here $\T_d$ denotes the discrete group underlying the
circle group $\T\subeq \C^\times$.
Consider the action
of the topological group $G = \T$ on $\cA$ by rotation, i.e.\
$\alpha_t(\delta_z) = \delta_{tz}$ for $t,z \in \T$.

Let $t_0 \in \T$ be an element of infinite order, so that
the the subgroup $T_0 \subeq \T_d$ generated by $t_0$ is infinite.
Define a character
\[ \xi_0 \: T_0 \to \T, \quad \xi_0(t_0^n) := (-1)^n\]
and let $\xi\: \T_d \to \T$ be an extension of this character to all of $\T_d$
(cf.\ \cite[Prop.~A1.35]{HM13}).
Then
\[ S := \xi^{-1}(\{ e^{it}  \: 0 \leq t < \pi\})\subset\T \]
satisfies $\T = S \dot\cup t_0 S$ because
$\xi(t_0) = -1$.

Since the abelian discrete group $\T_d$ is amenable,
there exists an invariant mean $\omega \in \cA^* = \ell^\infty(\T_d)^*$
which is an $\alpha$-invariant state on $\cA$. It corresponds to a finitely additive measure $\mu$ on $\T$ by
$\omega(\chi_E)=\mu(E)$. As $1=\omega(\1)=\omega(\chi_S+\chi_{t_0S}) = \mu(S) + \mu(t_0 S) = 2 \mu(S)$,
invariance of $\omega$ now implies that $\mu(S) = \shalf$, so that, for
$A = \chi_S \in \cA$, we obtain
\[ \omega(A \alpha_{t_0^n}(A)) = \mu(S \cap t_0^n S)
= \frac{1}{4}(1 + (-1)^n).\]
As there are elements from both $S$ and $t_0S$ arbitrarily close to $1\in\T$,
this implies that the function
$t \mapsto \omega(A \alpha_t(A))$ is not continuous on $\T$.
We conclude that
\[ \omega \in \fS(\cA)^\T \setminus \fS_\alpha(\cA).\]
In this example $\cA_c = C(\T)$ is strictly larger than $\C \1$.
\end{ex}

\subsection{Innerness for covariant representations}
\mlabel{InnCovRep}

Given a $C^*$-action $(\cA, G, \alpha),$ a desirable property for a covariant representation
 $(\pi,U)$ is that it is inner, i.e.\ $U_G\subset\pi(\cA)''$. This is desirable from a physical
 point of view, as it implies that the generators of the one-parameter groups in $G$ are
 affiliated with  $\pi(\cA)''$, hence are observables. It also is a peculiarly quantum requirement,
 as in the case that $\cA$ is commutative and $\alpha$ is nontrivial, then there are no inner
 covariant representations which are faithful on $\cA$. Below in Sect~\ref{BA-sect} we will see the surprising
 fact that certain spectral conditions guarantee the existence of inner
 covariant representations (Borchers--Arveson Theorem~\ref{BA-thm}).
Even in the absence of spectral conditions, the
 innerness of covariant representations have been analyzed. Moreover, we saw above in Kallman's Theorem~\ref{thm:Kallman}
 that in some circumstances, innerness of the action implies covariance.

 A great deal is known about $W^*$-actions on factors, starting from the simple observation that all  automorphisms of type~I factors are inner (as they are isomorphic to some $\cB(\cH)$). However,
 for the case where the von Neumann algebra is unrestricted, the best result seems to be
 the one in Kraus \cite[Theorem 3.2]{Kr79}, which we state below. We first need to fix some notation.

 Given a $W^*$-dynamical system $(\cM,G,\beta)$ such that $G$ is  locally compact and abelian,
 denote the fixed point algebra by $\cM^\beta$ and the center of $\cM$ by $Z(\cM)$.
 For any projection $P\in\cM^\beta$, the reduced von Neumann algebra $\cM_P=P\cM P$ is  left invariant by
 $\beta$, hence we can restrict the action to obtain a new action  $(\cM_P,G,\beta^P)$. Let $f\in L^1(G)$, which
 we recall has Fourier transform
 \[
 \hat{f}(\gamma):=\int_G f(g)\gamma(g)\,dg\qquad\hbox{for}\quad \gamma\in\hat G
 \]
 where $\hat G$ is the dual group. If we define
 $\beta_f\in\cB(\cM)$ by
 \[
 \beta_f(M):=\int_G f(g)\,\beta_g(M)\,dg\,,
 \]
 then ${\rm Spec}(\beta)$ denotes the support of the map $ \hat{f}\mapsto\beta_f$, i.e.
 \begin{equation}
 \label{Spec1def}
 {\rm Spec}(\beta)=\{\gamma\in\hat G\,\mid\,(\forall f \in L^1(G))\
\beta_f = 0 \ \Rarrow \ \hat{f}(\gamma)=0\}.
  \end{equation}
 Then for $P\in\cM^\beta$ we have ${\rm Spec}(\beta^P)\subset{\rm Spec}(\beta)$.
 \begin{thm}\mlabel{thm:LCInner}
For a $W^*$-dynamical system $(\cM,G,\beta)$
such that $G$ is connected, locally compact and abelian,
the following are equivalent:
\begin{itemize}
\item[\rm(i)] $\beta$ is inner.
\item[\rm(ii)] For every nonzero projection $P\in Z(\cM^\beta)$ and every
compact neighborhood $V$ of $0$
in the dual group~$\hat G$, there is a nonzero projection  $Q\in Z(\cM^\beta)$ such that
$Q\leq P$ and ${\rm Spec}(\beta^Q)\subset V$.
\item[\rm(ii)] For every nonzero projection $P\in Z(\cM^\beta)$,
 there is a nonzero projection  $Q\in Z(\cM^\beta)$ such that
$Q\leq P$ and ${\rm Spec}(\beta^Q)$ is compact.
\end{itemize}
\end{thm}
This is proven in \cite[Thm.~3.2]{Kr79}. Thus innerness is characterized by the
projections in the center of the fixed point algebra. Note that, by
\cite[Prop.~3.2.41]{BR02}, the condition that ${\rm Spec}(\beta^Q)$ is compact
is equivalent to the norm continuity of $\beta^Q:G\to\Aut(\cM_P)$.

\section{Standard and $P$-standard representations of $W^*$-algebras}
\mlabel{sec:3}

We saw above in Corollary~\ref{cor:bor1} of the Borchers--Halpern Theorem,
that a representation $\pi$ of $\cA$ is quasi-covariant if and only if
the $C^*$-action $(\cA,G, \alpha)$ extends to a  $W^*$-dynamical system
$(\pi(\cA)'',G,\beta)$.
However, no indication was  given on how to construct the covariant representation
needed for the quasi-covariance.
In this section  we want to address this question, i.e.
if one is given such a  $W^*$-dynamical system, how we can construct
a faithful normal representation of it which is covariant.
This will be done through standard forms of $W^*$-algebras (defined below).
Standard forms also occur naturally as the GNS representations of
KMS states. Below in Theorem~\ref{thm:Pinv} we will show that associated to
every invariant nonzero projection, there is a normal Stinespring dilation representation
which is covariant. This provides an interesting source of covariant representations.

\subsection{Standard forms of $W^*$-algebras}
\mlabel{SectStandardForm}

Given a  von Neumann algebra $\cM \subeq \cB(\cH)$, we recall next what its standard form representation is. This
can be defined either constructively, or by abstract characterization of its structure, and by uniqueness there
is only one such standard form representation, up to unitary equivalence.
We first state the structural definition.
\begin{defn} (\cite{Haa75}) \mlabel{def:a.1}
(a) A von Neumann algebra $\cM \subeq \cB(\cH)$ is said
to be {\it in standard form} if there exist an anti-unitary involution $J$ on $\cH$ and a cone $\cC \subeq \cH$
which is self-dual in the sense that
\[
\cC=\{\psi\in\cH\,\mid (\forall \xi \in \cC)\, \la \psi,\xi \ra \geq 0 \}
\]
and $\cM,\, J,\, \cC$ satisfy:
\begin{itemize}
\item[\rm(S1)] $J\cM J = \cM'$.
\item[\rm(S2)] $J\psi = \psi$ for every $\psi \in \cC$.
\item[\rm(S3)] $AJA\cC \subeq \cC$ for all $A \in \cM$.
\item[\rm(S4)] $JAJ = A^*$ for all $A \in \cM \cap \cM'$.
\end{itemize}
A von Neumann algebra in standard form is denoted by $(\cM, \cH, J,\cC)$.

(b) A normal representation $(\pi, \cH)$ of a $W^*$-algebra is called a
{\it standard (form) representation}
if there exist $J$ and $\cC$ such that $(\pi(\cM),\cH,J,\cC)$ is a von Neumann
algebra in standard form.

(c) From the definition it follows that  $(\cM, \cH, J,\cC)$ is in standard form
if and only if $(\cM', \cH, J,\cC)$ is in standard form.
\end{defn}

\begin{rem}
\label{OppAlg}
(1) For a von Neumann algebra in standard form, the map
\[ \cM \to \cM', \quad M \mapsto JM^* J \]
induces an isomorphism of $W^*$-algebras $\cM^{\rm op} \to \cM'$.
(The opposite algebra  $\cM^{\rm op}$
is the space $\cM$ equipped with
the previous multiplication but where the order of terms are reversed and all other
operations, including the scalar multiplication, are the same.
It is
isomorphic to the complex conjugate algebra, via the map $M\mapsto M^*$.)\\[2mm]
(2) We can also give a constructive definition of  a standard form representation
\break (cf.~\cite[Def.~III.2.6.1]{Bla06}).
It states that  $\cM $ is in standard form, if
 it is unitarily equivalent to the
GNS representation of a faithful normal semifinite weight on $\cM$
(here
normal means lower semicontinuous with respect to ~the ultraweak topology on $\cM_+$).
Recall that a {\it weight} $w$ on $\cM$ is
{\it semifinite} if the set
\[
\{M\in\cM_+\,\mid\,w(M)<\infty\}
\]
generates a *-algebra which is $\sigma(\cM,\cM_*)\hbox{--dense}$  in $\cM$.
Every von Neumann algebra has a faithful normal semifinite weight
(cf. \cite[III.2.2.26]{Bla06}),
in fact all normal semifinite weights are obtained as sums of normal positive forms
(cf. \cite{Haa75b}). As a consequence, standard form representations exist.
\end{rem}
We now state three equivalent characterizations of a standard form representation.

\begin{thm}
\mlabel{Def-SF} For a von Neumann algebra $\cM \subeq \cB(\cH)$, the
following are equivalent:
\begin{itemize}
\item[\rm(i)] $\cM$ is in standard form,
\item[\rm(ii)] There is a  faithful normal semifinite weight on $\cM$  such
that its GNS representation is unitarily equivalent to the inclusion
$\cM \into \cB(\cH)$,
\item[\rm(iii)] There is an anti-unitary involution $J$ on $\cH$ such that
$J\cM J=\cM'$ and $JZJ=Z^*$ for every $Z\in Z(\cM)$.
\end{itemize}
Moreover, a standard form representation exists for any von Neumann algebra,
and it is unique up to unitary equivalence.
\end{thm}
\begin{prf}
The equivalence of (i) and (ii) is in \cite[Thm.~IX.1.2]{Ta03}, and the
equivalence of (ii) and (iii) is in \cite[Thm.~III.4.5.7]{Bla06},
which also states the existence and uniqueness.
This is also in \cite[Thm.~1.6]{Haa75} and \cite[Thm.~IX.1.14]{Ta03}.
\end{prf}

\begin{rem}
  \mlabel{rem:2.6}
(a) It is of central importance that for a von Neumann algebra $\cM$, its
faithful standard form representations are unitarily equivalent.
Moreover, every cyclic normal representation of $\cM$ is contained in its
standard form representation (this has a generalization below in
Proposition~\ref{prop:3.13}). Note that, by condition (iii), the
additional structure in Definition~\ref{def:a.1} is automatic,
given $J$. Thus if $\cM \subeq \cB(\cH)$ is standard and commutative,
it is a maximal abelian subalgebra of~$\cB(\cH)$.

(b) For any von Neumann algebra in standard form $(\cM,\cH,J,\cC)$, the map
\[ \{ \xi \in \cC \,\mid\, \|\xi\| =  1\} \to \fS_n(\cM), \quad
\xi \mapsto \omega_\xi, \quad \omega_\xi(A) := \la \xi, \pi(A)\xi \ra \]
is a homeomorphism for the norm topology on $\fS_n(\cM)$ (\cite[Lemma~2.10]{Haa75}).

(c) By  \cite[Thm.~5.3.10]{BR96}, we see that a normal KMS state with respect to ~any $W^*$-dynamical
system is faithful, hence it is a faithful normal semifinite weight on $\cM$,
hence its GNS representation is in standard form. Below in Section~\ref{KMSsect}
 we will consider KMS states in greater detail.
\end{rem}

\begin{ex}\mlabel{ex:3.4}
(a) Let $(X,\fS,\mu)$ be a {\it semifinite measure space}, i.e.\
for each $E \in \fS$ with $\mu(E) = \infty$, there exists a measurable subset
$F \subeq E$ satisfying $0 < \mu(F) < \infty$.
Then the multiplication action of $\cM := L^\infty(X,\fS,\mu)$ on
$\cH := L^2(X,\fS,\mu)$ realizes $\cM$ as a von Neumann algebra in standard
form $(\cM,\cH,J,\cC)$. Here
$Jf =\oline f$ and $\cC = \{f \in L^2(X,\fS,\mu) \,\mid\, 0 \leq f \}$.
Any element of $\cH$ vanishing only in a zero set is a cyclic separating vector,
and such elements exist if and only if $\mu$ is $\sigma$-finite.

(b) If $\cK$ is a complex Hilbert space, $\cM = \cB(\cK)$ and $\cH := \cB_2(\cK)$
is the Hilbert space of Hilbert--Schmidt operators
on $\cK$, then the left multiplication representation of $\cM$ on $\cH$
yields a standard form $(\cM,\cH,J,\cC)$, where
$J(A) = A^*$ and $\cC = \{ A \in \cB_2(\cK) \,\mid\, 0 \leq A\}$.
A cyclic vector exists if and only if $\cK$ is separable.
\end{ex}

The faithful normal semifinite weights on $\cM$ include the faithful normal states, if any exist.
The existence of a faithful normal state, is equivalent to the property that $\cM\subset  \cB(\cH)$
is {\it countably decomposable}, i.e.\  every mutually orthogonal
family of projections in $\cM$ is at most countable (cf. \cite[Prop.~III.2.2.27]{Bla06}).
 This is the case if $\cH$ is separable.
 Given a faithful normal state, then
 its GNS representation has a cyclic separating vector,
hence we can apply the Tomita--Takesaki modular theory in the GNS representation.
This is directly connected to the standard form structures by:

\begin{prop} \mlabel{prop:2.9} {\rm(\cite[Thm.~1.1, Rem.~1.2]{Haa75})}
If $\cM \subeq \cB(\cH)$ is a von Neumann algebra with a cyclic separating vector
$\Omega\in \cH$, then the corresponding modular involution $J$ leads to a standard form
realization $(\cM,\cH,J,\cC)$, where $\cC\subeq \cH$ is the closed convex cone generated by
the elements $AJA\Omega$, $A \in \cM$.
\end{prop}

\begin{rem}
  \mlabel{rem:aut-w*}
(a) As we are concerned with $W^*$-dynamical systems $\beta:G\to\Aut(\cM)$, the following property
of the standard form representation $(\cM, \cH, J,\cC)$ is of central importance to us.
The group
\[ \cU(\cH)_\cM := \{ U \in \cU(\cH) \,\mid\, U\cC \subeq \cC, \;UJ=JU,\; U \cM U^* = \cM\}\]
has a natural homomorphism to $\Aut(\cM)$ by conjugation:
\begin{equation}
  \label{eq:standard-auto}
\Gamma \: \cU(\cH)_\cM \to \Aut(\cM), \quad
\Gamma(U)(M) := UMU^*,
\end{equation}
and this homomorphism is an isomorphism of topological groups with respect to the
$u$-topology on $\Aut(\cM)$ and the strong operator topology on $\cU(\cH)_\cM$
(\cite[Prop.~3.5]{Haa75}, \cite[Thm.~IX.1.15]{Ta03},
\cite[\S 2.23]{Str81}).\begin{footnote}{In \cite[Thm.~14]{Pi06}
it is asserted that the strong operator topology
on $\cU(\cH)_\cM$ corresponds to the $p$-topology on $\Aut(\cM)$,
but this is inconsistent with \cite[Rem.~3.9]{Haa75}
and contradicts Example~\ref{ex:linft}.}\end{footnote}
In particular, the whole group $\Aut(\cM)$ has a natural unitary representation
on $\cH$.

(b) As $\cU(\cH)_\cM$ commutes with $J$, we have
$\cU(\cH)_\cM\cap \cM\subset \cU(\cM\cap\cM')$, so that
$\cU(\cH)_\cM\cap \cM\subset \ker \Gamma = \{\be\}$. Therefore
$\cU(\cH)_\cM$ intersects $\cU(\cM)$ and $\cU(\cM')$ trivially.
\end{rem}

As any $W^*$-dynamical system is given by a $u$-continuous homomorphism
$\beta:G\to\Aut(\cM)$,  Remark~\ref{rem:aut-w*}(a) implies that:
\begin{prop} \mlabel{prop:2.9.1}
If $(\cM,G, \beta)$ is a  $W^*$-dynamical system, then the standard form representation
$(\cM, \cH, J,\cC)$ of $\cM$ is covariant for $\beta$. Moreover, the unitary implementers
for $\beta$ can be taken to be in  $\cU(\cH)_\cM$.
\end{prop}

This proposition is what makes  standard forms
particularly useful for physics (cf.~\cite{DJP03, Pi06}).
Note that from Remark~\ref{rem:aut-w*}(b), the implementers in $\cU(\cH)_\cM$ cannot be inner
for nontrivial automorphisms.
This proposition raises the question about how the Arveson spectrum of  $(\cM,\R, \beta)$
is related to the covariant implementers in $\cU(\cH)_\cM$. This will be considered in
Sect.~\ref{SpecCovRep}.

\begin{rem}
\label{SFPf}
As an application of Proposition~\ref{prop:2.9.1}, we
give an alternative proof of the Borchers--Halpern Theorem using
 standard forms (cf. Theorem~\ref{thm:bor}).
If $(\pi,U,\cH)$ is a covariant representation with $F = F(\pi)$, then
the $\alpha^*_G$-invariance of $F$ follows from
\[ \omega_S(\alpha_g(A)) = \tr(\pi(\alpha_g(A))S) = \tr(U_g A U_g^* S)
= \tr(A U_g^* S U_g) = \omega_{U_g^* S U_g}(A) \]
for $S \in \cB_1(\cH)$, $g \in G$ and $A \in \cA$. Lemma~\ref{lem:2.3} shows that $F \subeq (\cA^*)_c$.

Suppose, conversely, that $F$ is $\alpha^*_G$-invariant and contained in $(\cA^*)_c$.
We identify $\cA^*$ with the predual of enveloping $W^*$-algebra $\cA^{**}$ of $\cA$
and write $\alpha^{**}$ for the induced action of $G$ on $\cA^{**}$.
Then $F$ is a $G$-invariant folium of $\cA^{**}$.
Let $Z$ be its central support (cf. Remark~\ref{rem:2.x}(2)). Then
$\cM := Z \cA^{**}$ is a $G$-invariant weakly closed ideal of $\cA^{**}$ with
$\fS_n(\cM) = F$ for which we have a natural morphism of $C^*$-algebras
$\eta \:  \cA \to \cM, A \mapsto ZA$.

Let $\pi \: \cM \to \cB(\cH)$ be a standard form realization of $\cM$ on $\cH$
and observe that the action of $G$ on $\cM$ leads to a unitary representation
$U \: G \to \cU(\cH)$
by Prop.~\ref{prop:2.9.1}. Since $G$ acts continuously on
$\fS_n(\cM) \cong F$, this representation is continuous by \cite[Prop.~3.5]{Haa75}.
We thus have a faithful covariant representation $(\pi, U)$ of $(\cM,G,\alpha^{**}\res_\cM)$,
and by pullback via $\eta$
we obtain a covariant representation of $(\cA,G,\alpha)$ whose folium is
$F(\pi) = \fS_n(\cM) = F$.
\end{rem}

We will need the following lemma.

\begin{lem} {\rm(\cite[Lemma~2.6]{Haa75})} \mlabel{lem:2.7}
Let $(\cM,\cH,J,\cC)$ be a von Neumann algebra in standard form and
$p \in \cM$ a projection. Then $q := p JpJ$ is a projection in $\cB(\cH)$ and
$(q\cM q, q(\cH), qJq, q(\cC))$ is a von Neumann algebra in standard form.
\end{lem}

\begin{rem} \mlabel{rem:3.16}
If $\pi$ is a normal representation of $\cM$, then  $\ker \pi =(\1- p) \cM$ for a central projection $p$
(called the support projection of $\pi$ - cf. Def.~\ref{suppDef} below).
Thus we obtain a direct sum
of $W^*$-algebras
\[ \cM \cong p\cM \oplus (\1 -p)\cM \cong \cN\oplus\ker \pi  , \]
and Lemma~\ref{lem:2.7} asserts that the standard form
$(\cM,\cH,J,\cC)$ decomposes accordingly as a direct sum of standard form realizations
of $\cN$ and $\ker \pi$.
\end{rem}
In the following subsections we will analyze those projections $P\in\cM$
for which $P\cM P$ is standard on $P\cH$.

\subsection{Cyclic projections, reductions and dilations of $W^*$-algebras.}
\mlabel{subsec:3.2}

\begin{defn}
\label{DefRed}
 (i) Given a projection $P$ in a von Neumann algebra
$\cM \subeq \cB(\cH)$, we define the {\it reduced} von Neumann algebra by
 \[
 \cM_P:=P\cM\restriction P\cH\subeq \cB(P\cH),\quad\hbox{and the reduction map}\quad
M\mapsto M_P:=PM\restriction P\cH\,.\]
(ii)
  Given a projection  in the commutant $P\in\cM'$, then we will say that
 $\cM\s P.$ is the von Neumann algebra {\it induced} by $P$.\\[1mm]
 (iii) For a von Neumann algebra $\cM \subeq \cB(\cH),$ a subspace
 $S\subset\cH$ is called {\it $\cM$-generating} if  $\lbr\cM S\rbr = \cH$.
 A vector $\Omega \in \cH$ is {\it $\cM$-cyclic} if $\overline{\cM \Omega} = \cH$.
\end{defn}
Then $\cM_P$ is isomorphic to
 $P\cM P\subeq \cB(\cH)$ by restriction of the latter to $P\cH$. Henceforth
 we will not distinguish between $\cM_P$ and $P\cM P$.
 Note that every \strongH closed hereditary $C^*$-subalgebra of $\cM$ is of the form
 $P\cM P$ for some $P\in \cM$ (cf. \cite[Thm.~4.1.8]{Mu90}).

The next two lemmas recall some basic facts on reduction and induction for von Neumann algebras.

\begin{lem}\mlabel{lem:proj}
Let $\cM \subeq \cB(\cH)$ be a von Neumann algebra and $P \in \cM$ be a projection.
Then the following hold:
\begin{itemize}
\item[\rm(i)]
$(\cM_P)' = P \cM'\res_{P\cH}=(\cM')\s P.\subeq \cB(P\cH),$
\item[\rm(ii)] $Z(\cM)P = Z(\cM_P)$, where $\cM_P = P \cM P$,
\item[\rm(iii)]  $\cM_P$ is  in standard form if and only if $(\cM')\s P.=\cM'\restriction P\cH$ is in standard form.
\end{itemize}
\end{lem}

Proofs of (i) and (ii) are in  \cite[Thm.~3.13]{SZ79}  and \cite[Part~1, Ch.~2, Prop.~1]{Dix82}.
To see (iii), note that (i) implies that $\cM_P$ is standard if and only if $(\cM')\s P.=(\cM_P)'$ is standard, as
a von Neumann algebra is in standard form if and only if its commutant is in standard form (Theorem~\ref{Def-SF}(iii)).

\begin{lem}  \label{CP1cyclic}
Let $\cM \subeq \cB(\cH)$ be a von Neumann algebra and $P \in \cM$ be a projection.
Then the following are equivalent:
\begin{itemize}
\item[\rm(i)] The central support of $P$ i.e.\ $z(P):=\inf\{Z\in Z(\cM)\,\mid\,P\leq Z \}$ is $\1$.
\item[\rm(ii)] $\cH_P := P\cH$ is $\cM$-generating, i.e.\  $\lbr\cM \cH_P\rbr = \cH$.
\item[\rm(iii)] The ideal $\lbr\cM P \cM\rbr$ is weakly dense in $\cM$.
\item[\rm(iv)] The restriction map
$R \: \cM' \to \cM_P', R(M) := M\res_{\cH_P} = MP$
is an isomorphism of $\cM'$ onto~$(\cM_P)'$.
\end{itemize}
If these conditions are satisfied, then we further have:
\begin{equation}
  \label{eq:ker-rel}
 \cH_P = \ker(P\cM(\1-P)).\
\end{equation}
\end{lem}

\begin{prf} The projection $Z$ onto $[\cM \cH_P]$ is contained in
$\cM'$ because $Z\cH$ is $\cM$-invariant. Since $[\cM \cH_P]$ is also
$\cM'$-invariant, we likewise obtain $Z \in \cM'' = \cM$,
so that $Z$ is central in $\cM$.
It coincides with the central support of $P$. Therefore (i) and (ii) are equivalent.

The equivalence of (i) and (iii) follows from the fact that the central support
$Z$ of $P$ has the property that $Z \cM$ is the weakly closed
ideal of $\cM$ generated by $P$, i.e.\  the weak closure of $\cM P \cM$.

That (i) and (iv) are equivalent follows from \cite[Prop.~3.14]{SZ79}
or \cite[Part~1, Ch.~2, Prop.~2]{Dix82}  or \cite[Prop.~2.6.7]{Pe89}.

Now we assume that (i)-(iv) are satisfied.
Since $\cH_P$ is $\cM$-generating,
$(\1-P)\cM\cH_P$ is dense in $\cH_P^\bot = (\1-P)\cH$.
Therefore
\[ \cH_P = ((\cH_P)^{\bot})^\bot
= \lbr(\1-P)\cM)\cH_P\rbr^\bot= \lbr(\1-P)\cM P\cH\rbr^\bot \]
implies that $\cH_P = \ker((\1-P)\cM P)$.
\end{prf}

Below we will say that a projection  $P \in \cM$ is {\it generating}
if the subspace $\cH_P := P\cH$ is $\cM$-generating, i.e.\
$P$ has central support $\1$. This property is also equivalent to the
injectivity of the map $\cM' \to P\cM', M \mapsto MP$; in this sense
$P$ is {\it separating for $\cM'$}.
These will be very important below, e.g.
in Lemma~\ref{lem:mini}.

\begin{rem}  Let $\cM$ be a von Neumann algebra
and $P \in \cM$ be a projection with central support~$\1$.
Then the preceding lemma shows that
$\cM' \cong \cM_P'.$
In general, the complementary projection $\1-P$ need not have central support
$\1$. In fact, there may be a non-zero central projection $Z \leq P$.
Then $Z \cM = Z\cM_P$ is an ideal of $\cM$ contained in $\cM_P$.
If $\cM_P$ contains no proper ideals of $\cM$, then $\1-P$ also has central
support $\1$, so that we obtain
$\cM_P' \cong  \cM' \cong \cM_{\1-P}'.$
Therefore the von Neumann algebras $\cM_P$ acting on $\cH_P$ and the
von Neumann algebra \break $\cM_{\1-P}$ acting on $\cH_P^\bot$
have isomorphic commutants.
This is in particular the case if $\cM$ is a factor.
\end{rem}

\begin{ex} If the projection $P\in\cM$ is minimal with central support $\1$,
then $\cM_P \cong \C$ implies that $Z(\cM) \cong \C$, so that $\cM$
is a factor. Further, the existence of minimal projections implies that $\cM$
is of type I, hence isomorphic to some $B(\cK)$.
\end{ex}

Let $\cM \subeq \cB(\cH)$ be a von Neumann algebra and $P \in \cM$ be a projection.
Then for any normal representation of the reduced algebra $\pi_0:\cM_P\to\cB(\cH_0)$
there is a natural completely positive map $\phi:\cM\to\cB(\cH_0)$ defined by
\[ \phi \: \cM \to \cB(\cH_0), \quad \phi(M) := \pi_0(PMP) \]
which is a normal map. Thus there exists a normal
minimal Stinespring dilation
$(\pi_\phi, \cH_\phi,V_\phi)$, which is unique up to unitary equivalence (cf.~\cite[Thm.~IV.3.6]{Ta02},
\cite[Thm.~III.2.2.4]{Bla06}).
It consists of a normal representation $\pi_\phi$ of $\cM$ on $\cH_\phi$
and a continuous linear map $V_\phi \:  \cH_0 \to \cH_\phi$ with
\begin{equation}
\label{Stineprop}
 \pi_0(PMP) = V_\phi^* \pi_\phi(M) V_\phi \quad \mbox{ for } \quad M \in \cM
\quad\hbox{and}\quad
\lbr\pi_\phi(\cM)V_\phi\cH_0\rbr=\cH_\phi.
\end{equation}

We recall the construction for use below. There are several possible definitions,
which coincide by the uniqueness theorem  (cf.~\cite[Thm.~IV.3.6]{Ta02}).
\begin{defn}
\label{StineDef}
Given a von Neumann algebra $\cM \subeq \cB(\cH)$, a projection $P \in \cM$
and a normal representation of the reduced algebra $\pi_0:\cM_P\to\cB(\cH_0)$,
then the {\it minimal Stinespring dilation}
$(\pi_\phi, \cH_\phi,V_\phi)$ with respect to ~$\phi \: \cM \to \cB(\cH_0),$ $\phi(M) := \pi_0(PMP)$
is constructed as follows.
 Equip the algebraic tensor product
 \[
 \cM\otimes_{\cM_P} \cH_0:=\big(\cM\otimes \cH_0\big)\big/\cJ, \ \ \hbox{where}\ \
 \cJ:={\rm Span}\big\{MB\otimes\xi-M\otimes\pi_0(B)\xi
 \mid M\in\cM,\, B\in\cM_P, \,\xi\in\cH_0 \big\} \]
with a sesquilinear inner product, given on the elementary tensors by
\[
\la M\otimes\xi,N\otimes\eta\ra
:=\la \phi(N^*M)\xi,\eta\ra,\qquad M,\,N\in\cM,\;\xi,\,\eta\in\cH_0.
\]
This is well defined because $\phi(MB) = \phi(M) \pi_0(B)$ for $M \in \cM$ and $B \in \cM_P$.
Then factor out by the kernel $\cN:=\{\psi\in\cM\otimes_{\cM_P}
\cH_0\,\mid\, \la \psi,\psi\ra=0\}$ and complete
to obtain $\cH_\phi$.
Denote the factoring map by $\gamma:\cM\otimes_{\cM_P} \cH_0\to\cH_\phi$, and
define $\pi_\phi:\cM\to\cB(\cH_\phi)$  by
\[
\pi_\phi(A)\gamma(M\otimes\xi):=\gamma( AM\otimes\xi)
\]
and then extending it to $\cH_\phi$. Define
\[ V_\phi:\cH_0\to\cH_\phi, \quad V_\phi\xi:=\gamma(\1\otimes\xi) \]
which is an isometry as $\phi(\1)=\1$, which allows us
to identify $\cH_0$ with the subspace $V_\phi\cH_0$ in $\cH_\phi$.
\end{defn}
 It is easy to verify the
claimed properties of $(\pi_\phi, \cH_\phi,V_\phi)$ in
\eqref{Stineprop} from this construction.
Note that
$M\otimes\xi=MP\otimes\xi$ in $\cM \otimes_{\cM_P} \cH_0$.

For $P = \1$ we have $\phi = \pi_0$, which implies that
$\pi_\phi = \pi_0$ and that $V_\phi = \1$, where
we use the canonical identification of
$\cM \otimes_\cM \cH_0$ with $\cH_0$.

The given definition is a restriction of a more general definition for any
completely positive map $\phi$ (cf.~\cite[proof of Thm.~IV.3.6]{Ta02}). In this form,
 if $\phi$ is a state, then
$\pi_\phi:\cM\to\cB(\cH_\phi)$ is just the GNS-representation of the state.

\begin{defn}
\label{suppDef}
 If $(\pi, \cH)$ is a normal representation of the $W^*$-algebra $\cM$,
then we define the {\it support of $\pi$} as the unique central projection
$s(\pi)$ for which $\ker \pi = (\1 - s(\pi))\cM$ (cf.~\cite[Def.~1.21.14]{Sa71}).
\end{defn}
Then $\pi(\cM)\cong s(\pi)\cM$. (This definition of $s(\pi)$ is consistent
with the definition of it in the context of $\cA^{**}$ in Remark~\ref{rem:2.x}(3),
where it coincides with $z(\pi)$.)
\begin{lem} \mlabel{StineReduc}
Let $\cM \subeq \cB(\cH)$ be a von Neumann algebra and $P \in \cM$ be a projection.
 Fix a  normal representation of the reduced algebra $\pi_0:\cM_P\to\cB(\cH_0)$
and define
  $\phi:\cM\to\cB(\cH_0)$ by $\phi(M) := \pi_0(PMP).$ Then
the representation $(\pi_\phi,\cH_\phi)$ (cf. Def.~\ref{StineDef})
has the following properties:
\begin{itemize}
\item[\rm(i)] $s(\pi_\phi)=z(s(\pi_0))$, where $z(M)\in\cZ(\cM) $ denotes
the central support of $M\in\cM$
and $s(\pi_0) \in \cZ(\cM_P)$ is the central support of $\pi_0$.
\item[\rm(ii)] $V_\phi$ is $\cM_P$-equivariant,
i.e.\ $\pi_\phi(B)V_\phi=V_\phi\pi_0(B)$ for all $B\in\cM_P$. In particular,
$V_\phi\cH_0$ is $\pi_\phi(\cM_P)$-invariant.
\item[\rm(iii)] $V_\phi\cH_0= \pi_\phi(P)\cH_\phi$.
\end{itemize}
\end{lem}

\begin{prf} (i)
We have $A\in\ker(\pi_\phi)$ if and only if for all $M\in\cM$ and $\xi\in\cH_0$ we have
\[
0=\left\|\pi_\phi(A)\gamma(M\otimes\xi)   \right\|^2=\left\|\gamma(AM\otimes\xi)   \right\|^2
=\la \phi(M^*A^*AM)\xi,\xi\ra\,.
\]
As this holds for all $\xi$, it is equivalent to
\[
0=\phi(M^*A^*AM)=\pi_0(PM^*A^*AMP)
\quad\forall\,M\in\cM.
\]
Then the preceding is equivalent to
\[0=s(\pi_0) PM^*A^*AMP s(\pi_0)=(AMs(\pi_0))^*(AMs(\pi_0)),\quad
\mbox{ i.e.\ } \quad AM s(\pi_0) =0.\]
We conclude that
$A \in \ker(\pi_\phi)$ is equivalent to $A\cM s(\pi_0) = \{0\}$,
and this is equivalent to ${Az(s(\pi_0))=0}$.
This proves that $s(\pi_\phi) = z(s(\pi_0))$.

(ii) This follows from
\[ \pi_\phi(A) V_\phi\xi
=  \pi_\phi(A) \gamma(\1 \otimes \xi)
= \gamma(A \otimes \xi)
= \gamma(\1 \otimes \pi_0(A)\xi) = V_\phi\pi_0(A)\xi\]
for $A \in \cM_P$ and $\xi \in \cH_0$.

(iii)  For $M\in\cM$ and $\xi\in\cH_0$, we have
\[\pi_\phi(P)\gamma(M\otimes\xi)
=\gamma(PM\otimes\xi)=\gamma(PMP\otimes\xi)
=\gamma(\1\otimes\pi_0(PMP)\xi)
= V_\phi\pi_0(PMP)\xi.\]
This shows that $\pi_\phi(P)\cH_\phi =V_\phi\cH_0$
because $\pi_0(\cM_P)\cH_0 = \cH_0$.
\end{prf}

\begin{prop}
\mlabel{BijRedRep}
Let $\cM$ be a $W^*$-algebra and let $P\in\cM$ be a projection.
Given a normal representation $(\pi,\cH)$ of $\cM$
in which $\cH_{0} := \pi(P)\cH$ is $\cM$-generating,
construct the restricted
representation $(\pi_0, \cH_{0})$ of the reduction $\cM_{P}=P\cM P\subset\cM$
by $\pi_0(PMP):=\pi(PMP)\restriction\cH_0$, $M\in\cM$.
Then the map $\pi\to\pi_0$ defines a bijection between
unitary equivalence classes of normal
representations $(\pi,\cH)$ of $\cM$ generated by the spaces $\pi(P)\cH$,
and unitary equivalence classes of normal
representations $(\pi_0, \cH_0)$ of the reduction $\cM_{P}$.
\end{prop}

\begin{prf} Since the assignment
$\pi \mapsto \pi_0$ defines a functor from the category of normal
$\cM$-representations in which the range of $P$ is generating to the
category of normal $\cM_P$-representations, it induced a well-defined
map on the level of unitary equivalence classes.

To see surjectivity, let
 $(\pi_0, \cH_0)$ be a normal representation of $\cM_{P}$ and
define $\phi$ as $\phi(M) = \pi_0(P M P),$ $M\in\cM$, then the corresponding minimal dilation
$(\pi_\phi, \cH_\phi)$ is a normal representation of $\cM$
for which $\pi_\phi(P)\cH_\phi$ is generating (cf.~Lemma~\ref{StineReduc}(iii)).
The restriction $(\pi_\phi)_0$ of $\pi_\phi(\cM_P)$ to $\pi_\phi(P)\cH_\phi=V_\phi\cH_0$
is then unitarily equivalent to $(\pi_0, \cH_0)$ by  Lemma~\ref{StineReduc}(ii).

To verify injectivity, we have
to show that $\pi_0 \cong \pi_0'$ implies that $\pi \cong \pi'$.
Since $\cH_0\subset \cH$ is $\pi(\cM)$-generating,
by defining $V:\cH_0\to\cH$ to be the inclusion map, we can
verify the conditions~(\ref{Stineprop}) Thus
 the representation
$(\pi, \cH)$ is equivalent to the minimal Stinespring dilation
$(\pi_\phi, \cH_\phi,V_\phi)$ of the completely positive map
\[ \phi \: \cM \to \cB(\cH_0), \quad
\phi(M) := \pi(P) \pi(M) \pi(P)\res_{\cH_0} = \pi_0(PMP).\]
As the Stinespring construction is functorial from normal $\cM_P$-representations
to normal $\cM$-representations, it maps unitary equivalent  $\cM_P$-representations
to unitary equivalent $\cM$-representations.
\end{prf}

\subsection{Standard projections of $W^*$-algebras.}
\mlabel{subsec:3.3}

We now introduce the following key concept.

\begin{defn}
\label{standP}
Let $\cM \subeq \cB(\cH)$ be a von Neumann algebra.
\begin{itemize}
\item[\rm(i)] We call
a projection $P \in \cM$ {\it standard} if it is generating (i.e.\ its central support is
$\1$, cf. Lemma~\ref{CP1cyclic}), and  $\cM_P:=P\cH$ on $\cH_P$ is standard
(equivalently, the faithful representation of $\cM'$ on $\cH_P$ is standard
(Lemma~\ref{lem:proj}(iii)).
\item[\rm(ii)]
Let $\Omega \in \cH$, then the $\sigma(\cM,\cM_*)\hbox{--closed}$ left ideal $\{ M \in \cM \: M\Omega = 0\}$
can be written as
$\cM (\1-P)$ for a projection $P = s(\Omega)\in\cM$
(cf.~\cite[Prop.~1.10.1]{Sa71}), which we will call  the
{\it carrier projection} of $\Omega$.
This coincides with the carrier projection $s(\omega)$
for the vector state $\omega(M):=\la\Omega,M\Omega\ra$
as in Definition~\ref{def:carrier}.
\end{itemize}
\end{defn}

\begin{exs} The notion of a standard projection depends on the
realization of $\cM$ on some Hilbert space.

(a) For $\cM = \cB(\cH)$, we have $\cM' = \C \1$ and therefore
the rank-one projections are standard.

(b) For the representation of $\cM = \cB(\cK)$ by left multiplications
on the Hilbert space  $\cH := \cB_2(\cK)$, the commutant consists
of $\cB(\cK)^{\rm op}$ acting by right multiplications,
and a projection $P \in \cM$ is standard if and only if
$P = VV^*$ holds for an isometry $V \: \cH \to \cH$, i.e.,
if $P \sim \1$ (see Lemma~\ref{lem:3.2} below).
\end{exs}

 Below we shall see that
the carrier projection of a cyclic vector is standard.
Note that $\1 \in \cM$ is standard if and only if $\cM$ is in standard form.

\begin{lem} A von Neumann algebra $\cM \subeq \cB(\cH)$ contains a standard projection
if and only if there exists an $\cM'$-invariant subspace
$\cH_0 \subeq \cH$ on which the representation of $\cM'$ is faithful and standard.
\end{lem}

\begin{prf} If $P \in \cM$ is standard, then $\cH_P := P\cH$ is generating for $\cM$,
hence separating for $\cM'$. As the representation of $\cM'$ on $\cH_P$ thus leads
to an isomorphism $\cM' \cong (\cM_P)'$ (Lemma~\ref{lem:proj}), the
representation of $\cM'$ on $\cH_P$ is standard.

Suppose, conversely, that $\cH_0$ is a closed subspace of $\cH$ on which the
representation of $\cM'$ is faithful and standard and let $P \in \cM$ be
the orthogonal projection onto $\cH_0$. Then $\cH_0$ is $\cM$-generating because it separates
$\cM'$, and thus $z(P) = \1$. Further, the fact that $P\cM'\res_{\cH_0}$ is
the commutant of $\cM_P$ (Lemma~\ref{lem:proj}) implies that the representation
of $\cM_P$ on $P$ is standard.
\end{prf}

It is instructive to observe that there are von Neumann algebras containing no
standard projections. This happens if the representation is too large.
\begin{exs} (a) We consider the von Neumann algebra
$\cM = \C \1 \subeq \cB(\cH)$. Then $\cM$ contains a standard projection if and only if
$\1$ is standard, and this is equivalent to the representation of $\cM_P = \cM = \C \1$
on $\cH$ being standard. This is only the case for $\dim \cH = 1$.

(b) If $\cM \subeq \cB(\cH)$ is a commutative von Neumann algebra,
then $P = \1$ is the only projection with central support $\1$.
Then $\cH_P = \cH$ and $P$ is standard if and only if the representation of
$\cM$ on $\cH$ is. As $\cM$ is commutative, we then have
$\cM' = J\cM J = J Z(\cM) J = \cM$. In particular, the representation
must be multiplicity free.
For $\cM = L^\infty(X,\fS,\mu)$, where $\mu$ is a finite measure,
this means that the representation of $\cM$ is equivalent to the multiplication
representation on $L^2(X,\fS,\mu)$.

(c) If $\cM \subeq \cB(\cH)$ is a factor of type I, then
$\cH \cong \cK \otimes \cK'$ with $\cM = \cB(\cK) \otimes \1\cong \cB(\cK)$ and
$\cM' = \1 \otimes \cB(\cK') \cong \cB(\cK')$.
Let $P = Q \otimes \1\in \cM$ be a projection. As $\cM$ is a factor,
$z(P) = \1$ whenever $Q \not=0$. Further,
$\cM_P \cong \cB(\cK_Q)$ and $\cH_P = \cK_Q \otimes \cK'$.
The representation of $\cM_P \cong \cB(\cK_Q)$ on this space
is standard if and only if $\cH_P \cong \cB_2(\cK_Q)$
(with the left multiplication representation), and this is equivalent to
$\cK' \cong \cK_Q$. Therefore $\cM$ contains a standard projection
if and only if $\dim \cK' \leq \dim \cK$, i.e., if the multiplicity
space $\cK'$ is isomorphic to a subspace of $\cK$.
\end{exs}

The content of the following lemma can already
be found in St\o{}rmer's approach to modular invariants
of von Neumann algebras in \cite{St72}.

\begin{lem}
\mlabel{lem:3.1}
Let $\cM \subeq \cB(\cH)$ be a von Neumann algebra, let  $\Omega \in \cH$ be a unit vector and  $P = s(\Omega)$ be the corresponding carrier projection.
Then
\begin{itemize}
\item[\rm(i)]
$\lbr \cM'\Omega\rbr = P\cH,$
\item[\rm(ii)] If $\Omega$ is $\cM$-cyclic, then $\Omega \in \cH_P = P\cH$
is cyclic and separating for $\cM_P$. In particular,
$P = s(\Omega)$ is standard.
\end{itemize}
\end{lem}

\begin{prf} (i) Let $Q$ be the projection onto the closed subspace
$\lbr \cM'\Omega\rbr$. As $Q\cH$ is $\cM'$-invariant, the projection $Q$ is
contained in $\cM'' = \cM$.
For $M \in \cM$, the condition $M\Omega = 0$ is equivalent to $M\cM'\Omega= \{0\}$, resp.,
to $MQ = 0$. Therefore $\cM(\1- P) =\cM(\1-Q)$, and this implies that $P = Q$.

(ii)
 First we observe that $\Omega$ is $\cM_P$-cyclic because
$\cH_P = P\cH = \lbr P\cM \Omega\rbr = \lbr P\cM P\Omega \rbr = \lbr\cM_P\Omega\rbr.$
To see that $\Omega$ separates $\cM_P$, let
$M \in \cM_P$ satisfies $M\Omega = 0$, then
the definition of the carrier projection implies that
$M  \in \cM_P \cap \cM(\1-P) = \{0\}$.
From Proposition~\ref{prop:2.9} it now follows
that the representation of $\cM_P$ on $\cH_P$ is standard.
\end{prf}

For the next lemma we use the Murray--von Neumann equivalence
relation $\sim$ recalled in the lines just above Theorem~\ref{thm:cycl}.

\begin{lem}
  \mlabel{lem:3.2}
If $P$ is a standard projection in
the von Neumann algebra $\cM$, then a projection $Q \in \cM$ is standard if and only
if $P \sim Q$.
\end{lem}

\begin{prf} That $P$ is standard means that
the representation of $\cM'$ on $\cH_P = P \cH$ is standard
which implies in particular that $\cM'\cong \cM_P^{\rm op}$.
Since two standard representations of $\cM'$ are equivalent by
Remark~\ref{rem:2.6}(a), it follows from
\cite[Prop.~2.7.3]{Sa71}, applied to $\cA := \cM'$ and
$P,Q \in \cM'' = \cM$, that $Q$ is standard if and only if $P \sim Q$.
\end{prf}

\begin{defn}
\label{PstandRep}
($P$-standard representations)
Let $P$ be a projection in the $W^*$-algebra~$\cM$
and $\rho_P \: \cM_P = P\cM P \to \cB(\cH_0)$ be a faithful standard form representation
of $\cM_P$.
 Then
\[ \phi_P \: \cM \to \cB(\cH_0), \quad \phi_P(M) := \rho_P(PMP) \]
is a normal
completely positive function, so that there exists a normal
minimal Stinespring dilation
$(\pi_{\phi_P}, \cH_{\phi_P},V_{\phi_P})$, which is unique up to unitary equivalence (cf.~\cite[Thm.~IV.3.6]{Ta02},
\cite[Thm.~III.2.2.4]{Bla06}, Proposition~\ref{BijRedRep}).
It is called the {\it $P$-standard representation of $\cM$}.
If there is no risk of confusion, we will omit the subscript $P$ on $\phi_P$ and just use
$(\pi_{\phi}, \cH_{\phi},V_{\phi})$.

It consists of a normal representation $\pi_{\phi_P}$ of $\cM$ on $\cH_{\phi_P}$
and a continuous linear map \break $V_{\phi_P} \:  \cH_0 \to \cH_{\phi_P}$ with
\[ \rho_P(PMP) =  \phi_P(M) = V_{\phi_P}^* \pi_{\phi_P}(M) V_{\phi_P}
\quad \mbox{ for } \quad M \in \cM
\quad\hbox{and}\quad
\lbr\pi_{\phi_P}(\cM)V_{\phi_P}\cH_0\rbr=\cH_{\phi_P}.\]
\end{defn}
The construction and properties of $(\pi_{\phi_P}, \cH_{\phi_P},V_{\phi_P})$
was given above in the previous subsection, but we list the properties again below.

\begin{lem} \mlabel{lem:3.11} For
a projection $P$ in the $W^*$-algebra $\cM$, the Stinespring dilation
$(\pi_{\phi}, \cH_{\phi},V_{\phi})$
for $\phi(M) := \phi_P(M) := \rho_P(PMP)$
has the following properties:
\begin{itemize}
\item[\rm(i)] $s(\pi_\phi)= z(P)$ is the central support of $P$.
\item[\rm(ii)] $V_\phi$ is $\cM_P$-equivariant, i.e.\ $\pi_\phi(B)V_\phi=V_\phi\rho_P(B)$ for all $B\in\cM_P$.
Further, $V_\phi\cH_0$ is $\pi_\phi(\cM_P)$-invariant and the restriction of  $\pi_\phi(\cM_P)$
to this subspace is standard.
\item[\rm(iii)] $V_\phi\cH_0 = \pi_\phi(P)\cH_\phi$.
\item[\rm(iv)] If the central support of $P$ is $\1$, then
$\pi_\phi$ is a faithful normal representation for which the projection $\pi_\phi(P)$
onto $V_\phi\cH_0$ is standard.
\item[\rm(v)] If $\cM \subeq \cB(\cH)$ is a von Neumann algebra and $\cH_P$ is
$\cM$-generating, then the identity
representation of $\cM$ on $\cH$ is unitarily equivalent to $\pi_\phi$
if and only if $P$ is standard.
\end{itemize}
\end{lem}

\begin{prf}
 In Lemma~\ref{StineReduc}, replace $\pi_0$ with $\rho_P$ to obtain the
 $(\pi_\phi, \cH_\phi,V_\phi)$ here.

(i) As $\rho_P$ is faithful, $s(\rho_P) = P$, so that
this follows from Lemma~\ref{StineReduc}(i).

(ii) The equivariance was already proven in Lemma~\ref{StineReduc}(ii). As the restriction
of  $\pi_\phi(\cM_P)$ to $V_\phi\cH_0$ is unitarily equivalent to $\rho_P$ it is clearly standard.

(iii) This is Lemma~\ref{StineReduc}(iii).

(iv) If $z(P)=\1$ then by (i) $\pi_P$ is  faithful.
The rest is clear.

(v) In view of Proposition~\ref{BijRedRep}, the identical representation of
$\cM$ on $\cH$ is equivalent to $\pi_\phi$ if and only if the representation
of $\cM_P$ on $\cH_P$ is equivalent to $(\pi_\phi)_0 \cong \rho_P$, i.e.\
standard by (iii). This means that $P$ is standard.
\end{prf}

\begin{prop}  \mlabel{prop:3.11} For two projections $P,Q$ in the $W^*$-algebra $\cM$,
the representations $\pi_{\phi_P}$
and $\pi_{\phi_Q}$ are unitarily equivalent if and only if $P \sim Q$.
\end{prop}

\begin{prf} (a) Suppose first that $P \sim Q$. Then both
have the same central support.
As $P=z(P)P\in z(P)\cM\cong \cN := \pi_{\phi_P}(\cM)$
has central support $\1$ in $z(P)\cM$, it follows by
Lemma~\ref{lem:3.11}(iv) that $\pi_{\phi_P}(P)$ is a standard projection in
$\cN$. Now Lemma~\ref{lem:3.2} implies
that the projection $\pi_{\phi_P}(Q)$ is also standard in $\cN$
and Lemma~\ref{lem:3.11}(v) implies that the representations
 $\pi_{\phi_P}$ and $\pi_{\phi_Q}$ are unitarily equivalent.

(b) If, conversely, $\pi_{\phi_P} \cong \pi_{\phi_Q}$, then
$\pi_{\phi_P}(Q)$ is a standard projection in $\cN = \pi_{\phi_P}(\cM) \subeq \cB(\cH_{\phi_P})$,
hence equivalent to $\pi_{\phi_P}(P)$ by Lemma~\ref{lem:3.2}.
As $z(P) = s(\pi_{\phi_P})=s(\pi_{\phi_Q}) = z(Q)$, we have
$P, Q \in z(P)\cM \cong \cN$.
As $\pi_{\phi_P}$ is a faithful representation of $\cN$,
it follows that $P \sim Q$ in $z(P)\cM$, and hence that $P \sim Q$ in $\cM$.
\end{prf}

\begin{prop} \mlabel{prop:3.12} For a projection $P$  in the $W^*$-algebra $\cM$,
with central support $\1$, the
representation $(\pi_{\phi_P}, \cH_{\phi_P})$ has a cyclic vector  if and only if
the $W^*$-algebra $\cM_P$ is countably decomposable.
\end{prop}

\begin{prf} If $\cM_P$ is countably decomposable, then its standard representation
contains a cyclic vector $\Omega$ by Remark~\ref{rem:2.6}(d)
and therefore $\Omega$ is $\cM$-cyclic in $\cH_{\phi_P}$.

Suppose, conversely, that $\pi_{\phi_P}$ has a cyclic vector $\Omega$
and that $Q$ is its carrier
projection. Then $\pi_{\phi_P}(Q)$ is a standard projection by Lemma~\ref{lem:3.1}
and $\cM_Q$ is countably decomposable by  \cite[Prop.~III.2.2.27]{Bla06}.
Since the projection $\pi_{\phi_P}(P)$ is also standard,
$\pi_{\phi_P}(P) \sim \pi_{\phi_P}(Q)$ by Lemma~\ref{lem:3.2}, which in turn leads to
$P\sim Q$. We conclude that $\cM_P \cong \cM_Q$ is countably decomposable.
\end{prf}

The following proposition generalizes the observation that a standard form realization
contains all cyclic representations of~$\cM$.

\begin{prop} \mlabel{prop:3.13}
Let $(\cM,\cH,J,\cC)$ be a von Neumann algebra in standard form
and $P \in \cM$ be a projection with central support~$z(P) = \1$.
Then the representation $(\pi_{\phi_P}, \cH_{\phi_P})$ is unitarily
equivalent to the representation of $\cM$ restricted to the range of the
projection $JPJ \in \cM'$.
\end{prop}

\begin{prf}
Consider the projection $P' := JPJ \in \cM'$. It has the same central
support~$z(P') = \1$. This implies that, for the projection $Q := PP' = P'P$,
the map
\[ \Phi \: \cM_P \to \cM_Q := (\cM_P)_{P'} = P'\cM_P, \quad M \mapsto P' M \]
is an isomorphism of von Neumann algebras (cf. \cite[Prop.~2.6.7]{Ped79}).
In fact, since $P'$ is generating  for $\cM'$ because $z(P')=\1$, it is separating for $\cM$.
From Lemma~\ref{lem:2.7} (cf.~\cite[Lemma~2.6]{Haa75}), we know that
$(\cM_Q, Q\cH, QJQ, Q(\cC))$ is a von Neumann algebra in standard form.
Consider the linear map
\[ \gamma \: \cH_P \to \cH_Q = Q \cH = P'P\cH = P'\cH_P, \quad \xi \mapsto P'\xi = Q\xi.\]
For $M \in \cM_P$ we then have  $\gamma(M\xi) = P'M\xi = M P'\xi = \Phi(M)\gamma(\xi)$,
so that $\gamma$ intertwines the representation of $\cM_P$ on $\cH_P$
with the representation of $\cM_Q$ on $\cH_Q$.
This implies that the representation
$\rho_P(M) :=P'M$ of
$\cM_P$ on the subspace $\cH_Q = \cH_P \cap \cH_{P'}$ of $\cH_P$
is a faithful standard form representation of $\cM_P$.
As $z(P) = \1$, the subspace $\cH_{P} = P\cH$ is $\cM$-generating, so that
\[ \cH_Q = P' \cH_P \quad \mbox{ and } \quad
\lbr\cM \cH_Q\rbr = \lbr P' \cM \cH_P\rbr = P'\cH = \cH_{P'}.\]
Therefore $M \mapsto P'M$ defines a faithful representation
of $\cM$ on $\cH_{P'}$ (by \cite[Prop.~2.6.7]{Ped79}) in which the subspace $\cH_Q=P\cH_{P'}$ is
$\cM$-generating and carries
a  faithful standard representation of $\cM_P$. We conclude that  this representation
 is $P$-standard, hence
unitarily equivalent to $(\pi_{\phi_P},\cH_{\phi_P})$ (Lemma~\ref{lem:3.11}(v)).
\end{prf}

\subsection{Implementability for $W^*$-dynamical systems}
\label{ImpW-dynsys}

We reconsider Theorem~\ref{thm:cycl} above and we
give another proof based on standard representations.

\begin{thm}\mlabel{thm:cycl2} {\rm(Equivalence Theorem for cyclic representations)}
For two normal states $\omega, \eta$ of a $W^*$-algebra $\cM$, the corresponding
cyclic representations are equivalent if and only if $s(\omega) \sim s(\eta)$,
i.e.\ their carrier projections are equivalent.
\end{thm}

\begin{prf}
First we use Lemma~\ref{lem:3.1} and Lemma~\ref{lem:3.11}(v) to
see that, for the carrier projections $P := s(\omega)$ and $Q := s(\eta)$,
we have $\pi_\omega \cong\pi_{\phi_P}$ and $\pi_\eta \cong \pi_{\phi_Q}$.
Therefore Proposition~\ref{prop:3.11} implies that
$\pi_\omega \cong \pi_\eta$ is equivalent to $P \sim Q$.
\end{prf}

\begin{rem} \mlabel{rem:3.27} For a $W^*$-dynamical system $(\cM,G,\beta)$, we
obtain a similar picture to that in Subsection~\ref{cyclic-cov}
 if we replace the state $\omega$ by a
 projection $P$ and consider the corresponding $P$-standard representation
$(\pi_{\phi_P}, \cH_{\phi_P})$. A necessary condition for $(\pi_{\phi_P}, \cH_{\phi_P})$ to be covariant with respect to~$\beta$,
is that $\beta_G$ preserves the kernel of $\pi_{\phi_P}$, hence
the central support $z(P)$ of~$P$. If this is the case, then
we may replace $\cM$ by $\cM z(P)$, so that we may assume that
$z(P) =\1$ and that $\pi_{\phi_P}$ is faithful.

Another necessary condition is that $\beta_G$ preserves the equivalence class
$[P]$ of projections (Proposition~\ref{prop:3.11}), hence fixes its central
support $z(P)$. If this is the case, then
$\pi_{\phi_P} \circ \beta_g \cong \pi_{\phi_{\beta_g^{-1}(P)}}$
implies that each automorphism $\beta_g$ can be implemented in $\cH_{\phi_P}$.
This leads to a topological group extension
\[  \hat G_P := \{ (g,U) \in G \times \cU(\cH_{\phi_P}) |\
(\forall M \in \cM)\ U \pi_{\phi_P}(M) U^{-1} = \pi_{\phi_P}(\beta_g(M))\}\]
of $G$ by $N := \cU(\pi_{\phi_P}(\cM)') \cong \cU(\cM_P')$
and the covariance of the representation $\pi_{\phi_P}$ is equivalent to
the splitting of this extension of topological groups.

 This is closely related to the Lie group extensions
constructed in \cite{Ne08} for smooth actions of a Lie group
$G$ on a continuous inverse algebra $\cA$. For a projective $\cA$-right module of the form
$P\cA$, $\hat G_P$ is an extension of an open subgroup
\[ G_{[P]} := \{ g \in G \mid \beta_g(P) \sim P \} \]
of $G$ by the unit group $\cA_P^\times = (P\cA P)^\times$.
In the unitary context, which corresponds to Hilbert-$C^*$-modules,
where $\cA$ is a $C^*$-algebra, one expects extensions by
the unitary group $\cU(\cA_P)$.

For the required smoothness it may be enough that the orbit
of $P \in \cA$ is smooth in $\cA$; which is the case if $P$
is a smoothing operator for a unitary representation of
$G$, i.e., $P \cH \subeq \cH^\infty$
 (cf.\ \cite{NSZ17}).
\end{rem}

\begin{thm}\mlabel{thm:Pinv}
Given a $W^*$-dynamical system $(\cM,G,\beta)$ and a
 projection $P\in\cM$ such that
 $P$ is $\beta_G$-invariant, then
 $\beta$ can be continuously implemented in $(\pi_{\phi_P}, \cH_{\phi_P})$,
 i.e.\ $\pi_{\phi_P}$ is covariant. In particular,
the extension $\hat G_P$ of $G$ splits.
\end{thm}

\begin{prf} If $P$ is $\beta_G$-invariant, then $\beta_G$ preserves the
subalgebra $\cM_P$ and can be continuously implemented in the standard
representation $(\rho_P,\cH_0)$ of $\cM_P$ (cf.~Proposition~\ref{prop:2.9.1}).
Then the corresponding completely positive map
\[ \phi_P \: \cM \to \cB(\cM_P), \quad
M \mapsto \rho_P(PMP) \]
is $\beta_G$-equivariant, and the naturality of the Stinespring dilation
implies that $\beta_G$ can be continuously implemented in $(\pi_{\phi_P},\cH_{\phi_P})$.
Explicitly, fix the unitary implementing group $V:G\to\cU(\cH_0)$,
 $\rho_P(\beta_g(M))=V_g\rho_P(M)V_g^*$ for $M,N\in\cM_0$. Then
\begin{eqnarray*}
\Big\la\pi_\phi(\beta_g(A))\gamma(M\otimes\xi),\gamma(N\otimes\eta)\Big\ra &=&
\big\la\gamma(\beta_g(A)M\otimes\xi),\gamma(N\otimes\eta)\big\ra\\[1mm]
&=&\big\la\phi(N^*\beta_g(A)M)\xi,\eta\big\ra
=\big\la\rho_P(PN^*\beta_g(A)MP)\xi,\eta\big\ra\\[1mm]
&=& \Big\la V_g\rho_P(P\beta_{g^{-1}}(N)^*A\beta_{g^{-1}}(M)P)V_g^*\xi,\eta\Big\ra\\[1mm]
&=&\Big\la\pi_\phi(A)\gamma(\beta_{g^{-1}}(M)\otimes V_g^*\xi),\gamma(\beta_{g^{-1}}(N)\otimes V_g^*\eta)\Big\ra\\[1mm]
&=& \big\la U_g\pi_\phi(A)U^*_g\gamma(M\otimes \xi),\gamma(N\otimes \eta)\big\ra,
\end{eqnarray*}
where
\[ U_g\gamma(M\otimes \xi) := \gamma(\beta_{g}(M)\otimes V_g\xi)\quad
\mbox{implies} \quad \pi_\phi(\beta_g(A))=U_g\pi_P(A)U^*_g.\]
It is obvious that $U_g$ is a unitary group homomorphism, by letting $A=\1$ above,
and weak operator continuity is also easy to see.
\end{prf}

The following example shows that Theorem~\ref{thm:Pinv} does not extend
directly to the case where only $[P]$ is $G$-invariant. This case requires the
passage to possibly non-trivial central exensions.

\begin{ex} For $\cM = \cB(\cH)$ and $\dim \cH > 1$,
we consider a one-dimensional projection $P \in \cM$
and observe that it is standard by Lemma~\ref{lem:3.1}.
Thus the  representation $(\pi_{\phi_P}, \cH_{\phi_P})$ is unitarily equivalent to the
identical representation of $\cM$ on $\cH$ by Lemma~\ref{lem:3.11}(v). We
consider the action of $G := \PU(\cH)$ on $\cM$ induced by conjugation.
This action leaves the class $[P]$ of the projection $P$ invariant,
but to implement it on $\cH$, we have to pass to the non-trivial central extension
$\hat G = \cU(\cH)$ of $G$ by $\T \cong \cU(\cM_P)$.
That this central extension is non-trivial follows
for infinite dimensional Hilbert spaces from the fact that
every unitary operator is a commutator (\cite[Prob.~239]{Ha82}), and for
$\cH = \C^n$, the subgroup $\T\1\cap\SU_n(\C) \cong C_n \1$ (cyclic
group of order $n$) consists of commutators in $\SU_n(\C)$.
\end{ex}

\begin{rem} (i) If $P =\1$, then $\pi_{\phi_P}$ is the standard representation
of $\cM$ and Theorem~\ref{thm:Pinv} implies that $\Aut(\cM)$
can be implemented (which is already known from Proposition~\ref{prop:2.9.1}).

(ii) If $\cM$ is a von Neumann algebra, and the
$G$-invariant projection $P$ is standard,
then the covariant representation $\pi_{\phi_P}$ is faithful
(cf. Lemma~\ref{lem:3.11}(v)) and unitarily equivalent to
the identity representation of $\cM$. Hence the identity representation
of $\cM$ is covariant.

(iii) \cite[Thm.~8]{Hal72} describes criteria for the implementability in terms
of the $G$-action on $\cZ(\cM)$ and
\cite[Cor.~10]{Hal72} concerns semi-finite von Neumann algebras.

(iv) \cite[III.2.6.15/16]{Bla06} has a criterion for a von Neumann algebra
$\cM\subeq \cB(\cH)$ to be in standard form: If $\cH$ is separable and $\cM'$ is properly
infinite. In view of (i), this can be viewed as a sufficient condition for
unitary implementability of the $G$-action.
\end{rem}

\begin{rem} (Equivalence classes of projections for factors) \\
(a) If $\cM = \cB(\cH)$ is a factor of type I, then two projections
$P,Q \in \cM$ are equivalent if and only if $\dim P\cH = \dim Q\cH$, i.e.\  the
set of equivalence classes is parameterized by the Hilbert dimensions
of closed subspaces of $\cH$, which is the set of all cardinals $\leq \dim_{\rm Hilb} \cH$.

In this case $\Aut(\cM) = \PU(\cH)$ acts by conjugation, so that every class
$[P]$ is invariant under $\Aut(\cM)$.

(b) If $\cM$ is a factor of type II$_1$, then the set of equivalence classes
of finite projections (this means that $P \sim Q \leq P$ implies $P = Q$)
can be identified with the unit interval $[0,1]$ because any
normalized trace $\tau \:  \cM\to \C$ provides a complete invariant.
Since $\tau$ is $\Aut(\cM)$-invariant, the automorphism group also preserves
all equivalence classes of projections.

(c) From \cite[Thm.~III.1.7.9]{Bla06} we recall that, the set $[\Proj(\cM)]$
for a countably decomposable factor can be described as:
\begin{itemize}
\item $\{0, 1, \ldots, n\}$ if $\cM$ is of type I$_n$, $n\in \N$,
\item $\{0, 1, 2, \ldots, \infty\}$ if $\cM$ is of type I$_\infty$,
\item $[0, 1]$ if $\cM$ is of Type II$_1$.
\item $[0, \infty]$ if $\cM$ is of Type II$_\infty$.
\item $\{0,\infty\}$ if $\cM$ is of Type III.
\end{itemize}

This shows that only for type II$_\infty$, there is no a priori reason for
$\Aut(\cM)$ to preserve all equivalence classes of projections.
Let $\cM$ be a  factor of type II$_\infty$.
Let $\cM_+$ be its cone of positive elements
and assume that $\tau_0\colon\cM_+\to[0,\infty]$ is a semi-finite faithful normal trace.
Then, for any $P,Q\in\Proj(\cM)$ with $\min\{\tau_0(P),\tau_0(Q)\}<\infty$,
we have $\tau_0(P)=\tau_0(Q)$ if and only if $P\sim Q$
(see \cite[Part III, Ch.~2, \S~7, Prop.~13(iii)]{Dix82}).
The trace $\tau_0$ is unique up to multiplication by a positive scalar
by \cite[Part I, Ch. 6, \S 4, Cor.]{Dix82}, hence
there exists a group homomorphism $\mu\colon\Aut(\cM)\to \R^\times_+$
depending on $\tau_0$,
with $\tau_0\circ\theta=\mu(\theta)\tau_0$
for every $\theta\in\Aut(\cM)$.
Thus, if $\theta_0\in\Aut(\cM)$ satisfies $\mu(\theta_0)\ne 1$, then for every $P\in\Proj(\cM)\setminus\{0\}$ with $\tau_0(P)<\infty$
we have $\tau_0(\theta_0(P))\ne \tau_0(P)$, hence $\theta_0(P)\not\sim P$.
Specific examples of such automorphisms of factors of type II$_\infty$
occur in connection with the structure of factors of type III$_1$;
see \cite[Ch. XII, Th.~1.1(ii) and Def. 1.5(iii)]{Ta03}.
In particular,
the hyperfinite factor $\cR_{0,1}$ of type II$_\infty$ admits automorphisms $\theta_0$ as above,
because $\cR_{0,1}$ is involved in the decomposition of the hyperfinite factor of type III$_1$ as the crossed product of a $W^*$-dynamical system $(\cR_{0,1},\R,\alpha)$.

It is easy to construct a concrete example of such automorphisms.
We consider the hyperfinite type II$_1$-factor
$\cN = \oline{\bigotimes}_{n \in \N} M_2(\C)$. For the infinite dimensional
Hilbert space $\cH$, the tensor product
$\cM := \cB(\cH) \bar\otimes \cN$ is then a factor of type II$_\infty$.
From any unitary operator $U \: \cH \to \cH \oplus \cH$, we obtain
an isomorphism
\[ \Phi_0 \: \cB(\cH) \to \cB(\cH \oplus \cH) \cong \cB(\cH) \otimes M_2(\C),
\qquad \Phi_0(A) = U A U^{-1}.\]
Now
\[ \Phi \: \cM \to \cM, \qquad \Phi(A \otimes B) := \Phi_0(A) \otimes B, \qquad
B \in \cN  \]
is an automorphism of $\cM$.
On $\cM$ we consider the tensor product trace
$\tau = \tr \otimes \tau_\cN$, where $\tau_\cN$ is the normalized
trace on $\cN$. For a minimal projection $P$ on $\cH$, we have
\[ \tau(\Phi(P \otimes \1)) = \tau(U P U^{-1} \otimes \1)
= (\tr \otimes \tau_{M_2(\C)})(UPU^{-1})
= \frac{1}{2} \tr(UPU^{-1})
= \frac{1}{2} \tr(P)  = \frac{1}{2} \tau(P \otimes \1).\]
This means that $\mu(\Phi) = \frac{1}{2}$.
\end{rem}

\section{Spectral theory for covariant representations}
\label{SpecCovRep}

In this section we will assume that $G=\R$ for simplicity, i.e.\ we have the one-parameter case.
The Arveson spectrum is defined for any locally compact abelian group.

\subsection{Arveson spectrum and spectral conditions for a $C^*$-action $(\cA, \R, \alpha)$.}

\begin{defn}
\label{DefSpecCon}
For a covariant representation $(\pi,U)$ of a $C^*$-action $(\cA, \R, \alpha)$ on $\cH$
we have $U_t=\exp(-itH),$ $t\in\R,$ for some selfadjoint operator $H$ on $\cH$.
In this case,
for  a subset $C \subeq \R$, a
{\it $C$-spectral condition} will mean that
the spectrum $\Spec(H)$ is contained in $C$.
We will mostly be interested in the case that $C=[0,\infty),$
i.e.\ $H\geq 0,$ in which case we will   say that $U:\R\to\cU(\cH)$  has {\it positive spectrum.}                
A covariant representation $(\pi,U)\in{\rm Rep}(\alpha,\cH)$ will be said to have {\it positive spectrum}
if $U:\R\to\cU(\cH)$ has positive spectrum.  
\end{defn}

\cite{Bo84} seems to be the first paper where the spectrum
condition is studied in a context where $\alpha$ is not \strong
continuous.
Note that by adding a real multiple of the identity to $H$ we can trivially convert
a unitary one-parameter group with positive spectrum to one satisfying a
$[\lambda,\infty)$-spectral condition, for any $\lambda\in\R$. So the important property
here is that $H$ is bounded below. However, by the next Proposition, this property need not hold
for all implementing unitary groups.

\begin{prop} Let $(U_t)_{t \in \R}$ be a  \strongH continuous unitary
one-parameter group  with positive spectrum in the von Neumann algebra $\cM\subeq \cB(\cH)$.
Then  $\cM'$ is finite dimensional if and only if
 for any \strongH continuous
unitary one-parameter group $(W_t)_{t \in \R}\subset\cM'$ the spectrum of the
one-parameter group $(U_t W_t)_{t \in \R}$ is also bounded from below.
\end{prop}

\begin{prf}
It is clear that if  $\cM'$ is finite dimensional, then the right hand side follows.
We prove the converse.

(a)  We first deal with the special case where $(U_t)_{t \in \R}\subset\cM$ is norm continuous.
Thus $U_t=\exp(-itH)$ where $H\in\cM$ and $H\geq 0$. Let $(W_t)_{t \in \R}\subset\cM'$
be a \strongH continuous unitary one-parameter group, hence $W_t=\exp(-itB)$ for $B$ a selfadjoint
operator, possibly unbounded. Then $U_t W_t=\exp(-it(H+B)),$ and the assumption is that
$\Spec(H+B)$ is bounded from below. If  $E$ is the spectral measure of $B$, then  the
subspaces $E[n,n+1)\cH$, $n\in\Z$ are all preserved by $H$ and $B$, and
$H+B$ restricted to such a subspace has
spectrum in $[n, n+1+\|H\|]$. Thus if $\Spec(H+B)$ is bounded from below, then
there is a $K$ such that $E[n,n+1)=0$ for $n<K$. Hence
 $\Spec(B)$ is bounded from below. Thus
  the spectrum of every \strongH continuous one-parameter group
$(W_t)_{t \in \R}\subset\cM'$ is bounded from below.
Since this also applies to
$(W_{-t})_{t \in\R}$, it follows that $(W_t)_{t \in \R}$ is norm continuous.

If all \strongH continuous unitary
one-parameter groups in $\cM'$ are norm continuous, then every orthogonal family of
projections in $\cM'$ must be finite (or else from an infinite sequence of projections in $\cM'$
we can define an unbounded selfadjoint operator which generates a one-parameter unitary group in
$\cM'$ which is not norm continuous).
Thus $\cM'$ is finite dimensional by \cite{Og54}
(see also Lemma~\ref{lem:fin}).

(b) Now we turn to the general case.
For $a < b$, let $P[a,b)$ denote the corresponding spectral projection of~$U$. Then
the subspace $\cH[a,b) := P[a,b)\cH$ is invariant under $\cM'$ and
$U$, and since the restriction of $U$ to $\cH[a,b)$ is norm continuous,
 (a) implies that the subalgebra $\cM'[a,b) := P[a,b) \cM'$ of $\cM'$ is finite dimensional.
Let $Z_j \in \cM'$ be the central support of $\cM'[0,j)$, $j \in \N_0$.
If the set $\{ Z_j \: j \in \N_0\}$ is infinite, then there exists a
subsequence $(Z_{j_k})_{k \in \N}$ for which
$Q_k := Z_{j_{k+1}} - Z_{j_k} \not=0$. Then
$B := \sum_{k = 1}^\infty j_{k+1}^2 Q_k$ has the property that $H - B$ is not
bounded from below. Hence there are only finitely many $Z_j$. In particular,
there is a maximal one $Z_N$ which must be $\1$. Therefore the representation
of $\cM'$ on $\cH[0,N)$ is faithful, and this implies that $\cM'$ is finite
dimensional.
\end{prf}

Thus in general, for an action $\alpha:\R\to \Aut(\cM),$
given one implementing unitary group $(U_t)_{t \in \R}$
 with positive spectrum, then
other implementing unitary groups  need not have generators bounded from below,
except if $\cM'$ is finite dimensional.

We will follow the convention of \cite{BR02} that a unitary one-parameter
group $(U_t)_{t \in \R}$ is related to its spectral measure $E$ by
\[ U_t = e^{-itH}= \int_\R e^{-itp}\, dE(p)\quad\hbox{where}\quad \quad H = \int_\R p\, dE(p).\]
In this picture, for $f\in L^1(\R)$ we have
\begin{equation}
  \label{eq:ftra}
 U_f
= \int_\R f(t)U_t\, dt
= \int_\R \int_\R e^{-itp} f(t)\, dt\, dE(p)
= \int_\R \hat f(p)\, dE(p) = \hat f(H).
\end{equation}
Thus if $H\geq 0$ then $U_f=0$ whenever $\supp\hat{f}\subset(-\infty,0)$.\\

Given a covariant representation $(\pi,U)$, there are two spectral theories which we will use;-  that of
$U$ (i.e.\ of $H$), and the Arveson spectral theory for $\alpha$ (cf.~\cite{Arv74}).
 The relation between them will be made explicit.
Arveson's spectral theory was motivated by the search for a constructive proof
of Borchers' theorem (cf.~Theorem~\ref{BA-thm} below; see \cite[Ch.~XI]{Ta03}).
We first define
the Arveson spectral subspaces $\cM^\alpha(S)$
(cf.~\cite[Def~3.2.37]{BR02}- this can be done for any locally compact abelian group):

\begin{defn} \mlabel{def:arveson}
Let $(\cM,\R, \alpha)$ be a $W^*$-dynamical system on a von Neumann algebra
$\cM\subseteq\cB(\cH)$. For $f \in L^1(\R)$, we write
\[ \alpha_f(A):=\int f(t)\alpha_t(A)\, dt, \quad A \in \cM \]
for the corresponding integrated representation (\cite[Lemma~7.5.1]{Pe89}), where
$\alpha_f(A)$ is a weak integral with respect to  the weak operator topology. We define
\begin{itemize}
\item[\rm(1)] the {\it spectrum} of an $A\in\cM$ with respect to  $\alpha$ as
\[
\Spec_\alpha(A):=\big\{p\in \R\mid (\forall f \in L^1(\R))\,
\alpha_f(A) = 0 \ \ \Rarrow \ \ \hat{f}(p)=0\big\},\]
where $\hat f(p) = \int_\R e^{-ixp} f(x)\, dx$ is the Fourier transform.
Then the {\it Arveson spectrum of }$\alpha$, denoted ${\rm Spec}(\alpha)$, is the closure of the union of the sets
$\Spec_\alpha(A)$ for all $A\in\cM$.
(This agrees with the generalization to arbitrary locally compact groups in (\ref{Spec1def}) above.
Useful equivalent definitions are listed in \cite[Prop.~3.2.40]{BR02}).
\item[\rm(2)] For a subset $S\subseteq\R$,
the {\it Arveson spectral subspace} of $\alpha$ is
\[ \cM^\alpha(S):=\oline{\{A\in\cM\,\mid\,\Spec_\alpha(A)\subseteq S\}}^{\sigma},\]
where the closure is with respect to  the $\sigma(\cM,\cM_*)$-topology. The subspace
\[
\cM^\alpha_0(S):=\oline{\Spann\big\{\alpha_f(A)\;\mid\;A\in\cM,\;f\in L^1(\R)\quad\hbox{such that}\quad\supp(\hat{f})\subseteq S\big\}}^{\sigma}\]
is contained in $\cM^\alpha(S)$ and,
if $S$ is open, then $\cM^\alpha(S)=\cM_0^\alpha(S)$ (cf.~\cite[Lemma~3.2.39(4)]{BR02}).
\end{itemize}
\end{defn}
By the definition of the Arveson spectrum ${\rm Spec}(\alpha)$, and the fact that
$\Spec_\alpha(A^*)=-\Spec_\alpha(A)$ \cite[Prop.~3.2.42(1)]{BR02}, it follows that
${\rm Spec}(\alpha)$ is a symmetrical set.

 The basic algebraic structure of the Arveson spectral spaces for $(\cM,\R, \alpha)$ which we will need is:
 \begin{itemize}
 \item[\rm(1)] $\cM^\alpha(S)^*=\cM^\alpha(-S)$ for all subsets $S\subseteq\R$
(cf.~\cite[Lemma~3.2.42(2)]{BR02}),
  \item[\rm(2)] $\cM^\alpha(S_1)\cM^\alpha(S_2)\subseteq \cM^\alpha(\overline{S_1+S_2})$ for all closed subsets
  $S_1, S_2\subseteq\R$ (cf.~\cite[Lemma~3.2.42(4)]{BR02}).
  \item[\rm(3)] The union of the spaces $\cM^\alpha[t,\infty)$ for $t\in\R$ is weak operator dense in $\cM$
  (cf.~Lemma~\ref{Hinfty}(1) below).
 \end{itemize}
The space $\cM^\alpha(\{0\})=\cM^{\R}$ is the von Neumann algebra of
 invariant elements, and if  $U:\R\to\cU(\cH)$ is a strong operator continuous unitary implementing group for $\alpha$, then
clearly $\cM^\alpha(\{0\})=U_{\R}'\cap\cM$. If $U_{\R}\subset\cM$ then $U_{\R}''\subset\cM^\alpha(\{0\})$.

The Arveson spectral spaces determine uniquely the action $\alpha:\R\to{\rm Aut}(\cM)$
by the following  (cf.~\cite[Prop.~3.2.44]{BR02}):
\begin{prop}\mlabel{ArvSpecAct}
Let $(\cM,\R, \alpha)$ and $(\cM,\R, \beta)$ be two $W^*$-dynamical systems on a von Neumann algebra
$\cM\subseteq\cB(\cH)$ such that
\[
\cM^\alpha[t,\infty)\subseteq\cM^\beta[t,\infty) \quad\mbox{ for } \quad t\in\R.
\]
Then $\alpha_t=\beta_t$ for all $t\in\R$.
\end{prop}

One can obtain the Arveson spectral spaces from the spectral projections
$E[t,\infty)$ of a unitary group implementing $\alpha$ by
\begin{equation}
\label{MalphChar}
\cM^\alpha[t,\infty)=\big\{A\in\cM\;\mid (\forall s \in \R)\
\;A\,E[s,\infty)\cH\subseteq E[s+t,\infty)\cH \big\}
\end{equation}
(cf.~\cite[Lemma~3.2.39(3), Prop.~3.2.43]{BR02}). Such an implementing unitary group will exist
if we choose e.g.\ $\cM=\cM_{co}$ as above for a given action  $(\cA, \R, \alpha)$.
This suggests that $\cM^\alpha[t,\infty)$ consists of ``shift operators,'' and indeed, we can write
$\cM$ in terms of ``matrix'' expansions w.r.t $E$ (or equivalently
$U(C^*(\R)) = \oline{U_{L^1(\R)}}$),
and characterize the Arveson spectral subspaces $\cM^\alpha(S)$
in these terms:

 \begin{ex}
 \label{Zmatrix}
In the case that the generator $H$ of $U$ has spectrum only in $\Z$,
\eqref{MalphChar} above shows that
 with respect to  the matrix decomposition of $A$ with respect to
the eigenspaces of $H$, an $A\in\cM^\alpha[t,\infty)$ must consist of an upper triangular (infinite) matrix, cf.~\cite[Rem.~C.4]{GrN14}.

Specifically, let $\alpha:\R\to\Aut\cB(\cH)$ be the conjugation $\alpha_t(A) = U_t A U_{-t}$,
where $U_{2\pi} = \1$, so that it actually defines a representation
of the circle group $\T \cong \R/2\pi \Z$.
Denote by
\[ \cB(\cH)_n := \{ A \in \cB(\cH)\mid (\forall t \in \R)\
\alpha_t(A)  = e^{int} A\}\]
its eigenspaces in $\cB(\cH)$ and similarly let
$\cH_n$ be the eigenspace of $U$ in $\cH$ with the projection $P_n$ onto it.
Note that $\cB(\cH)_n=\cB(\cH)^\alpha\{n\}$, i.e.\ it coincides with the Arveson
spectral subspace for $\{n\}$.
The Peter--Weyl Theorem generalizes to continuous Banach representations of $G$
(cf.~\cite[Thm.~2]{Sh55} and~\cite[Thm.~3.51]{HM13}), hence an application of it
to $\alpha\restriction\cB(\cH)_c$ implies that
\begin{equation}
\label{BHc}
 \cB(\cH)_c = \oline{\Spann\Big(\mathop{\bigcup}_{n \in \Z} \cB(\cH)_n\Big)}.
 \end{equation}
Write $A = (A_{jk})_{j,k\in \Z}$ as a matrix with
$A_{jk} \in \cB(\cH_k, \cH_j)$, and keep in mind that the convergence
$A=\sum\limits_{j\in\Z}\sum\limits_{k\in\Z}A_{jk}=\sum\limits_{j\in\Z}\sum\limits_{k\in\Z}P_jAP_k$ is in general
with respect to  the strong operator topology.
We have
\[ \alpha_t(A) = (e^{it(j-k)}A_{jk})_{j,k\in \Z},\]
so that
\[ A \in \cB(\cH)_n \quad \Longleftrightarrow \quad
(j-k \not=n \Rarrow A_{jk} =0).\]
For $A = (A_{jk})_{j,k\in \Z} \in \cB(\cH)$, let
$A_n := (A_{jk} \delta_{j-k,n})_{j,k\in \Z}$ and observe that $A_n$ defines a bounded
operator on $\cH$, hence an element of $\cB(\cH)_n$. In fact, all elements of the
Arveson spectral space
$\cB(\cH)^\alpha\{n\}$ must be of this type, i.e.\ consist of a single diagonal
in the $n$th position above the main diagonal.
As $\cB(\cH)^\alpha[t,\infty)$ is the strong operator closed span of all $\cB(\cH)^\alpha\{n\}$ with $n\geq t$,
we see that the matrix decomposition of an $A\in\cB(\cH)^\alpha[t,\infty)$ consists of upper triangular matrices for which the
$n$th diagonal is zero if $n<t$.

Consider the invariance subalgebra $\cB(\cH)_0=\cB(\cH)^\alpha\{0\}$, which we note from the matrix decomposition must
consist of elements of the form $A=\sum\limits_{k\in\Z}A_{kk}=\sum\limits_{k\in\Z}P_kAP_k$ (strong operator convergence).
We may therefore define a projection $p_0:\cB(\cH)\to\cB(\cH)_0$ onto the invariant algebra by
\[
p_0(A):=\sum\limits_{k\in\Z}P_kAP_k\in \cB(\cH)_0\quad\hbox{for}\quad A\in\cB(\cH).
\]
As the maps $A\to P_kAP_k$ are completely positive, it is clear that $p_0$ is a strong operator limit of completely positive maps
(the finite partial sums) hence it is completely positive. It coincides with the usual group-averaging projection onto
$\cB(\cH)_0$ by:
\begin{eqnarray*}
\int_\T \alpha_z(M)\, dz&=&p_0\Big(\int_\T \alpha_z(M)\, dz    \Big)=\sum\limits_{k\in\Z}P_k\int_\T \alpha_z(M)\, dzP_k\\[1mm]
&=&\sum\limits_{k\in\Z}\int_\T \alpha_z(P_kMP_k)\, dz=\sum\limits_{k\in\Z}P_kMP_k=p_0(M).
\end{eqnarray*}
 \end{ex}

In this example, we obtained a completely positive projection $p_0:\cB(\cH)\to\cB(\cH)_0$. By applying the
Stinespring Dilation Theorem (or more precisely the generalized GNS construction in its proof), any representation $(\rho_0,\cK_0)$ of
$\cB(\cH)_0$ leads to a new representation $(\rho, \cK)$ of $\cB(\cH)$ with
$\cK_0 \subeq \cK$ for which $\rho_0(p_0(A)) = P^*\rho(A) P$ holds for the
orthogonal projection $P \: \cK \to \cK_0$
(cf.~\cite[Thm~IV.3.6]{Ta02}).
The question now arises whether we have such a map $p_0$ in the general case.
In fact we do by the following (cf.~\cite[Lemma~1.4]{EW74}):

\begin{prop}\mlabel{average}
Let $(\cM,\R, \alpha)$ be $W^*$-dynamical system for a von Neumann algebra $\cM\subseteq\cB(\cH)$ and let
$\eta$ be an invariant mean on $C_b(\R)$. For each $M\in\cM$ define $\hat\eta M\in\cM=(\cM_*)^*$ by
\[
(\hat\eta M)(\varphi)
:=\eta\big(\varphi(\alpha^M)\big) \quad\hbox{for all}\quad
\varphi\in\cM_* \quad \mbox{ and } \quad
\alpha^M(t) := \alpha_t(M).
\]
Then the map $\hat\eta:\cM\to\cM$ is an $\alpha_\R\hbox{--invariant}$ conditional expectation onto
the fixed point algebra $\cM^\alpha(\{0\})=\cM^\R$.
 \end{prop}

As conditional expectations are completely positive, it follows that the maps
 $\hat\eta:\cM\to\cM^\R$ are always completely positive (cf.~\cite{NTU60}).
 Under specific additional assumptions, the maps $\hat\eta$ can even be independent of the
 choice of $\eta$ (cf.~\cite{EW74}). Moreover, if the
completely positive map is normal,
 then there is a normal version of the Stinespring
Theorem which guarantees that the new representation
must be normal (cf.~\cite[Thm~III.2.2.4]{Bla06}).
We note however that there may exist no invariant mean $\eta$ on $C_b(\R)$
for which the map $\hat\eta:\cM\to\cM$ from Proposition~\ref{average} is normal,
as the following example shows:

\begin{ex} \mlabel{ex:4.6}
Let $\cH=L^2(\R)$, $\cM=\cB(\cH)$, and for every $f\in L^\infty(\R)$ let $M_f\in\cM$ be the operator
defined by multiplication by~$f$.
Also, for every $t\in\R$, let  $\chi_t\in L^\infty(\R)$ be given by $\chi_t(x):=e^{itx}$ for all $x\in\R$.
Defining $\alpha_t(A):=M_{\chi_t}AM_{\chi_t}^*$ for all $A\in\cM$ and $t\in\R$,
we claim that
$(\cM,\R, \alpha)$ is a $W^*$-dynamical system with the property that, for every
invariant mean $\eta$ on $C_b(\R)$,
the conditional expectation $\hat\eta$ fails to be normal.
In fact, as Proposition~\ref{average} shows that $\hat\eta$ is a conditional expectation onto
$\cM^\alpha(\{0\})$, it suffices to check that there exists no normal conditional expectation
from $\cB(\cH)$ onto~$\cM^\alpha(\{0\})$.
To this end, first note that $\cM^\alpha(\{0\})=\{M_{\chi_t}\mid t\in\R\}'$.
As the $\sigma$-algebra of Borel subsets of $\R$ is the smallest one for
which all functions $\chi_t$ are measurable, they generate
the von Neumann algebra $\cD := L^\infty(\R)$ by
Corollary~\ref{cor:genvonneumann}.
As $L^\infty(\R)$ is a maximal abelian self-adjoint subalgebra of $\cM$ (see for instance \cite[Part I, Ch. 7, no. 3, Th. 2]{Dix82}),
it follows that $\cM^\alpha(\{0\}) = \cD' = \cD$.

On the other hand, for every conditional expectation $E\colon \cB(\cH)\to\cD$
one has $\cK(\cH)\subseteq\ker E$ by \cite[Rem. 5]{KS59}, hence $E$ cannot be $\sigma$-weakly continuous,
because $\cK(\cH)$ is $\sigma$-weakly dense in $\cB(\cH)$.
This shows that our claim above holds true.
\end{ex}

\begin{rem}
(i) This example can be easily generalized to $\cH=L^2(G)$ for any non-discrete locally compact abelian group $G$ instead of $\R$, using the same averaging procedure
(see also \cite{BP07}).
If $G$ is a discrete abelian group, its dual $\widehat{G}$ is a compact abelian group and one has a normal conditional expectation from $\cB(\cH)$ onto its maximal abelian subalgebra consisting of the multiplication operators by functions in $L^\infty(G)=\ell^\infty(G)$, just as in the special case discussed in Example~\ref{Zmatrix}, where $G=\Z$ and $\widehat{G}=\T$.
\\[2mm]
(ii)
By the Kovacs \& Sz\"ucs Theorem (cf.~\cite[Prop.~4.3.8, p.~383]{BR02}),
the statement of Proposition~\ref{average} can be strengthened to give
a normal invariant conditional expectation. For this, we need to assume in addition,  that the subspace
of invariant vectors is $\cM'$-generating, and that the given representation $\cM\subseteq\cB(\cH)$
is covariant for $\alpha$.
\end{rem}
In the case that we have a representation in standard form, the connection between
the Arveson spectrum of the $W^*$-dynamical system $(\cM,\R, \beta)$ and the
one-parameter group of unitary implementers (cf. Proposition~\ref{prop:2.9.1})
is more direct:
\begin{prop}
\label{SpectraNormalF}
For any  $W^*$-dynamical system $(\cM,\R, \beta)$ such that $\cM$ has an invariant
faithful normal weight, then in the standard form representation of $\cM$,
the Arveson spectrum ${\rm Spec}(\beta)$ coincides with the spectrum of
the one-parameter group $U:\R\to \cU(\cH)_\cM$ which implements $\beta$.
\end{prop}

Recall Proposition~\ref{prop:2.9.1} which follows from the fact that
$\Aut(\cM) \cong \cU(\cH)_\cM$ in any standard form realization of $\cM$.
The proof of Proposition~\ref{SpectraNormalF}
 is in \cite[Prop.~XI.1.24]{Ta03}. Note that by uniqueness of the standard
form, the existence of an invariant faithful weight (or an invariant faithful normal state)
is enough.
\begin{rem}
As a selfadjoint operator $A$ on a Hilbert space $\cH$ has a division of its spectrum
\[
\Spec(A)=\overline{\Spec_{pp}(A)}\cup\Spec_{ac}(A)\cup\Spec_{sing}(A)\qquad
\hbox{with  decomposition}\qquad
\cH=\cH_{pp}\oplus\cH_{ac}\oplus\cH_{sing}
\]
one may look for a similar decomposition of the Arveson spectrum
of a $C^*$-action, and to relate this
to the decomposition of the spectrum of its implementing groups. This has indeed
been done for the \usual case with additional assumptions (cf.~\cite{Dy10}),
but thus far not for our case.
\end{rem}

 We also have:

\begin{lem} \mlabel{lem:4.11}
Let $(\cM,\cH,J,\cC)$ be a standard form realization of a von Neumann algebra,
and $(\beta_t)_{t \in \R}$ be a \strongH continuous one-parameter group
of $\cU(\cH)_\cM \cong \Aut(\cM)$. Then the following assertions hold:
\begin{itemize}
\item[\rm(i)] If $\beta$ is implementable on $\cM$ by a unitary one-parameter
group $(U_t)_{t \in\R}$ in $\cU(\cM)$
and $V_t := J U_t J$ is the corresponding one-parameter
group of $\cU(\cM')$, then $\beta_t = U_t V_t$ for all $t \in \R$.
\item[\rm(ii)] If $\Spec(U) \subeq [0,\infty)$, then
$\Spec(V) \subeq (-\infty,0]$ and the factorization of $\beta$ corresponds to the
factorization into the negative and positive spectral part.
\item[\rm(iii)] If $\beta_t = e^{-itH}$, then $JHJ = - H$.
In particular, the spectrum of $H$ is symmetric.
\end{itemize}
\end{lem}

\begin{prf} (i) Since $U_t$ implements the conjugation with $\beta_t$,
both commute for every $t$. The same holds for $V_t$ because
\[ \beta_t M \beta_t^{-1} = V_t M V_t^{-1} \quad \mbox{ for } \quad t \in \R, M \in \cM'\]
follows fron $J\beta_t = \beta_t J$.
Therefore $W_t := U_t V_t$ is a \strongH continuous one-parameter group
of $\cU(\cH)$. It satisfies
$JW_t J = V_t U_t = U_t V_t = W_t.$
Further $Z_t := \beta_t W_t^{-1}$ commutes with $\cM$ and $\cM'$, hence is contained
in the center of $\cM$. We conclude that
$Z_t = J Z_t J = Z_t^* = Z_t^{-1}$, and thus $Z_t^2 = \1$, which in turn implies that
$Z_t = \1$.

(ii) is clear from the definitions.

(iii) follows immediately from $J \beta_t J = \beta_t$
(Remark~\ref{rem:aut-w*}) because $J$ is antilinear.
\end{prf}

\begin{prop}
For any  $W^*$-dynamical system $(\cM,\R, \alpha)$, the subspace
$\cM_c \subeq \cM$ is the closed subalgebra generated by the
elements with bounded Arveson spectrum.
\end{prop}

\begin{prf}
Every element $M$ with bounded spectrum lies in a closed subspace on
which the action is norm continuous (\cite[Prop.~3.2.41]{BR02}), so that $M \in \cM_c$.
Conversely,  hitting an element $M \in \cM_c$ with an approximate identity
$(u_n)_{n \in \N}$ of $L^1(\R)$ for which the supports $\supp(\hat{u}_n)$
are compact leads to elements $\alpha_{u_n}(M)$ with bounded spectrum
converging to $M$.
\end{prf}

\subsection{The Borchers--Arveson Theorem and minimal implementing groups.}
\label{BA-sect}

We first consider an easily proven result which
shows a connection between spectral properties
and innerness of covariant representations.
For a locally compact abelian group $G$ and
a continuous unitary representation $(U,\cH)$ of $G$,
we write $\Spec(U) \subeq \hat G$ for its spectrum, i.e., the
support of the corresponding spectral measure on~$\hat G$.

\begin{lem} \mlabel{lem:longo} {\rm(Longo's Lemma)}
Let $(U,\cH)$ be a continuous unitary representation
of the abelian locally compact group $G$ and
$\cM \subeq \cB(\cH)$ be a von Neumann algebra normalized by
$U_G$. Suppose that
\begin{itemize}
\item[\rm(i)] there exists an $\cM$-cyclic unit vector
$\Omega\in\cH$ fixed by $U_G$, and that
\item[\rm(ii)] $\Spec(U) \cap \Spec(U)^{-1} \subeq \{e\}$.
\end{itemize}
Then $U_G \subeq \cM$.
\end{lem}

\begin{prf} We consider the action of $G$ on the commutant $\cM'$
defined by $\beta_g(M) := U_g M U_g^*$. We have to show that
$\beta$ is trivial; then $U_G \subeq \cM'' = \cM$.
As $\Omega$ is cyclic for $\cM$, it is separating for $\cM'$.
Let $\cE:=\lbr \cM'\Omega\rbr$, with projection $E:\cH\to\cE$ onto it,
and note that $U_G$ preserves $\cE$.
As $\cE\ni\Omega$ is $\cM\hbox{-generating,}$ it follows from
Lemma~\ref{CP1cyclic}(iv) that the restriction map $\cM'\mapsto\cM'\restriction \cE$
is an isomorphism. Thus it suffices to prove that the $W^*$-dynamical system
$(\cM'E,G, \beta^E)$ is trivial, where $\beta_g^E:=\Ad U_g^E$
and $U_g^E:=U_g\restriction\cE$.
As $\Omega$ is cyclic, separating and invariant for this $W^*$-dynamical system,
it follows from \cite[Prop.~XI.1.24]{Ta03} that
$\Spec(U^E) = \Spec(\beta^E)$. By $\Spec(U^E)\subseteq\Spec(U)$,
and condition (ii) we conclude that
$\Spec(\beta^E) \cap \Spec(\beta^E)^{-1}\subeq \{e\}$.
However,
the Arveson spectrum of an automorphic action is symmetric, i.e.
 $\Spec(\beta^E) = \Spec(\beta^E)^{-1}$ (Lemma~\ref{lem:4.11}(iii)), hence
 $\Spec(\beta^E) =\{e\}$, i.e. $\beta^E$ is trivial.
\end{prf}

Note that the preceding lemma applies in particular
to  covariant representations of actions of $G = \R$ with positive spectrum.

For covariant representations with positive spectrum of von Neumann algebras,
we have the stronger, and very important Borchers--Arveson Theorem
(cf.~\cite[Thm.~3.2.46]{BR02}), which conversely, gives us a way of constructing the spectral projections  of
an implementing unitary group from the Arveson spectral subspaces.
\begin{thm} \mlabel{BA-thm} {\rm(Borchers--Arveson)}
Let $(\cM,\R, \alpha)$ be a $W^*$-dynamical system on a von Neumann algebra
$\cM\subseteq\cB(\cH)$. Then the following are equivalent.
\begin{itemize}
\item[\rm(i)] There is a strong operator continuous unitary one-parameter group
$U:\R\to\cU(\cH)$ with positive spectrum, such that $\alpha_t=\Ad U_t$ on $\cM.$
\item[\rm(ii)] There is a  strong operator continuous unitary one-parameter group
$U:\R\to\cM$ with positive spectrum, such that $\alpha_t=\Ad U_t$ on $\cM.$
\item[\rm(iii)] Let $\cM^\alpha(S)$ denote the Arveson spectral subspace
for $S\subseteq\R$. Then
\[
\bigcap_{t\in\R}\br\cM^\alpha[t,\infty)\cH.=\{0\}.
\]
\end{itemize}
If these conditions hold, then we may take  $U:\R\to\cM$ to be
$U_t=\int_{\R}e^{-itx}dP(x),$ where $P$ is the projection-valued measure uniquely determined by
\[
P[t,\infty)\cH=\bigcap_{s<t}\br\cM^\alpha[s,\infty)\cH..
\]
\end{thm}
\begin{prf}
(i)$\Rightarrow$(iii): Let $P$ denote the projection valued measure of $\U$. As $\U$ has
positive spectrum, $P[0,\infty)=\1$, hence, using (\ref{MalphChar}), we obtain
\[
\cM^\alpha[t,\infty)\cH=\cM^\alpha[t,\infty)P[0,\infty)\cH\subseteq P[t,\infty)\cH\,.
\]
Thus, as $P$ is a projection-valued measure,
\[
\bigcap_{t\in\R}\br\cM^\alpha[t,\infty)\cH.\subseteq\bigcap_{t\in\R}P[t,\infty)\cH =\{0\}
\]
which proves (iii).

(iii)$\Rightarrow$(ii): In this proof we will let $[S]\in \cB(\cH)$ denote the
orthogonal projection onto the space $\br S.$. For $t\in\R$ define
\[
Q_t:=\Big[\bigcap_{s<t}{\br\cM^\alpha[s,\infty)\cH.}\Big]\in \cB(\cH).
\]
Then the map $t\mapsto \1- Q_t$ is a spectral family, as it is an increasing,  left \strongH continuous
map such that $\1-Q_t=0$ if $t\leq 0$ and it increases \strongH to $\1$ as $t\to\infty$
(cf. \cite[Def.~7.11]{Wei80}). Thus  there is a unique  projection valued measure $P$ such that
$P[t,\infty)=Q_t$ for all $t\in\R$. As the subspaces $\br\cM^\alpha[t,\infty)\cH.$ are invariant
with respect to  $\cM'$ their projections are in $\cM''=\cM$ and hence $P[t,\infty)\in\cM$ for all $t\in\R$.
Define
\[
U_t:=\int e^{-itp}dP(p) \in\cM
\]
then by $P[0,\infty)=\1$ it has positive spectrum. Define $\beta_t:={\rm Ad}(U_t)\in{\rm Aut}\cM$. As
\[
\cM^\alpha[s,\infty)P[t,\infty)\cH=\bigcap_{r<t}{\br\cM^\alpha[s,\infty)\cM^\alpha[r,\infty)\cH.}
\subseteq \bigcap_{r<t}{\br\cM^\alpha[s+r,\infty)\cH.}=P[s+t,\infty)\cH
\]
we obtain from (\ref{MalphChar}) that $\cM^\alpha[s,\infty)\subseteq\cM^\beta[s,\infty)$ for all $s\in\R$.
Thus by Proposition~\ref{ArvSpecAct} we get that $\alpha_t=\beta_t$.

(ii)$\Rightarrow$(i) is trivial.
\end{prf}

\begin{rem}
\label{RemBorAr}
(a) The theorem gives a sharp criterion stating when we have a covariant representation
with positive spectrum. It
states that amongst the implementing unitary one-parameter groups with positive spectrum, we can find
one which is inner, and it selects one by construction.
Hence by (ii) in Theorem~\ref{BA-thm},
 every normal representation of $\cM$ is covariant.

Moreover, given a covariant representation $(\pi,U)$ with positive spectrum of a
$C^*$-action $(\cA, G, \alpha)$, we can always find
a new  covariant representation $(\pi,V)$ with positive spectrum,
such that its generator is affiliated with
$\pi(\cA)$.

(b) An important consequence of the Borchers--Arveson Theorem~\ref{BA-thm}
is that for any covariant representation $(\pi, U)$ of $(\cM,\R,\alpha)$ for which
$\pi$ is faithful and $U$ has positive spectrum,
the action $\alpha$ is trivial on the center of $\cM.$
Hence, $\cM$ must be non-commutative in order to admit non-trivial
actions and covariant representations with positive spectrum.
Moreover, as  any  commutative $C^*$-subalgebra of $\cM$ preserved by the action $\alpha$
must be in  $\cM^{\R}$, it follows that $\alpha$
cannot have any normal eigenvectors except for the identity eigenvalue.
It seems that a
$[0,\infty)$-spectral condition is a quantum mechanical phenomenon, which cannot occur in
classical systems. It is now easy to give examples of actions for which there are
covariant representations, but no
covariant representations with positive spectrum, e.g.\ the translation action of $\R$ on $C_0(\R)$.

(c) The Borchers--Arveson Theorem has been generalized by Kishimoto to $\R^n$ \cite[Thm.~2]{Ki79},
and further to  connected locally compact abelian groups in \cite[Cor.~8.4.12]{Pe89}.
\end{rem}

Apart from the observations in Remark~\ref{RemBorAr}(b),
the existence of a covariant representation with positive spectrum places
strong algebraic restrictions on the  $C^*$-action
 $(\cA, G, \alpha)$.
This is explored in Section~\ref{PosObs}.

By the Borchers--Arveson Theorem~\ref{BA-thm}, if we have an implementing  unitary one-parameter group $U:\R\to\cU(\cH)$
with positive spectrum, we may
take it to be inner, and then  $U(C^*(\R))\subset U_{\R}''\subset\cM^\alpha\{0\}=\cM^{\R}.$
Above we saw that the Arveson spectral subspaces can be written in terms of ``matrix decompositions''
with respect to  $C^*(\R)$ (cf.~Example~\ref{Zmatrix} and preceding discussion).
Thus the subalgebra $\cM^\alpha\{0\}$
already contains the spectral information of  $(\cM,\R, \alpha)$
because it contains all the spectral projections of~$U$.

The set $\cU$ of unitary one-parameter groups with positive spectrum implementing $\alpha$,
has an interesting structure. It has a partial order, obtained from the generators of the groups;-
let  $(U_t)_{t \in \R},\;(V_t)_{t \in \R}\subset\cB(\cH)$ be in $\cU$,
and write $U_t=\exp(itA)$ and $V_t=\exp(itB)$ for $A,\;B\geq 0$. Then
$U_t\leq V_t$ iff $A^n\leq B^n$ for all $n=1,\,2,\,3,\ldots$ iff $P^U[t,\infty)\leq P^V[t,\infty)$
for all $t\in\R$ where $P^U$ denotes the projection valued measure associated with $U$ (cf.~\cite[(3.2), p235]{Arv74}).
Then $\cU$ with this partial ordering has minima. Those minima in $\cM$ are particularly interesting, in that
they are least elements over all of $\cU,$ hence there can be at most be one minimum in $\cM.$
 We will show below in Lemma~\ref{lem:minimal} the existence of it, using the Borchers--Arveson theorem.
 We concentrate only on this inner minimal positive group in $\cM.$

\begin{defn} \mlabel{def:minimal}
Let $\cM \subeq \cB(\cH)$ be a von Neumann algebra, and
 $(\cM,\R, \alpha)$ be a $W^*$-dynamical system.

(a) Let $(U_t)_{t \in \R}\subset\cM$ be a \strongH continuous
unitary one-parameter group with positive spectrum implementing $\alpha$.
We say that $(U_t)_{t \in \R}\subset\cM$ is {\it  minimal} if,
for all other one-parameter groups $(V_t)_{t \in \R}\subset \cB(\cH)$ with positive
spectrum implementing the same automorphisms, i.e.\
$\Ad(U_t) = \Ad(V_t)=\alpha_t$ for $t \in \R$,
the corresponding one-parameter group
$Z_t := V_t U_t^*\in\cM'$ has positive spectrum
(note that $(U_t)_{t \in \R}\subset\cM$ ensures that $Z_t$ is a one-parameter group).
A~minimal one-parameter group in $\cM$ is clearly unique, if it exists, and will be called the
{\it inner minimal positive} one-parameter group.

(b) The set of projections $(Q_t)_{t \in \R}\subset \cB(\cH)$ defined by
\begin{equation}
  \label{eq:Qa}
Q_t \cH := \bigcap_{s < t} \lbr\cM^\alpha[s,\infty) \cH\rbr
\end{equation}
are called {\it Borchers--Arveson projections}. We also put
\begin{equation}
  \label{eq:Q2}
Q_\infty := \lim_{t\to \infty} Q_t
\end{equation}
and observe that the limit exists because $Q_s \leq Q_t$ for $s \geq t$.
(Note that $Q_t\in\cM\ni Q_\infty$ by the bicommutant theorem, as
$\cM'$ preserves $\lbr\cM^\alpha[s,\infty) \cH\rbr$.)\\
If $\cH_\infty := Q_\infty\cH= \{0\}$, then
the unitary one-parameter group
$(U_t)_{t\in \R}\subset\cM$
whose spectral measure~$P$ is determined by
$P[t,\infty) = Q_t$ for $t \in \R$ is called  the
{\it Borchers--Arveson group} for the
$W^*$-dynamical system $(\cM,\R, \alpha)$
(cf. Theorem~\ref{BA-thm}).
\end{defn}

In terms of its generator, the unitary group  $U_t=\exp{(-itA)}\in \cM$, $A\geq 0$, is minimal if
for all other one-parameter groups $(V_t)_{t \in \R}=\exp{(-itB)}\subset \cB(\cH)$,
$B\geq 0$ such that
$\Ad(U_t) = \Ad(V_t)$ for $t \in \R$, we have
$B\geq A$.

The first part of the following lemma is \cite[Prop.~p.~235]{Arv74}.

\begin{lem} \mlabel{lem:minimal}
Suppose that $Q_\infty = 0$. Then the
Borchers--Arveson subgroup $(U_t)_{t \in \R}$ in $\cM$
is minimal and, for every $\eps > 0$ the projection
$P[0,\eps) = Q_0 - Q_\eps\in\cM$ has central support~$\1$.
\end{lem}

\begin{prf}
From the formula
\begin{equation}
  \label{eq:can-spec}
P[t,\infty)\cH = \bigcap_{s < t} \lbr\cM^\alpha[s,\infty)\cH\rbr
\end{equation}
for the spectral projections of $U$, one derives as follows
that $P[0,\eps) \not=0$ for any $\eps > 0$:
If $m := \inf \Spec(U)$ and $0 < s < \eps$, then $\cH = \cH^U[m,\infty):=P[m,\infty)\cH$, so that
\[  P[\eps, \infty)\cH \subeq \cM^\alpha[s,\infty)\cH
= \cM^\alpha[s,\infty)\cH^U[m,\infty) \subeq \cH^U[m+s,\infty),\]
is a proper subspace of $\cH$.
Since this remains valid for every subrepresentation,
the central support of the projections $P[0,\eps)$, $\eps > 0$, is $\1$.

To see that $U$ is minimal
(\cite[Thm.~II.4.6]{Bo96} or \cite[Prop.~p.~235]{Arv74}), let $(\tilde U_t)_{t \in \R}$ be another
strongly continuous unitary one-parameter group implementing
the same automorphisms, i.e.\  $\Ad(U_t) = \Ad(\tilde U_t)$.
As $( U_t)_{t \in \R}\subset\cM,$ we have that $ ( U_t)_{t \in \R}$ and
$(\tilde U_t)_{t \in \R}$ commute, hence they can be diagonalized simultaneously.
Then the spectral measure
$\tilde P$ of $\tilde U$ satisfies
\[ \cM^\alpha[s,\infty)
\tilde P[t,\infty)\cH \subeq \tilde P[t+s,\infty)\cH \quad \mbox{ for } \quad
t,s \in \R,\]
so that, for $t \in \R$,
\begin{align*}
P[t,\infty)\cH
&= \bigcap_{s < t} \lbr\cM^\alpha[s,\infty) \cH\rbr
= \bigcap_{s < t} \lbr\cM^\alpha[s,\infty) \tilde P[0,\infty)\cH\rbr
\subeq \bigcap_{s < t} \tilde P[s,\infty)\cH
= \tilde P[t,\infty)\cH.
\end{align*}
We conclude that $P[t,\infty) \leq \tilde P[t,\infty)$ for $t \in \R$.
We prove that this implies that $H \leq \tilde H$ holds for the infinitesimal generators of
$U$ and $\tilde U$, respectively
(see \cite[Prop. 6.3]{PS12} for the inclusion $\cD({\tilde H})\subseteq\cD(H)$), and hence that $U$ is minimal.

We first show for the domains that $\cD(\tilde H)\subseteq\cD(H).$
For $n\in\N$, we approximate $H$ from below by using the step function
$f_n := \frac{1}{n}\sum_{k=1}^\infty \chi_{[\frac{k}{n},\infty)}$
to get
$H\geq H_n:=f_n(H)$ and note that the operators $H$ and $H_n$ have
the same domain, as they differ by a bounded operator, and $\|H-H_n\|\leq\frac{1}{n}$.
Likewise $\tilde H\geq \tilde H_n:=f_n(\tilde H),$ $\cD(\tilde H) = \cD(\tilde H_n)$
and $\|\tilde H-\tilde H_n\|\leq\frac{1}{n}$. It suffices to show that $\cD(\tilde H_1)\subseteq\cD(H_1).$
Now
$H_n=\frac{1}{n}\sum_{k=1}^\infty P[\frac{k}{n},\infty),$ where the convergence of the sum is on
vectors $\xi \in \cD( H) = \cD(H_n).$ Moreover,
\[
\xi \in \cD( H_1)\quad\hbox{if and only if}\quad\lim_{K\to\infty} \|A_K\xi\|^2<\infty,\qquad\hbox{with}\quad
A_K:=\sum\limits_{k=1}^KP[k,\infty).
\]
We have likewise expressions for $\tilde H_n$. For any $\xi\in\cH$,
\begin{eqnarray*}
\|A_K\xi\|^2&=& \big(A_K\xi,\,A_K\xi\big)=\big(\xi,\,(A_K)^2\xi\big)
= \sum\limits_{k=1}^K\big(\xi,\,   (2k-1)\,P[k,\infty)\xi\big)\\[1mm]
&\leq& \sum\limits_{k=1}^K\big(\xi,\,   (2k-1)\,\tilde P[k,\infty)\xi\big)
=\|\tilde A_K\xi\|^2,
\end{eqnarray*}
using   $P[k,\infty)P[j,\infty) =  P[\max\{k,j\},\infty)$ to simplify $(A_K)^2$. Thus
\[
\lim_{K\to\infty}\|\tilde A_K\xi\|^2<\infty\quad\Rightarrow\quad
\lim_{K\to\infty}\| A_K\xi\|^2<\infty
\]
i.e. $\cD(\tilde H_1)\subseteq\cD(H_1),$ hence $\cD(\tilde H)\subseteq\cD(H).$

Let $\xi \in \cD(\tilde H)\subseteq\cD(H)$. Then  $H \leq \tilde H$ follows from
\[ \la \xi,H\xi \ra
= \lim_{n \to \infty} \la \xi, H_n \xi \ra
\leq  \lim_{n \to \infty} \la \xi, \tilde H_n \xi \ra = \la \xi, \tilde H \xi \ra.\qedhere\]
\end{prf}
Thus for a covariant representation with positive spectrum, the
inner minimal positive one-parameter group exists, and coincides with the Borchers--Arveson group in $\cM.$
\begin{lem}
  \mlabel{lem:mini}
Let $\cM \subeq \cB(\cH)$ be a von Neumann algebra and
$(\cM,\R,\alpha)$ a $W^*$-dynamical system.
A unitary one-parameter subgroup $(U_t)_{t \in \R}\subset\cM$ with non-negative spectrum
implementing $\alpha$ on $\cM$ is minimal if and only
if, for every $\eps > 0$, the central support of $P[0,\eps)$ is~$\1$.
\end{lem}

\begin{prf} If $U$ is minimal, then it coincides with the Borchers--Arveson
subgroup in a faithful normal representation of $\cM$. Hence the central
support of every $P[0,\eps)$, $\eps > 0$, is $\1$ by Lemma~\ref{lem:minimal}.

Assume, conversely, that the central support of every $P[0,\eps)$, $\eps > 0$,
is $\1$. As  $(U_t)_{t \in \R}\subset\cM$ has positive spectrum, $\cM$ also contains the inner minimal
positive implementing group, and we only need to compare $(U_t)_{t \in \R}$ with that.
Thus we have to show that, for every central
subgroup
$Z_t = e^{itW}\in Z(\cM)$ for which $(U_t Z_t)_{t \in \R}$ has non-negative spectrum, we have
$W \geq 0$. We argue by contradiction. If $W$ is not positive, then the
corresponding spectral projection $P^W((-\infty,-2\eps])$ is non-zero for some $\eps > 0$.
Our assumption implies that $P^W((-\infty,-2\eps]) P[0,\eps) \not=0$
 in any normal representation,
hence $H+W$ is negative on the range of this projection, where $H$ is the infinitesimal implementer of
 $(U_t)_{t \in \R}$.
Therefore $H + W$ is not positive, which contradicts the assumption that
$(U_t Z_t)_{t \in \R}$ has non-negative spectrum.
\end{prf}

From this we obtain that normal representations take \minimal groups to \minimal groups:
\begin{lem}
\mlabel{mintomin}
Let $\cM \subeq \cB(\cH)$ be a von Neumann algebra and
$(\cM,\R,\alpha)$ a $W^*$-dynamical system. Let  $(U_t)_{t \in \R}\subset\cM$
be the \minimal implementing unitary group  for $\alpha$. If $\pi:\cM\to\cB(\cH_\pi)$
is a normal representation, then $(\pi(U_t))_{t \in \R}\subset\pi(\cM)$ is the \minimal
implementing unitary group for $(\pi(\cM),\R,\alpha_\pi)$, where
$\alpha_\pi(t)A = \pi(U_t)A\pi(U_t)^*$.
\end{lem}

\begin{prf}
By the previous lemma, it suffices to prove that for every $\eps > 0$, the central support of $\pi(P[0,\eps))$ is~$\1$. If $Z := s(\pi)$ is the support of
$\pi$, then $\pi(\cM) \cong Z\cM$ and $\pi(P[0,\eps))$ corresponds to
$Z P[0,\eps) \in Z\cM$. If $Z' \in Z\cM$ is a central projection
with $0 = Z' Z P[0,\eps) = Z' P[0,\eps)$, then $Z' = 0$ follows from
the fact that $Z'$ is also central in $\cM$ and
the central support of $P[0,\eps)$ is~$\1$.
\end{prf}

We will use these lemmas in the next subsection when we study the structure of covariant
representations with positive spectrum.
We next show that every covariant representation contains a maximal subrepresentation which
satisfies the Borchers--Arveson criterion (Theorem~\ref{BA-thm}),
which we then apply to the universal  covariant representation
$(\pi_{co},U_{co})\in{\rm Rep}(\alpha,\cH_{co})$. This is in fact already known through the ``minimal
covariant subrepresentation with positive spectrum''
constructed in either \cite[Thm.~II.4.6]{Bo96} or \cite[Thm~8.4.3]{Pe89}, but we will need to make some
of its details explicit.

\begin{Lemma}
\mlabel{Hinfty}
Let $(\cM, \R,\alpha)$ be a $W^*$-dynamical system on a von Neumann algebra
$\cM\subseteq\cB(\cH)$. Then
\begin{itemize}
\item[\rm(i)] $\cM=\overline{\bigcup\{\cM^\alpha[s,\infty)\mid s\in\R\}}^{\,\rm w-op}$.
\item[\rm(ii)] The space $\cH_\infty := Q_\infty \cH
= \bigcap\limits_{s\in\R}\br\cM^\alpha[s,\infty)\cH.$
is an invariant subspace for  $\cM\cup\cM'$, i.e.\  $Q_\infty \in \cM \cap \cM'$,
and in the case that $(\cM,U)$ is covariant, $\cH_\infty$
 is also $U$-invariant.
\end{itemize}
\end{Lemma}

 \begin{prf} (i)
 By \cite[Lemma~3.2.38(3)]{BR02}, we know that, for
$f \in L^1(\R)$ such that ${\rm supp}\widehat{f}\subseteq [s,\infty)$
and $A \in \cM$, we have $\alpha_f(A) \in \cM^\alpha[s,\infty)$.
 Let $f\in L^1(\R)$ be such that $\widehat{f}$ is a smooth function with support in $[-1,1]$, and normalized
 such that $\int_\R |f|dt =1$ (note that both $f$ and $\widehat{f}$ are Schwartz functions).
  Let $f_n(x):=nf(nx)$. Then  $\int_\R |f_n|dt =1$ and $\widehat{f_n}(p)=\widehat{f}(p/n) $
  which has support in $[-n,n]$. Moreover the $f_n$ are progressively narrower concentrated around $0$,
  i.e.\ given any $a>0$ and an $\varepsilon>0$, then there is an $N\in\N$ such that
 $\int_{-a}^a|f_n|dt>1-\varepsilon$   for $n>N$.
Note that all $\alpha_{f_n}(A)\in \bigcup\{\cM^\alpha[s,\infty)\mid s\in\R\}$ for every $n \in \N$.
   We want to show that $\alpha_{f_n}(A)\to A$ in the weak operator topology for all $A\in\cM$.

   Let $\omega$ be a vector state on $\cM$, and fix $A\in\cM$, so that
$t\to\omega(\alpha_t(A))$ is continuous. For $\varepsilon>0$,
there is a $\delta>0$ such that $|t|<\delta$ implies that
   ${\big|\omega(A-\alpha_t(A))\big|}<\varepsilon$. Note that
${\big|\omega(A-\alpha_t(A))\big|}\leq 2\|A\|$
   for all $t$.
   Then there exists an  $N\in\N$ such that
  $\int_{-\delta}^\delta|f_n|dt>1-\varepsilon$   for all $n>N$. Then
  \begin{eqnarray*}
\big|\omega(A-\alpha_{f_n}(A))\big|  &\leq& \int|f_n(t)|\,\big|\omega(A-\alpha_t(A))\big|\,dt\\[1mm]
&=&\Big(\int_{(-\delta,\delta)}+\int_{\R\backslash(-\delta,\delta)}\Big)|f_n(t)|\,\big|\omega(A-\alpha_t(A))\big|\,dt <\varepsilon+2\|A\|\varepsilon\,.
   \end{eqnarray*}
   Thus $\alpha_{f_n}(A)\to A$ in the weak operator topology, which proves part~(i).

   (ii) According to \cite[Lemma~3.2.39(2)]{BR02}, we have
   $U_t\cM^\alpha[s,\infty)=\cM^\alpha[s,\infty)U_t$ for all $s,\,t\in\R$, hence the last claim is clear.
   As
   \[
   \cM'\br\cM^\alpha[s,\infty)\cH.=\br\cM^\alpha[s,\infty)\cM'\cH.=\br\cM^\alpha[s,\infty)\cH.
   \]
   it is also clear that $\cM'\cH_\infty\subseteq\cH_\infty$. Finally, by
\cite[Lemma~3.2.42(4)]{BR02}, we have
\[    \cM^\alpha[s,\infty)\cdot\cM^\alpha[t,\infty)\subseteq
\cM^\alpha[s+t,\infty)\]
and hence
\[   \cM^\alpha[s,\infty) \cH_\infty\subseteq
\bigcap_{t\in\R}\br\cM^\alpha[s+t,\infty)\cH.
  =\cH_\infty. \]
   As $\cH_\infty$ is closed, it follows from part~(i) that $\cM\cH_\infty\subseteq
   \overline{\bigcup\limits_{s\in\R}\cM^\alpha[s,\infty)\cH_\infty}\subseteq\cH_\infty$.
 \end{prf}

The Borchers--Arveson Theorem~\ref{BA-thm} states that $(\cM,\R,\alpha)$ has a
 strong operator continuous unitary one-parameter implementing group
with positive spectrum if and only if $\cH_\infty=\{0\}.$
 This indicates how to select a state for which its GNS representation has a implementing unitary group
with positive spectrum for $\alpha$ (see below).

In the context of this lemma,
let $Q_\infty$ be the orthogonal projection onto $\cH_\infty$. It follows from
Lemma~\ref{Hinfty}(ii) that $Q_\infty\in\cM'\cap\cM''=Z(\cM)$, hence $\cM$ is
diagonal with respect to  the decomposition
 $\cH=\cH_\infty\oplus\cH_\infty^\perp=:\cH_\infty\oplus\cH^{(+)} $. Let
$P^{(+)}:={\bf 1}-Q_\infty$. Then $\cM$ is the direct  sum of the two ideals
 $\cM_\infty:= \cM Q_\infty$ and $\cM^{(+)}:=\cM P^{(+)}$.
 Define the {\it  subrepresentation with positive spectrum} of $\cM$ to be the representation
 $\pi^{(+)}:\cM\to\cB(\cH^{(+)})$ by $\pi^{(+)}(A):=A\restriction\cH^{(+)}$, $A\in\cM$, then clearly
 $\pi^{(+)}(\cM)\cong \cM^{(+)}$. Its name is justified by the following proposition:
 \begin{prop}
 \mlabel{posrep}
 Let  $U:\R\to\cU(\cH)$  be a strong operator continuous unitary one-parameter group
such that $\alpha_t:=\Ad U_t$ defines  an action
 $\alpha:\R\to \Aut(\cM)$  on  a given von Neumann algebra
$\cM\subseteq\cB(\cH)$. Then  its subrepresentation with positive spectrum
 $\pi^{(+)}:\cM\to\cB(\cH^{(+)})$ has the following properties:

\begin{itemize}
\item[\rm(i)] There is a strong operator continuous unitary one-parameter group
with positive spectrum ${V:\R\to\cU(\cH^{(+)})}$ such that $\alpha_t=\Ad V_t$ on $\pi^{(+)}(\cM)$.
 This unitary implementing group may be chosen to be inner.
\item[\rm(ii)] $\pi^{(+)}$ is maximal, in the sense that
any  subrepresentation of $\cM$ to which $U_t$ restricts, and which has a implementing
unitary group with positive spectrum, must be contained in the subrepresentation with positive spectrum.
\end{itemize}
 \end{prop}
 \begin{prf}
 We first need to prove that if
 $\cH_1\subset\cH$ is a subspace invariant with respect to  $\cM$ and $U_{\R}$, then the spectral subspaces restrict.
That means, if we label the subrepresentation by \break
 $\pi_1:\cM\to\cB(\cH_1)$, $\pi_1(A):=A\restriction\cH_1$, $A\in\cM$, then
 $\pi_1(\cM^\alpha[s,\infty))=\pi_1(\cM)^\beta[s,\infty)$ for all $s\in\R$,
 where $\beta_t:={\rm Ad}(U_t\restriction\cH_1)$. But this follows from
 the characterization~(\ref{MalphChar}) since the spectral projection of $U_t$ commutes
 with the projection onto $\cH_1$. If we let $\cH_1=\cH^{(+)}$, then the spectral subspaces
 of $\beta_t$ are the projections of the spectral subspaces of $\alpha_t$ by $P^{(+)}$,
 hence by construction $\beta_t$ satisfies the condition of the
 Borchers--Arveson Theorem~\ref{BA-thm}, and this proves (i).  Then it follows  that the subrepresentation
 with positive spectrum of
its orthogonal subrepresentation is zero, which is equivalent to (ii) by the
Borchers--Arveson Theorem~\ref{BA-thm}.
 \end{prf}

Given a $C^*$-action  $(\cA, \R, \alpha)$,
consider the universal  covariant representation\break
$(\pi_{co},U_{co})\in{\rm Rep}(\alpha,\cH_{co})$ with associated $W^*$-dynamical system
$(\cM_{co}, \R, \alpha^{co})$. Then the  subrepresentation with positive spectrum
${\pi^{(+)}_{co}}:\cA\to\cB(\cH^{(+)}_{co})$ has the universal property that every
cyclic
covariant representation with positive spectrum of $(\cA, \R, \alpha)$
is unitarily equivalent to a subrepresentation of it.
Moreover, it is also unitarily equivalent to
 the ``minimal representation with positive spectrum''
constructed in \cite[Thm~II.4.6]{Bo96} and \cite[Thm.~8.4.3]{Pe89}.

 Consider a state $\omega\in\fS(\cA)$ which is quasi-invariant, i.e.
  $\pi_\omega$ is quasi-covariant (cf. Def.~\ref{eq:covar}(c)). Then
$\alpha$ induces a $W^*$-dynamical system $(\cM,\R,\beta)$, where
$\cM := \pi_\omega(\cA)''$. Moreover $(\cM,\R,\beta)$ has a strong operator continuous unitary
one-parameter implementing group  with positive spectrum if and only if $(\pi_\omega(\cA),\R,\alpha)$ has.
In view of the Borchers--Arveson Theorem~\ref{BA-thm}, this is equivalent to
 \[ \{0\}=\cH_{\omega,\infty}:=\bigcap\limits_{s\in\R}\br\cM^\beta[s,\infty)\cH_\omega.\]
Thus any equivalent condition to $\cH_{\omega,\infty}=\{0\}$
would characterize the set of such states  with implementing group with positive spectrum:

 \begin{prop}
 \mlabel{posrep1}
For a $C^*$-action $(\cA, \R, \alpha)$,  define
\[\fS_{co}^{(+)}(\cA):=\{\omega\in\fS(\cA)\mid (\pi_\omega,V)\in{\rm Rep}(\alpha,\cH_\omega)\quad\hbox{for some}\quad
V:\R\to\cU(\cH_\omega)\quad\hbox{with positive spectrum}\}.
\]
For a quasi-invariant state $\omega\in\fS(\cA)$,
let  $Q_\infty^\omega \in \pi_\omega(\cA)''$ be the orthogonal projection onto
$\cH_{\omega,\infty}$, and let $\omega$ also
denote its extension  to $\pi_\omega(\cA)''$ as the vector
state ${(\Omega_\omega,\cdot\,\Omega_\omega)}.$ Then
\[
\omega\in\fS_{co}^{(+)}(\cA)\quad \Longleftrightarrow
\quad\omega(Q_\infty^\omega)=0\quad \Longleftrightarrow \quad
\omega(Q_\infty^\omega A)=0\quad \mbox{ for all } \quad A\in\cA.
\]
 \end{prop}

 \begin{prf} Let $\cM := \pi_\omega(\cM)''$.
By Lemma~\ref{Hinfty}(ii) we have that $Q_\infty^\omega\in\cM'\cap\cM''=Z(\cM)$.
From the Cauchy--Schwartz inequality
 \[
 |\omega(Q_\infty^\omega A)|^2\leq \omega(Q_\infty^\omega)\omega(A^*A)\quad \mbox{ for } \quad A\in\cA, \]
we get that  $\omega(Q_\infty^\omega )=0$  implies $\omega(Q_\infty^\omega A)=0$ for all $A\in\cA$.
Conversely, as $\pi_\omega(\cA)$ acts non-degenerately on $\cH_\omega$, $\pi_\omega(E_\lambda)\to\1$ in strong operator topology
for any approximate identity $(E_\lambda)_{\lambda\in\Lambda}$ in $\cA$,
hence if $\omega(Q_\infty^\omega A)=0$ for all $A\in\cA$
then $\omega(Q_\infty^\omega )=\lim\limits_\lambda\omega(Q_\infty^\omega E_\lambda)=0$ which gives the converse implication,
and hence the second equivalence is established.
Moreover, if $\omega(Q_\infty^\omega )=0$ then also all $\omega_B(Q_\infty^\omega )=0$
 where $\omega_B(A):=\omega(B^*AB)$ for $A,\,B\in\cA$, $\|B\|=1$, and hence all vector states of $\pi_\omega$ will also satisfy it.
The vector state of any vector orthogonal to $\cH_{\omega,\infty}$ clearly satisfies the condition, whereas any nonzero vector
$\psi\in\cH_{\omega,\infty}$ produces $\omega_\psi(Q_\infty^\omega)=\|\psi\|^2\not=0$. Thus
the condition $\omega(Q_\infty^\omega)=0$ is equivalent to $\cH_{\omega, \infty}=\{0\}$, which
by the Borchers--Arveson Theorem~\ref{BA-thm} characterizes
$\fS_{co}^{(+)}(\cA)$.
 \end{prf}

This condition looks different in  Borchers approach (cf.~\cite[Def.~II.4.3(i)]{Bo96})
as his selection condition is
 $\omega E(\infty)=E(\infty)\omega=\omega$ where $E(\infty)$ is  the projection
 onto the subspace $\cH_{co}\cap\cH_{co,\infty}^\perp$ in the universal
representation on $\cH_{co}$.
However, this condition clearly coincides with the condition above in the given context.

\subsection{Covariant representations with positive spectrum and obstruction results.}
\label{PosObs}

The Borchers--Arveson Theorem produces several obstruction results for
covariant representations with positive spectrum. By Remark~\ref{RemBorAr}(b),
for any covariant representation $(\pi, U)$ of $(\cM,\R,\alpha)$ for which
$\pi$ is faithful and $U$ has positive spectrum, if the action is nontrivial,
the algebra $\cM$ must be noncommutative. This obstruction result leads to
further obstructions, which we now discuss.

\begin{prop} \mlabel{prop:nogo1}
Let $\cB$ be a $C^*$-algebra and $\cA := C_0(\R,\cB)$, endowed with the
automorphisms $(\alpha_t f)(x) := f(x-t)$.
Then all  covariant representations with positive spectrum  $(\pi,U)$ of $(\cA,\R,\alpha)$
satisfy $\pi = 0$.
\end{prop}

\begin{prf}
Writing $\cA \cong C_0(\R) \otimes \cB$, we see that every non-degenerate representation
of $\cA$ can be written as
$\pi(A_1 \otimes A_2) = \pi_1(A_1) \pi_2(A_2)$, where
$\pi_1 \: C_0(\R) \to \cB(\cH)$ and
$\pi_2 \: \cB \to \cB(\cH)$ are commuting representations
(cf. \cite[Prop. 4.7, Lemma 4.18]{Ta02}). Hence every
covariant representation $(\pi, U)$ with positive spectrum of $(\cA,\R,\alpha)$ leads to a covariant
representation of $(C_0(\R), \R, \alpha)$, so that the Borchers--Arveson Theorem implies that
$U_\R$ commutes with $\pi_1(C_0(\R))$ and $\pi_2(\cB)$,
and this implies that $U_\R$ commutes with $\pi(\cA)$.

A function $f \in C_c(\R)$ is a derivative of a compactly supported function $F$ if and only
if $\int_\R f(x)\, d x= 0$. Then
\[ f(x) = \lim_{h \to 0} \frac{F(x+h) - F(x)}{h} \]
shows that we must have $\pi_1(f) =0$ for all these functions.
Now the density of
\[ \Big\{ f \in C_c(\R) \: \int_\R f(x)\, d(x)= 0\Big\}
= \{ f' \: f  \in C_c(\R)\} \]
in
$C_0(\R)$ implies $\pi_1=0$. This in turn leads to $\pi= 0$.
\end{prf}
The translation action can be twisted by a cocycle without affecting
the obstruction. To see this,
modify the construction as follows. On
$\cA := C_0(\R,\cB)$, we consider the automorphisms
\begin{equation}
\label{eq:twist}
(\alpha_t f)(x) = \beta_t(x)(f(x-t)),
\end{equation}
where $\beta \: \R \to C_b(\R,\Aut(\cB))$ is a cocycle in the sense that the translation automorphism
$(\alpha^0_t f)(x) := f(x-t)$ satisfies $\alpha_t = \beta_t \cdot \alpha_t^0$.
Then
\[\beta_{t+s} \alpha_{t+s}^0 =  \alpha_{t+s} = \alpha_t \alpha_s = \beta_t \alpha_t^0\beta_s \alpha_s^0
= \beta_t (\alpha_t^0 \beta_s \alpha_{-t}^0) \alpha_{t+s}^0\]
leads to the cocycle relation
\[\beta_{t+s}= \beta_t \cdot (\alpha_t^0 \beta_s \alpha_{-t}^0).\]
This means that
\[ \beta_{t+s}(x) = \beta_t(x) \beta_s(x-t)\quad \mbox{ for } \quad t,s,x \in \R.\]

\begin{cor}  Let $\alpha:\R\to\Aut(cA)$ for $\cA = C_0(\R,\cB)$ be defined as in \eqref{eq:twist}
for a cocycle  $\beta \: \R \to C_b(\R,\Aut(\cB))$.
Then all  covariant representations with positive spectrum $(\pi,U)$ of $(\cA,\R,\alpha)$
satisfy $\pi = 0$.
\end{cor}

\begin{prf} Let $(\pi, U)$
be a covariant representation of $(\cA, \R, \alpha)$ and observe that
it extends to a covariant representation of the multiplier algebra
$(M(\cA), \R, \alpha)$. In
$M(\cA)$ we have the subalgebra $C_0(\R,\C)$ obtained from the functions whose values are
multiples of $\1$. On this subalgebra the
$\R$-action takes the form $(\alpha_t f)(x) =f(x-t)$
because $f(\R) \subeq \cB$ is fixed by all automorphisms.
Then Proposition~\ref{prop:nogo1} implies that
$\pi(C_0(\R,\C))=\{0\}$. This in turn yields $\pi(\cA) = \{0\}$.
\end{prf}

\begin{rem} \mlabel{rem:twist}
If $\hat G \cong \T \rtimes_\gamma G$ is a central $\T$-extension of $G$, for a given 2-cocycle
$\gamma:G\times G\to\T$, then we associate
the corresponding twisted group $C^*$-algebra $\cA :=C^*_\gamma(G_d)$
defined by the unitary generators $(\delta_g)_{g \in G}$ satisfying the relations
\[ \delta_g \delta_h = \gamma(g,h) \delta_{gh}.\]

Any $\R$-action by automorphisms on $\hat G$ fixing the central subgroup
$\T$ pointwise induces a homomorphism
$\alpha \: \R \to \Aut(\cA)$.
Now covariant projective unitary representations
for the cocycle $\gamma$ correspond to covariant
representations of $(\cA,\R,\alpha)$.
\end{rem}

\begin{ex}
\label{WeylPos} {\rm (The Weyl algebra)}
Let $ \Heis(V,\sigma) = \T \times V$ be the Heisenberg group
of the real symplectic topological vector space $(V,\sigma)$ with the
multiplication
\[  (z,v) (z',v') := (zz' e^{-\frac{i}{2}\sigma(v,v')}, v + v'), \quad
z\in\T,\; v\in V \]
and let $\cA :=\overline{ \Delta(V,\sigma)}$ be the corresponding Weyl algebra,
which is the discrete twisted group algebra
$C^*_\gamma(V_d)$, where $\gamma(v,w) = e^{-\frac{i}{2}\sigma(v,v')}$.
We consider a smooth one-parameter group $(\tau_t)_{t \in \R}=e^{tY}\in\Sp(V,\sigma)$,
$Y \in \sp(V,\sigma)$. Here smoothness refers to the smoothness of the
$\R$-action $\R \times V \to V$.
This defines
an action $\alpha_0:\R\to \Aut(\Heis(V,\sigma))$ by $\alpha_{0,t}(z,v):=(z,\tau_t(v))$, and as
$\alpha_{0,t}$ fixes all $(z,0)$, it also defines an automorphic $\R$-action
$\alpha$ on $\cA$ which is singular, as it is not \strong continuous.
Now $\cA$ has many representations which are not continuous with respect to ~the underlying group
$ \Heis(V,\sigma)$ (nonregular representations), so to avoid these, we consider the
associated Lie groups.
\end{ex}
As the action $\alpha$ on $ \Heis(V,\sigma)$ is smooth, we  form the corresponding
{\it oscillator group}
\[ G:= \Heis(V,\sigma) \rtimes_\alpha \R.\]
It is a Lie group because the $\R$-action on $\Heis(V,\sigma)$ is smooth.
Now any smooth unitary representation $(\pi, \cH)$ of $G$ for which $\pi(z,0,0)=z\1$
will define a covariant
representation of $(\cA,\R, \alpha)$, where the unitary implementers of $\alpha_t$ are $U_t := \pi(0,0,t)$.
 We analyze positivity for these covariant representations.
\begin{prop} \mlabel{prop:osc} If $(\pi, \cH)$ is a smooth unitary representation
of $G= \Heis(V,\sigma) \rtimes_\alpha \R$ for which the one-parameter group
$U_t = \pi(0,0,t)$ has positive spectrum and
$\pi(z,0,0) = z\1$, then
the infinitesimal generator $Y\in \sp(V,\sigma)$ of $\tau$ satisfies
\[   \sigma(Yv,v) \geq 0 \quad \mbox{ for } \quad v \in V.\]
\end{prop}

\begin{prf}
Let $d := (0,0,1) \in \g$, then $U_t = \pi(0,0,t)=\exp(t\dd\pi(d))=\exp(-itH)$.
Let $\xi \in \cH^\infty$ be a smooth vector of $\pi$, let  $v \in V$, then by assumption
we have for every $t \in \R$ the inequality
\[ 0
\leq \la \pi(\exp tv) H \pi(\exp -tv)\xi, \xi \ra=
 i \la \pi(\exp tv) \dd\pi(d) \pi(\exp -tv)\xi, \xi \ra
=  i \la \dd\pi(e^{t \ad v}d)\xi, \xi \ra.\]
Now
\begin{eqnarray*}
(\ad v)d &=& (0, - (\ad d)v,0)= -\frac{d}{dt}(0,\tau_t(v),0)\Big|_{t=0}=(0,-Yv,0), \\[1mm]
(\ad v)^2d &=& -[v, Yv]=(\sigma(v,Yv),0,0),
\end{eqnarray*}
hence
\[  e^{t \ad v}d = d + t [v,d] + \frac{t^2}{2} [v,[v,d]]
= \Big(\frac{t^2}{2} \sigma(v,Yv), - t Yv, 1\Big) \]
and so
\[ 0\leq  i \la \dd\pi(e^{t \ad v}d)\xi, \xi \ra
= \frac{t^2}{2} \sigma(Yv, v) \la \xi, \xi \ra
 -  i t\la \dd\pi(Yv)\xi,\xi\ra+ i \la \dd\pi(d)\xi, \xi \ra.\]
Since this holds for all $t \in \R$, we  obtain
$\sigma(Yv,v) \geq 0$.
\end{prf}

In the special case that $V$ is a complex pre-Hilbert space $\cD$,
$\sigma(v,w) = \Im \la v, w\ra$ and
$\tau_t \in \cU(\cD)$, then $\la Yv,v \ra \in i\R$, so that
\[ 0\leq\sigma(Yv,v) = \Im \la Yv, v \ra = -i \la Yv,v \ra
= \la -i Yv,v\ra \]
implies that the infinitesimal generator $- iY$ of the unitary
one-parameter group $(\tau_t)_{t \in \R}$ is non-negative if
there exists a covariant representation with positive spectrum for $(\cA, G, \alpha)$.
In this case, as Fock representations exist and the second quantization of
a positive operator is positive, we also have the converse implication
(cf.\ \cite{NZ13}, \cite{Ze13}).

\begin{ex}
\label{cocycleXmp}
 Let $G = \cU_2(\cH) := \cU(\cH) \cap (\1 + \cB_2(\cH))$, where
$\cB_2(\cH)$ is the ideal of Hilbert--Schmidt operators.
Then $G$ is a Banach--Lie group with Lie algebra
$\g = \fu_2(\cH) = \{ X \in \cB_2(\cH) \,\mid\, X^* = - X\}$.
It is an interesting problem to determine
all projective unitary representations of~$G$.
That this problem is naturally linked to covariant
representations is due to the fact that
every continuous cocycle $\omega \: \g \times \g\to \R$ is of the form
\[ \omega(X,Y) = \tr([D,X]Y) = \tr(D[X,Y]) \]
for some $D \in \fu(\cH)$ (see \cite[Prop.~III.19]{Ne03} and its proof).
Then
$\alpha_t(g) :=\exp(tD) g \exp(-tD)$ is a continuous one-parameter group
of automorphisms of $G$ acting naturally on the central extension
$\R \oplus_\omega \g$ with the bracket
$[(z,x),(z',x')] := (\omega(x,x'), [x,x'])$ by
$\alpha_t(z,x) := (z, \alpha_t(x))$ and this
action lifts to the corresponding simply connected group
$\hat G$, which leads to a Lie group
$G^\sharp := \hat G \rtimes_\alpha \R$.
Its Lie algebra is the double extension
\[ \g^\sharp = \R \oplus \g \oplus \R, \quad
[(z,X,t), (z',X',t')] = (\omega(X,X'), [X,X'] + t [D,X'] - t'[D,X], 0).\]
Presently, the classification of all corresponding
projective covariant representations with positive spectrum is still open. However,
the case where $D$ is diagonalizable and the representation
is a highest weight representation has been treated fully in \cite{MN16};
see also \cite{Ne17} for more complete results.

Since projective covariant representations with positive spectrum of $G$ lead to
unitary representations $(U,\cH)$
of the corresponding doubly extended group $G^\sharp$
for which the convex cone
\[ W := \{ x \in \g^\sharp \: -i\dd U(x) \geq 0\} \]
has interior points, the method developed in \cite{NSZ17} provides a natural
$C^*$-algebra whose representation corresponds to these
representations of $G^\sharp$. From the perspective of Remark~\ref{rem:twist},
these representations correspond as well to
covariant representations with positive spectrum of $(\cA,\R,\alpha)$
for $\cA = C^*_\gamma(\cU_2(\cH)_d)$, where this denotes the twisted group
algebra corresponding to a central extension $\hat G$
of $\cU_2(\cH)$ by $\T$ corresponding to the Lie algebra extension
defined by the cocycle $\omega$
(see also \cite{Ne14}).
\end{ex}

The Borchers--Arveson Theorem
 also produces obstructions for various actions of groups on $C^*$-algebras, as in
the following framework:
\begin{itemize}
\item{} There is a unital $C^*$-algebra $\cA$, and two actions
$\alpha \: \R \to \Aut(\cA)$,  $\beta\: G\to \Aut(\cA)$ for a topological group $G$
and a nontrivial group action $\gamma:\R\to \Aut(G)$ which intertwines $\alpha$ and
$\beta$, i.e. $\beta(\gamma_t(g))=\alpha_t\circ \beta(g)\circ \alpha_{-t}$ for all $t\in\R$, $g\in G$.
\item{} Given this setting, then a covariant representation is a triple
 $(\pi, U, V)$, where $\pi \: \cA \to \cB(\cH)$ is a nondegenerate representation,
$U \: \R \to \cU(\cH)$ is a unitary one-parameter group, and
$V \: G \to \cU(\cH)$ is a continuous unitary representation such that
\[ U_t \pi(A) U_{-t} = \pi(\alpha_t(g)),
\quad V_g \pi(A) V_{g^{-1}} = \pi(\beta_g(A)),
\quad  U_t V_g U_{-t} =V_{\gamma_t(g)}
\]
for all $A\in\cA,$ $g \in G$ and  $t \in \R.$ We will say it
{\it has positive spectrum} if $U$ has positive spectrum.
\end{itemize}
This framework will occur for example if one tries to quantize Lagrangian
classical gauge theory
on Minkowski space (cf.~\cite{Ble81}).
In such a quantum gauge theory, $\cA$ will be the algebra of  observables,
$\alpha$ is time evolution, and $\beta$ gives the gauge transformations.
As the base space of the gauge theory is Minkowski space,
$G$ can be matrix-valued functions on the base space, and $\gamma$ will consist of translations along
the time coordinate.

The important action in this setting which will prohibit covariant representations, is
$\gamma:\R\to \Aut(G)$. We give a class of relevant examples where no nontrivial
covariant representations with positive spectrum are possible.
\begin{prop} \mlabel{prop:nogo0}
Let $X$ be a locally compact Hausdorff space, and let $F\subseteq U(n)$ be a closed subgroup containing $\T\1$,
and let $G\subset C_b(X,F)$ be a subgroup with respect to  pointwise multiplication.
Let $t\mapsto\varphi_t\in{\rm Homeo}(X)$ be a one-parameter group of homeomorphisms, and assume
 that $g\circ\varphi_t\in G$ for all $t\in\R$ and $g\in G$. Consider the action
\[\gamma:\R\to \Aut(G), \qquad
\gamma_t(g)(x):=\lambda_t(x)(g(\varphi_{-t}(x))), \quad
\mbox{ where } \quad  \lambda_t(x)\in\Aut(F),\]
so that   $\lambda_t(x)$ fixes  $\T\1$ pointwise.
 If $(V, U)$ is a  covariant representation with positive spectrum of $\gamma$,
 then for any $g\in C_b(X,\T\1)\cap G$ we have $V_g=V_{\gamma_t(g)}$ for all $t\in\R$.

Assume as above, the two actions $\alpha \: \R \to \Aut(\cA)$,  $\beta\: G\to \Aut(\cA)$ for  $G\subset C_b(X,F)$
where $\cA$ is a simple unital $C^*$-algebra, and that
$\gamma:\R\to \Aut(G)$ above intertwines $\alpha$ and
$\beta.$
 Assume that $\alpha_\R$ and $\beta_g$ do
 not commute in $\Aut\cA$ for some $g\in C_b(X,\T\1)\cap G$.
 Then the only covariant representation with positive spectrum is
 the zero representation.
\end{prop}

\begin{prf}
Assuming a covariant representation with positive spectrum $(V, U)$ of $\gamma$, note
that $(V_{\gamma_t(g)})_{t\in\R}$
is commutative if $g\in C_b(X,\T\1)\cap G$. Thus $\cN:=(V_{\gamma_\R(g)})''\subset\cB(\cH)$ is commutative,
and by construction the action of $\Ad(U_t)=:\tilde\alpha_t$ will
preserve~$\cN$ (Remark~\ref{RemBorAr}(b)).
Thus by the Borchers--Arveson
theorem $\cN$ contains the minimal unitary implementers for $\tilde\alpha_t$
which therefore commutes with $V_g\in\cN$ and so by covariance $V_g=V_{\gamma_t(g)}$ for all $t\in\R$.

For the second part, let  $(\pi, U, V)$ be a covariant representation with positive spectrum. By the previous part we have
that  $V_g=V_{\gamma_t(g)}$ for all $t\in\R$ and for $g\in C_b(X,\T\1)\cap G$. Now
\[
\pi(\beta_g(A))=\Ad(V_g)\pi(A)=\Ad(V_{\gamma_t(g)})\pi(A)=\pi(\beta_{\gamma_t(g)}(A))
=\pi(\alpha_t\circ\beta_g\circ\alpha_{-t}(A)),
\]
hence $\alpha_t\circ\beta_g(A)-\beta_g\circ\alpha_t(A)\in\Ker\pi$ for all $A\in\cA$.
By hypothesis there is an $A\in\cA$ and $t\in \R$ for which this is nonzero, hence as
$\cA$ is simple, $\pi$ must be the zero representation.
\end{prf}

\begin{rem}
One way to circumvent the obstruction from Proposition~\ref{prop:nogo0}, is to
ask instead for a covariant representation with positive spectrum $(\pi, U, V)$, where
$V \: G \to \cU(\cH)$ is a continuous {\it projective} unitary representation.
It is interesting that even in the Hamiltonian approach to quantum gauge theory
(where $\gamma$ is trivial),
 projective gauge transformations occur naturally. These are obtained e.g. by using
a quasi-free Fock representation of the CAR-algebra to produce an implementing unitary group with positive spectrum for
the time evolutions (cf. \cite{CR87,La94}).

In this context we also mention that the method to relate
covariant representations with positive spectrum
to positivity of a Lie algebra cocycle that we have seen
in Example~\ref{WeylPos} has been put to work extensively
in the context of  covariant representations with positive spectrum for
gauge groups corresponding to semi-simple structure groups in \cite{JN17}.
\end{rem}

Given the obstruction in Proposition~\ref{prop:nogo0}, one strategy is to weaken the requirements
on the representation. Starting with the actions
$\alpha \: \R \to \Aut(\cA)$,  $\beta\: G\to \Aut(\cA)$
and  $\gamma:\R\to \Aut(G)$ such that
 $\beta(\gamma_t(g))=\alpha_t\circ \beta(g)\circ \alpha_{-t}$ for all $t\in\R$, $g\in G$,
 one considers triples $(\pi, U, V)$, where $\pi \: \cA \to \cB(\cH)$ is a nondegenerate representation,
$U \: \R \to \cU(\cH)$ is a  unitary one-parameter group with positive spectrum, and
$V \: G \to \cU(\cH)$ is a map (not necessarily a representation) such that
\[ U_t \pi(A) U_{-t} = \pi(\alpha_t(g)),
\quad V_g \pi(A) V_{g^{-1}} = \pi(\beta_g(A)),
\quad  U_t V_g U_{-t} =V_{\gamma_t(g)}\]
Then it follows that $V \: G \to \cU(\cH)$  must be a cocycle representation, i.e.
\[
V_gV_h=\mu(g,h)V_{gh}\quad\hbox{where}\quad \mu(g,h)\in\pi(\cA)'\cap (U_{\R})'
\]
for $g,\,h\in G$.
By Proposition~\ref{prop:nogo0}, we know that the cocycle $\mu$
must be nontrivial.

\section{Ground states and their covariant representations}
\mlabel{sec:5}

Recall from Lemma~\ref{lem:minimal} that if
a unitary one-parameter subgroup with positive spectrum $(U_t)_{t \in \R}\subset\cM$  is minimal,
then for every $\eps > 0$, the central support of $P[0,\eps)$ is $\1$. Below,
in a suitable subrepresentation, we will find
a similar property for $P(\{0\})$,  the projection onto the space of invariant vectors.
 In the next two theorems, we  first investigate structures
associated with projections of central support $\1$.

\subsection{Ground states of a $C^*$-action  $(\cA,\R,\alpha)$}
\mlabel{subsec:5.1}

\begin{defn} \mlabel{def:groundvec}
Let $(\cM,\R, \alpha)$ be a concrete $W^*$-dynamical system on $\cH$, i.e.\
$\cM \subeq \cB(\cH)$ is a von Neumann algebra.
The {\it ground state vectors}
of a covariant representation with positive spectrum are the $U$-invariant
elements of $\cH$ with respect to \
the \minimal one-parameter group from the Borchers--Arveson Theorem~\ref{BA-thm}.
(This should be distinguished
from the ground states defined in Definition~\ref{defgroundst0} below; but see
Corollary~\ref{groundchar}).
\end{defn}

In the
physics literature, the ground state vectors are defined as the eigenvectors of the Hamiltonian
corresponding to the lowest value of its spectrum. As is well-known, for
e.g.~the quantum oscillator
in the Schr\"odinger representation, this lowest spectral value can be nonzero. However, this
definition coincides with our
definition, as we took the \minimal one-parameter group, and for this, the lowest spectral value of its generator is zero. In the example of the quantum oscillator, the generator of the minimal group
is the usual Hamiltonian plus the multiple of the identity needed to shift the lowest value of its spectrum to zero.

\begin{thm} \mlabel{thm:8.7} Let $\cM$ be a von Neumann algebra
and $(U_t)_{t \in \R}\subset\cM$ be a one-parameter group with positive spectrum in $\cM$
which is minimal.
We write $P$ for the spectral measure of $U$ for which $U_t = \int_\R e^{-itp}\, dP(p)$
and put $P_\eps := P[0,\eps]$ for $\eps \geq 0$. Let
$Z_0$ be the central support of $P_0$. Then we obtain a direct sum decomposition
\[ \cM = Z_0 \cM \oplus (\1 - Z_0) \cM.\]
Moreover, the following assertions hold:
\begin{itemize}
 \item[\rm(i)] For all normal representations $(\pi, \cH)$ of the ideal $Z_0 \cM$,
the subspace $\cH_0 := \pi(P_0)\cH$ of ground state vectors is $\pi(\cM)\hbox{-generating.}$
\item[\rm(ii)] All normal representations $(\pi, \cH)$ of $(\1-Z_0) \cM$ are  covariant
representations with positive spectrum with respect to $(\pi(U_t))_{t \in \R}$,
but they contain no non-zero ground state vectors.
\item[\rm(iii)] For all normal representations $(\pi, \cH)$ of $\cM$ and
$\eps > 0$, the subspace $\pi(P_\eps)\cH$ is $\pi(\cM)\hbox{-generating.}$
\end{itemize}
\end{thm}

\begin{prf} (i) Let $(\pi, \cH)$ be a normal representation of $Z_0\cM$, which corresponds to a normal
representation of $\cM$ with $\pi(Z_0) = \1$. Then the central support of
$\pi(P_0)$ is $\1$, so that the assertion follows from Lemma~\ref{CP1cyclic}.

(ii) As $(U_t)_{t \in \R}$ has positive spectrum,
 the one-parameter group $(\pi(U_t))_{t \in \R}$
has positive spectrum.
If $\pi(Z_0) = 0$, then also $\pi(P_0) = 0$, so that $\pi(U_t)$ has no non-zero fixed vectors.
The minimality of $U$ implies that $\pi \circ U$ is minimal in $\pi(\cM)$
(Lemma~\ref{mintomin}),
so that $\inf \Spec(\pi \circ U) = 0$. Hence there are no ground state vectors for $\alpha$
in $\cH$.

(iii) follows immediately from $\inf \Spec(\pi \circ U) = 0$ in every normal representation
$\pi$ of~$\cM$, Lemmas~\ref{CP1cyclic} and ~\ref{lem:minimal}.
\end{prf}

\begin{rem} Suppose, in the context of Theorem~\ref{thm:8.7}, that
$0$ is isolated in the spectrum of the implementing unitary group $U$ with positive spectrum.
Then the central support of $P_0$ is $\1$, hence $\cM = Z_0 \cM$.
This is clearly an important subcase, which we will analyze in detail in Subsect.~\ref{Isolated0}
below.
\end{rem}

\begin{prop} \mlabel{prop:swallow} Let $(\cA,\R,\alpha)$ be a $C^*$-action and
let $(\pi,U)$ be a  covariant representation with positive spectrum
for which
the subspace $\cH_0$ of $U$-fixed vectors
is generating. Then $(U_t)_{t \in \R}$ is
the Borchers--Arveson minimal group, hence in particular
$U_\R \subeq \pi(\cA)''$.
\end{prop}

\begin{prf} Let $\cM:= \pi(\cA)''$. From the Borchers--Arveson Theorem~\ref{BA-thm}, we obtain the
unique \minimal  one-parameter
group $(V_t)_{t\in \R}$ in $\cM$ implementing the automorphisms $\Ad(U_t)$.
Then $W_t := U_t V_t^*\in \cM'$ is a one-parameter group with positive spectrum
(Lemma~\ref{lem:minimal}).
Let $H$, $H_1$ and $H_2$ denote the infinitesimal generators  of $U$, $V$ and $W$,
respectively. All these operators have non-negative spectrum, so that
Lemma~\ref{lem:UVW} implies that $H = H_1 + H_2$. Therefore
$\cH_0 \subeq \cD(H) = \cD(H_1) \cap \cD(H_2)$ and, for every
$\Omega \in \cH_0$, we have
\[ 0 = \la H \Omega, \Omega \ra
=  \la H_1 \Omega, \Omega \ra + \la H_2 \Omega, \Omega \ra. \]
This implies $H_2 \Omega =0$, so that $\Omega$ is fixed by $W$.
As $\cH_0$ is $\cM$-generating, it is separating for $\cM'$,
 which leads to $W_t = \1$ for $t \in \R$. This proves that $U_t = V_t
\in \pi(\cA)''$.
\end{prf}

Recall for an invariant state $\omega$, the
GNS unitary group $U^\omega$ from above (preceding Proposition~\ref{fixstatecov}).
We define:

\begin{defn}\mlabel{defgroundst0}
Given  a $C^*$-action $(\cA,\R,\alpha),$ then a {\it ground state} is
 an invariant state $\omega\in\fS(\cA)$
for which its GNS unitary group $(U^\omega_t)_{t \in \R}$ is continuous and has positive spectrum
(cf. \cite[Def.~IV.4.9]{Bo96}).
Then $\Omega_\omega$ is a ground state vector
in the GNS representation by the next corollary.
\end{defn}

\begin{cor}
\mlabel{groundchar}
Assume a $C^*$-action $(\cA,\R,\alpha)$ and an invariant state $\omega
\in\fS(\cA)^\R$.
\begin{itemize}
\item[\rm(i)]
If $\omega$ is a ground state, i.e.\ $(U^\omega_t)_{t \in \R}$ is continuous and has positive spectrum, then $U^\omega$  is
the Borchers--Arveson minimal group, hence
$U^\omega_\R \subeq \pi_\omega(\cA)''$, and the GNS cyclic vector $\Omega_\omega$ is a
ground state vector for~$U^\omega$.
\item[\rm(ii)] If there is a
 Borchers--Arveson minimal group  $(V_t)_{t \in \R}$ on $\cH_\omega$
 and $\Omega_\omega$ is a ground state vector,
then $U^\omega$ has positive spectrum and coincides
with the Borchers--Arveson minimal group. Hence $\omega$ is a ground state.
\end{itemize}
\end{cor}

\begin{prf} (i) follows from $\Omega_\omega\in\cH_0$ and Proposition~\ref{prop:swallow}.

For (ii), by assumption we have $V_t\Omega_\omega=\Omega_\omega$ for all $t\in\R$.
 Together with covariance, this implies that $V_t=U^\omega_t$ for all~$t$, so that
by the definition $\omega$ is a ground state.
\end{prf}

In Subsection~\ref{SectGrndSt} below we will  study existence of ground states.

\begin{ex} A case of an invariant state $\omega$ for which $(\pi_\omega,U^\omega)$
does not have positive spectrum
but $\Spec(U^\omega)$ is bounded from below, so that there exists a positive implementation,
can be obtained as follows.

We consider $\cA = M_2(\C)$ with elements $A = (a_{ij})_{1 \leq i,j\leq 2}$, and let
$\alpha_t(A)=
\left({\hfill a_{11} \atop e^{it}a_{21}}{e^{-it}a_{12}\atop \hfill a_{22}} \right)
$.
Define the state $\omega$ by $\omega(A) = a_{11}$, which is a vector state invariant with respect to ~$\alpha_t$.
Then $U_t^\omega = \diag(1,e^{it})$, but  $U_t = \diag(e^{-it},1)$ also implements $\alpha_t$.
Then
$\Spec(U^\omega) = \{0,-1\}$ is not positive, and
$\Spec(U) = \{1,0\}$ is  positive.
\end{ex}

\begin{rem}
(a)
For the case where $(\cA,\R,\alpha)$ is a $C^*$-dynamical system, i.e.\ $\alpha$
is \strong continuous, then the analog of Corollary~\ref{groundchar} follows from \cite[Thm.~8.12.5]{Pe89}.

(b) The properties of ground states listed above in  Corollary~\ref{groundchar} are in the literature,
though with more restrictive assumptions than ours.
E.g in the \usual case, for a ground state $\omega$, we know from Araki \cite{Ar64} (cf.~\cite[Cor.~2.4.7]{Sa91})
that  $U^\omega_\R\subset\cM=\pi_\omega(\cA)''$.
If we do not have the \usual case,
but $\cA$ is assumed to have a local net structure as in \cite[Sect,~1.1]{Bo96},
then one obtains from \cite[Cor.~IV.4.11(2)]{Bo96}
that the GNS unitary group $U^\omega:\R\to\cU(\cH_\omega)$ of a ground state $\omega$
coincides with the minimal representation with positive spectrum $V:\R\to\cU(\cM)$. The main assumption
for a local net of observables is that $\cA$ is an inductive limit of ``local''
$C^*$-algebras $\cA(O)$ indexed by the bounded open sets $O$ in $\R^4$ such that $O_1\subset O_2$ implies
 $\cA(O_1)\subseteq \cA(O_2),$ and $\alpha$ is covariant with respect to time translations acting on the regions
 $O\subset\R^4$.

(c) In general, the projection onto a generating subspace
as in Proposition~\ref{prop:swallow} need not be contained in
$\pi(\cA)''$.
A typical example can be obtained for
$\cA = \cB(\cH) \oplus \cB(\cH)$ and the canonical representation on
$\cH \oplus \cH$. For any unit vector $v \in \cH$, the element
$(v,v) \in \cH \oplus \cH$ is cyclic, but the projection onto
$\C(v,v)$ is not contained in the von Neumann algebra $\cA$.
\end{rem}

\begin{prop}
\label{BijGrStRep}
Let $(\cM,\R, \alpha)$ be a $W^*$-dynamical system and let
$(U_t)_{t \in \R}\subset\cM$ be the unique \minimal one-parameter group
 such that $\alpha_t = \Ad(U_t)$ for $t \in \R.$
Given a normal representation $(\pi,\cH)$ of $\cM$
in which $\cH_{0} := \pi(P_0)\cH$ is generating,
i.e.\ which is generated by the ground state vectors,
construct the restricted
representation $(\pi_0, \cH_{0})$ of the reduction $\cM_{P_0}=P_0\cM P_0\subset\cM$,
i.e.\ $\pi_0(P_0MP_0):=\pi(P_0MP_0)\restriction\cH_0$, $M\in\cM$.
Then the map $\pi\to\pi_0$ is a bijection between
unitary isomorphism classes of normal
representations of $\cM$ generated by ground state vectors and
unitary isomorphism classes of normal
representations $(\pi_0, \cH_0)$ of the reduction~$\cM_{P_0}$.
\end{prop}

\begin{prf}
This is an application of Proposition~\ref{BijRedRep}.
\end{prf}

\begin{ex} Let $P$ be a projection in the $W^*$-algebra $\cM$ and consider the corresponding
unitary one-parameter group
\[ U_t := P + e^{-it}(\1-P) = e^{-itH} \quad \mbox{ for } \quad H = \1-P.\]
We assume that the central support of $P$ is $\1$, so that $(U_t)_{t \in \R}$
is minimal (Lemma~\ref{lem:mini}).
For any normal representation $(\pi, \cH)$ of $\cM$, the
subspace $\cH_P := \pi(P)\cH$ of ground states for~$U$ is generating. It carries a
representation of the ideal $\cM_P = P \cM P$ of $\cM_0$ which determines it uniquely.
\end{ex}

\subsection{Existence of ground states}
\label{SectGrndSt}

Recall that above in Definition~\ref{defgroundst0} we defined a ground state for a
given   $C^*$-action $(\cA,G,\alpha),$ as
 an invariant state $\omega\in\fS(\cA)$ for which
$(U^\omega_t)_{t \in \R}$ is continuous and has positive spectrum.
In this case $\Omega_\omega$ is a ground state vector for $(\pi_\omega,U^\omega)$.
Denote the set of ground states by $\fS_\alpha^0(\cA)$.

\begin{rem}
\label{DopDef}
In the \usual case of a
given   $C^*$-action $(\cA,G,\alpha), $(i.e.\ $\alpha$ is \strong continuous), the
left ideal
 $\cL=\br \cA\cdot\cA^\alpha_0(-\infty,0). $
generated by the subspace $\cA^\alpha_0(-\infty,0)$ (cf.\ Definition~\ref{def:arveson})
 selects the ground states by $\omega(\cL)=\{0\}$,
 i.e.\ $\omega \cA^\alpha_0(-\infty,0)=0$. The left ideal
$\cL$ is the well-known {\it Doplicher  ideal}
used for algebraic characterization of a ground state (cf.~\cite{Dop65}),
and leads to an alternative definition of a ground state
 (cf.~\cite[Def.~4.3, p.~82]{Ar99} and \cite[Prop.~5.3.19]{BR96}).
Then  $\fS_\alpha^0(\cA)\not=\emptyset$ if and only if $\cL$ is proper in $\cA$.
\end{rem}

 In our case, we need not have
 that $\alpha$ is \strong continuous,
  hence we need to deal with $\cM^\alpha_0(-\infty,0)\subset\cM=\pi_\omega(\cA)''$,
hence the condition $\cM^\alpha_0(-\infty,0)\Omega_\omega=0$ is external to $\pi_\omega(\cA)$.
We first make our condition explicit in the next proposition.

 \begin{lem}\mlabel{lem:3.58}
 Given  a $C^*$-action $(\cA,\R,\alpha),$ consider the
associated $W^*$-dynamical system \break
$\alpha^{co}:\R\to \Aut(\cM_{co})$, where $\cM_{co}:=\pi_{co}(\cA)''$ and
$\alpha^{co}(t)=\Ad U_{co}(t)$. Then $\omega\in\fS(\cA)$ is a ground state for
$(\cA,\R,\alpha)$ if and only if it has a normal extension to $\cM_{co}$ which
is a ground state for $(\cM_{co},\R,\alpha^{co})$.
  \end{lem}

\begin{prf}
Let $\omega\in\fS(\cA)$ be a ground state of $(\cA,\R,\alpha)$. Then the
GNS covariant representation ${(\pi_\omega, U^\omega,\Omega_\omega)}$
extends to a cyclic representation of $\cA \rtimes_\alpha \R_d\supset\cA$ for which $(\pi_\omega,U^\omega)\in{\rm Rep}(\alpha,\cH_\omega)$.
Thus $(\pi_\omega,U^\omega)$ is a subrepresentation of $(\pi_{co},U_{co})$ hence $\omega$ has the  normal extension
$\tilde\omega(M):=(\Omega_\omega,\pi_{co}(M)\Omega_\omega)$ for $M\in\cM_{co}.$
It is clear that it is invariant with respect to  $\alpha^{co}=\Ad U_{co}(t)$ as $U_{co}(t)$ coincides on $\cH_\omega$
with $U^\omega_t$.
To see that $\tilde\omega$ is a ground state for $(\cM_{co},\R,\alpha^{co})$, note first that the GNS representation
$(\pi_{\tilde\omega},\Omega_{\tilde\omega})$ is just the restriction of $\pi_{co}$ to $\cH_\omega$,
where
$\Omega_{\tilde\omega}=\Omega_{\omega}$ on which we have
$\pi_{\tilde\omega}(\cM_{co})=\pi_\omega(\cA)''.$ It suffices to show that $U^{\tilde\omega}=U^\omega,$
as we know that $U^\omega$ has positive spectrum.
This follows from
\[
U^\omega_t M\Omega_\omega = (U^\omega_t MU^\omega_{-t})\Omega_\omega=\alpha_t^{co}(M)\Omega_{\tilde\omega}=U^{\tilde\omega}_tM\Omega_{\tilde\omega}
=U^{\tilde\omega}_tM\Omega_\omega
\]
 for all
$M\in\pi_\omega(\cA)''.$

Conversely, let  $\nu\in\fS(\cM_{co})$ be a normal ground state for $(\cM_{co},\R,\alpha^{co})$.
Then $\nu\restriction \pi_{co}(\cA)$ is an invariant state for $\alpha$. As $\nu$ is normal, it follows that
$\pi_\nu(\pi_{co}(\cA))$ is strong operator dense in $\pi_\nu(\pi_{co}(\cM_{co})),$ hence $\Omega_\nu$ is cyclic with respect to  both algebras, and $(\pi_\nu\circ \pi_{co})\restriction\cA= \pi_{(\nu\restriction\cA)}$.
  Furthermore
\[ U^\nu_t (\pi_\nu\circ\pi_{co})(A)\Omega_\nu
= (U^\nu_t (\pi_\nu\circ\pi_{co})(A)U^\nu_{-t})\Omega_\nu
=\pi_\nu\big(\pi_{co}(\alpha_t(A))\big)\Omega_\nu \quad \mbox{ for } \quad A\in\cA,\]
 hence $U^\nu$ is the GNS implementing unitary group for both $\nu$ and $\nu\restriction\cA$, and it is clear that
it has positive spectrum and leaves $\Omega_\nu$ invariant. Thus $\nu\restriction\cA$ is a ground state of $(\cA,G,\alpha).$
\end{prf}
By the preceding lemma, the next
proposition also covers $C^*$-actions $(\cA,G,\alpha)$.
\begin{prop}
\label{GroundCond}
 Let $\cM \subeq \cB(\cH)$ be a von Neumann algebra
and $(\cM,\R,\alpha)$ a $W^*$-dynamical system. Then the following conditions
are equivalent for a normal state $\omega$ of $\cM$:
\begin{itemize}
\item[\rm(i)] $\omega$ is a ground state.
\item[\rm(ii)] $\omega\cM^\alpha(-\infty,0) = \{0\}$.
\end{itemize}
If these conditions are satisfied, then the corresponding
GNS representation $(\pi_\omega, \cH_\omega)$ is covariant.
\end{prop}

\begin{prf}
Let ${(\pi_\omega,\cH_\omega,\Omega_\omega)}$ be the  GNS representation of a given normal state $\omega$ and $\cN := \pi_\omega(\cM)$.\\[1mm]
(i) $\Rarrow$ (ii): Let  $\omega$ be a normal ground state of $(\cM,\R,\alpha)$ and let
 $(U^\omega_t)_{t \in \R}$ be the \minimal
one-parameter group  implementing
$\beta_t:=\pi_\omega\circ\alpha_t = \Ad(U^\omega_t)$ for $t \in \R$ on $\cM.$
Then
$\Omega_\omega \in \cH_\omega^{U^\omega}(\{0\})$, so that
\[ \cN^\beta(-\infty,0)\Omega_\omega \subeq \cH^{U^\omega}(-\infty,0) = \{0\} \]
follows from $\Spec(U^\omega) \subeq [0,\infty)$.
For $M \in \cM$ and $f\in L^1(\R)$ with  $\supp(\hat{f})\subseteq(-\infty,0)$,
we have
\[ \beta_f(\pi_\omega(M)):=\int f(t)\beta_t(\pi_\omega(M))\, dt=\pi_\omega(\alpha_f(M))\in
\cN^\beta(-\infty,0). \]
As $(-\infty,0)$ is open,  it follows by
\cite[Lemma~3.2.39(4)]{BR02} that all elements in $\cM^\alpha(-\infty,0)$ are  $\sigma(\cM,\cM_*)$-limits
of such $\alpha_f(M)$. As $\pi_\omega$ is normal, we thus have
\[ \pi_\omega(\cM^\alpha(-\infty,0)) \subeq \cN^\beta(-\infty,0) \]
from which it follows by the first part that $\omega\cM^\alpha(-\infty,0) = \{0\}$.

(ii) $\Rarrow$ (i):
Assume that $\omega\cM^\alpha(-\infty,0) = \{0\}$, then we first prove that $\omega$
is $\alpha\hbox{--invariant}$ (using a short argument from \cite{Ped79}). As $\1\in\cM$
we see that $\omega(\cM^\alpha(-\infty,0))=\{0\}$, and this yields
\[
\omega(\cM^\alpha(0,\infty))=\overline{\omega(\cM^\alpha(0,\infty)^*)}
=\overline{\omega(\cM^\alpha(-\infty,0))}=\{0\}
\]
i.e.\ $\cM^\alpha(-\infty,0)\cup\cM^\alpha(0,\infty)\subset\ker\omega$. However $\cM$ is  the
$\sigma(\cM,\cM_*)$-closure of the span of the
$\alpha\hbox{--preserved}$ spaces $\cM^\alpha(-\infty,0),$ $\cM^\alpha(0,\infty)$
and $\cM^\alpha\{0\}$, and these spaces only intersect in zero. Thus $\omega$ is only nonzero
on $\cM^\alpha\{0\}$ and as the action of $\alpha$ on this space is trivial, it follows that
$\omega$ is $\alpha$-invariant. Thus the GNS unitary group $U^\omega$ implements
$\alpha$ in $\pi_\omega$ (but at this point we do not know that it is continuous).
Then $\beta_t:={\rm Ad}(U^\omega_t)$ defines a  $W^*$-dynamical system
 $(\pi_\omega(\cM),\R,\beta)$, because $\beta_t(\pi_\omega(M))=\pi_\omega(\alpha_t(M))$ and the right hand side
 is a composition of the $\sigma(\cM,\cM_*)$-continuous map $t\mapsto \alpha_t(M)$ with the normal map
 $\pi_\omega$.

Now by \cite[Lemma~3.2.39(4)]{BR02} it also follows for the  $W^*$-dynamical system
$(\cN,\R,\beta)$ where $\cN = \pi_\omega(\cM)$ and  $\beta_t \circ \pi_\omega
=\pi_\omega\circ\alpha_t$, that
all elements in $\cN^\beta(-\infty,0)$ are  $\sigma(\cN,\cN_*)$-limits
of elements $\beta_f(\pi_\omega(M))=\pi_\omega(\alpha_f(M))$ for $M\in\cM$ and
 $f\in L^1(\R)$ such that $\supp(\hat{f})\subseteq(-\infty,0)$.
Thus we get from our assumption that $\cN^\beta(-\infty,0)\Omega_\omega = \{0\}$.
For every
$s > 0$ we have
\[ \Omega_\omega \in \ker(\cN^\beta(-\infty,-s]) =
\ker(\cN^\beta[s,\infty)^*)
= \lbr\cN^\beta[s,\infty)\cH_\omega\rbr^\bot\]
(for the first equality, see \cite[Lemma~3.2.42(ii)]{BR02}).
Thus
\[
\Omega_\omega\perp\bigcup_{t>0}\br\cN^\beta[t,\infty)\cH_\omega.
\supset\bigcap_{t>0}\br\cN^\beta[t,\infty)\cH_\omega.
\supset\bigcap_{t\in\R}\br\cN^\beta[t,\infty)\cH_\omega. .\]
Hence
\[ \Omega_\omega\in \Big(\bigcap_{t\in\R}\br\cN^\beta[t,\infty)\cH_\omega.\Big)^\perp.\]
The closed space on the
right hand side is $\cN$-invariant,
so, as it contains a cyclic vector, it must be all of $\cH_\omega$. Thus
\[
\bigcap_{t\in\R}\br\cN^\beta[t,\infty)\cH_\omega.=\{0\}
\]
and so we may apply the Borchers--Arveson Theorem~\ref{BA-thm}
to conclude that
the minimal Borchers--Arveson subgroup $(U_t)_{t \in \R}\subset\cN$
exists and that its spectral measure $P$ is given by
\[
P[t,\infty)\cH_\omega = \bigcap_{s < t} \br\cN^\beta[s,\infty)\cH_\omega..
\]
As $P(\R) = P[0,\infty)$, it follows
that $\Omega_\omega \in P(\{0\})$, so that $\Omega_\omega$ is  $U$-invariant,
hence $\omega$ is a ground state for $(\cM,\R,\alpha)$, and $U$ coincides with $U^\omega$
by Corollary~\ref{groundchar}(ii).
Now the covariance of $\pi_\omega$ follows from  Proposition~\ref{fixstatecov}.
\end{prf}

In analogy to the Doplicher existence criterion for ground states in the \usual case, we then have:
\begin{cor}\mlabel{groundExist}
\begin{itemize}
\item[\rm(i)] Let
$(\cM,\R,\alpha)$ be a $W^*$-dynamical system. Then  a normal ground state  of $\cM$
exists if and only if the $\sigma(\cM,\cM_*)$-closed left ideal generated in $\cM$ by $\cM^\alpha(-\infty,0)$
is not all of $\cM$.
\item[\rm(ii)] For a given   $C^*$-action $(\cA,\R,\alpha),$ a ground state exists if and only if
for the associated $W^*$-dynamical system $(\cM_{co},\R,\alpha^{co})$,
the $\sigma(\cM_{co},(\cM_{co})_*)$-closed left ideal generated in $\cM_{co}$ by $\cM_{co}^{\alpha^{co}}(-\infty,0)$
is not all of $\cM_{co}$.
\end{itemize}
\end{cor}

\begin{prf} As (ii) is obvious, we only prove (i).
By  Proposition~\ref{GroundCond} the ground states are precisely the states in the annihilator
in  $\cM_*$ of the $\sigma(\cM,\cM_*)$-closed left ideal  in $\cM$ generated by $\cM^\alpha(-\infty,0)$.
As the predual $\cM_*$ separates the $\sigma(\cM,\cM_*)$-closed left ideals
of $\cM$ by \cite[Thm~3.6.11, Prop.~2.5.4]{Pe89}, we conclude that the annihilator in  $\cM_*$
of a $\sigma(\cM,\cM_*)$-closed left ideal in $\cM$
is nonzero if and only if this left ideal is not all of $\cM$.
\end{prf}
If the  $C^*$-action $(\cA,\R,\alpha),$ is not a $C^*$-dynamical system, it seems very difficult to
obtain a similar internal criterion on $\cA$ alone for the existence of ground states.

\begin{rem} {\rm(Weak clustering)}
If $\omega$ is a ground state of a  $C^*$-action $(\cA,\R,\alpha),$
 then the question arises whether its ground state vectors in its GNS representation
are unique (up to multiples) or not. Let $P_0$ be the projection onto the fixed points of $U^\omega$, so $\Omega_\omega\in P_0\cH_\omega$.
If ${\rm dim}(P_0\cH_\omega)=1$
(the ground state vector is unique) then $\pi_\omega$ is irreducible
(cf.~\cite[Prop.~2.4.9]{Sa91}). By Theorem~\ref{thm:b.3} below, this will be the case if
$(\pi_\omega(\cA)'')^{\R} = \C \1$.

Otherwise, if ${\rm dim}(P_0\cH_\omega)>1$ and $\cM_{P_0} := P_0\cM P_0$ is abelian,
 then $\cM$ is type~I (cf.~\cite[Prop.~2.4.11]{Sa91}).
The condition that $\cM_{P_0}$ is abelian  will
be guaranteed in a local net of $C^*$-algebras
as for the Haag--Kastler axioms (cf.~\cite[Prop.~3]{Ar64}).

Recall from Proposition~\ref{posrep1} the definition of $\fS_{co}^{(+)}(\cA)$.
If one assumes the Haag--Kastler axioms, then all states in
$\fS_{co}^{(+)}(\cA)$ are
ground states (cf.~\cite[Thm~IV.4.10]{Bo96})
and for these the GNS unitary group $U^\omega:\R\to\cU(\cH_\omega)$
coincides with unique \minimal representation  $V:\R\to\cU(\cM)$ (cf.~\cite[Cor.~IV.4.11]{Bo96}).
\end{rem}

\subsection{The case where $0$ is isolated in $\Spec_\alpha(\cM)$}
 \label{Isolated0}

We now take a closer look at ground states of a given
 $W^*$-dynamical system $(\cM,\R,\alpha)$
under the
assumption that $0$ is isolated in $\Spec_\alpha(\cM)\subeq\R$
(this includes the case of $\T$-actions). In the physics literature this
is discussed as the ``spectral gap,'' and this is well-studied, e.g.\ in lattice systems
\cite{HK06}, or the mass gap in quantum field theory \cite[Sec.~4.4]{Ar99}.
We assume that there exists an $\eps > 0$ such
\begin{equation}
  \label{eq:gap}
 \Spec_\alpha(\cM) \cap [-\eps,\eps] = \{0\}.
\end{equation}
Accordingly, we write
\[ \cM_0 := \{ M \in \cM \mid (\forall t \in \R)\,
\alpha_t(M) = M\} = \cM^\alpha(\{0\}), \]
\[  \cM_+ := \cM^\alpha(0,\infty) = \cM^\alpha[\eps, \infty) \quad \mbox{ and } \quad
\cM_- :=
\cM^\alpha(-\infty,0) = \cM^\alpha(-\infty,-\eps].\]
These are weakly closed subalgebras with
$\{ M^* \mid M \in \cM_\pm\} = \cM_\mp$.
For any $f \in L^1(\R)$ with $\supp(\hat f) \subeq (-\eps,\eps)$ and $\hat f(0) =1$,
we then have $\alpha_f (\cM_\pm)= \{0\}$ and $\alpha_f(M) = M$ for $M \in \cM_0$, so that
this element defines a weakly continuous projection
\[ p_0 = \alpha_f \: \cM \to \cM_0\quad \mbox{ with } \quad
\ker p_0 \supeq \cM_+ + \cM_-.\]
Further, any $f \in \cS(\R)$ can be written as a sum of three Schwartz functions
$f = f_- + f_0 + f_+$ with
\[ \supp(\hat f_0) \subeq (-\eps,\eps), \quad
\supp(\hat f_-) \subeq (-\infty, -\eps/2)\quad \mbox{ and } \quad
\supp(\hat f_+) \subeq (\eps/2, \infty).\]
Then $\alpha_f = \alpha_{f_+} + \alpha_{f_0} + \alpha_{f_-}$ with
$\alpha_{f_\pm}(\cM) \subeq \cM_\pm$, so that
$\cM_- + \cM_0 + \cM_+$ is weakly dense in $\cM$, resp.,
$\cM_- +  \cM_+$ is weakly dense in $\ker p_0$.
In general we cannot expect that $\cM = \cM_- + \cM_0 + \cM_+$, as
the example
$\cM = \cB(\ell^2(\N))$ and $\alpha_t((M_{jk}) = (e^{it(j-k)}M_{jk})$ shows
(cf.~Example~\ref{ex:bel} below).
We also note that
\begin{equation}
  \label{eq:tridiag}
\cM^\alpha[0,\infty) = \cM_0 \oplus \cM_+  \quad \mbox{ and } \quad
\cM^\alpha(-\infty,0] = \cM_- \oplus \cM_0.
\end{equation}

\begin{rem}
If $(\cA,\R,\alpha)$ is a $C^*$-dynamical system for which
$0$ is isolated in $\Spec_\alpha(\cA)$, then
we would like to have a direct decomposition into three closed subalgebras
\begin{equation}
  \label{eq:split-con}
\cA = \cA_- \oplus \cA_0 \oplus \cA_+,
\end{equation}
defined by the spectral projections corresponding to
$(-\infty,0)$, $\{0\}$ and $(0,\infty)$.
Such a decomposition  always exists if $\alpha$ is
norm continuous (cf.~\cite{Ne10}), but if the generator
$D := -i\alpha'(0)$ is unbounded, then the situation is more complicated.
In any case we know from \cite[Thm.~XI.1.23]{Ta03}
that
\[ \Spec_\alpha(\cA) = -i \Spec(\alpha'(0)).\]
\end{rem}

\begin{rem} In \cite[Prop.~15.12]{Str81} it is shown that, for a $W^*$-dynamical system
$(\cM,\R,\alpha)$, the existence of an $\eps > 0$ with
\begin{equation}
  \label{eq:gap2}
\Spec_\alpha(\cM) \cap ([-2\eps,-\eps]\cup [\eps,2\eps]) = \eset
\end{equation}
implies the existence of a hermitian element $A \in Z(\cM^\alpha)$ with $\|A\| \leq \eps/2$, such that
the modified action $\tilde\alpha_t := \Ad(e^{itA}) \alpha_t$ satisfies
\[  \Spec_{\tilde\alpha}(\cM) \cap (-\eps,\eps) = \{0\}.\]
Here the main point is that \eqref{eq:gap2} implies that
$\cN := \cM^\alpha[-\eps,\eps]$ is a subalgebra of $\cM$ on which $\alpha$ is uniformly
continuous, hence of the form $\Ad(e^{-itA})$.
\end{rem}

\begin{ex} \mlabel{ex:bel} Define $U_t \in \cU(L^2(\T))$ by
$(U_t f)(z) = f(e^{it}z)$ and
$\alpha_t(A) := U_t A U_t^*$ for $A \in \cA := \cB(L^2(\T))$.
Then $\Spec_\alpha(\cA) \subeq \Z$ but there is no splitting
as in \eqref{eq:split-con} (\cite[Prop. 1.1]{Be09}).
\end{ex}

We now assume that $\cM \subeq \cB(\cH)$ is a von Neumann algebra
and that $\alpha_t$ is implemented by a Borchers--Arveson
one-parameter group $(U_t)_{t \in \R}$ of $\cM$ with non-negative
spectrum. Let $P$ denote the $\cM$-valued spectral measure of $U$
with $U_t = \int_\R e^{-itx}\, dP(x)$.
For the spectral projections of $U$, we then have
\[ P[t,\infty)\cH = \bigcap_{s < t} \lbr\cM^\alpha[s,\infty)\cH\rbr \]
(cf.~Theorem \ref{BA-thm}).
For $0 < t < \eps$, this leads with \eqref{eq:tridiag} to
\[ P[t,\infty)\cH = \lbr\cM^\alpha[\eps,\infty)\cH\rbr = \lbr\cM_+\cH\rbr. \]
We conclude that, for $0 < t < \eps$,
\[ P[0,t)\cH = (P[t,\infty)\cH)^\bot
= \ker \cM^\alpha(-\infty, -\eps]
= \ker \cM_- = \ker \cM^\alpha(-\infty, 0) \]
consists of ground state vectors (\cite[Prop.~5.3.19(4)]{BR96}). By minimality of $U$,
ground states are contained in $P_0\cH$ for $P_0 := P(\{0\})$ which leads to
\[ P_0 = P[0,\eps).\]
This leads to:

\begin{lem}
Under the assumption above, that  $\cM \subeq \cB(\cH)$ is a von Neumann algebra
and $\alpha_t\in\Aut(\cM)$ is implemented by a Borchers--Arveson
one-parameter group $(U_t)_{t \in \R}$ of $\cM$ with non-negative
spectrum, if
 $0$ is isolated in $\Spec_\alpha(\cM)$, then $0$ is isolated in
$\Spec(U)$.
\end{lem}

\begin{rem}
With $P_+ := \1 - P_0 = P[\eps,\infty)$, we now obtain
\[ \cM = P_0 \cM P_0 + \underbrace{P_+ \cM P_0}_{\subeq \cM_+}
+ \underbrace{P_0 \cM P_+}_{\subeq \cM_-} + P_+ \cM P_+ \]
and as $P_0$ has central support $1$ (since $0$ is isolated in the spectrum of $U$),
it follows from Lemma~\ref{CP1cyclic} that
\[ (\1-P_0) \cH = \lbr P_+\cM \cH_0\rbr.
\]
\end{rem}

\begin{rem} Let $\pi_0$ be the representation of $P_0 \cM P_0$ on $\cH_0$.
Then
\[ \phi(M) := \pi_0(P_0 M P_0) \]
is a completely positive linear map
vanishing on the subspace $\cM_- + \cM_+$ and its restriction to $\cM_0$ is a representation.
Further,
\[ \phi(M^*M) = 0 \quad \mbox{ for } \quad M \in \cM_-,\]
which is equivalent to
\begin{equation}
  \label{eq:kerprod2}
\phi(\cM \cM_-) = \{0\}.
\end{equation}

If, conversely, $(\pi_0, \cH_0)$ is a normal representation of
$P_0 \cM P_0 = \cM_{P_0}$, then $\phi(M) := \pi_0(P_0 M P_0)$ is a completely positive function
on $\cM$ with $\phi(\1) = \1$, so that dilation leads to a representation
$(\pi, \cM)$ containing $(\pi_0, \cH_0)$ as a subrepresentation with respect to
$\cM_{P_0}$.
Clearly, $\pi(U_t)$ defines a unitary one-parameter group with non-negative
spectrum and $\pi(\cM)$-generating space of ground states.
\end{rem}

\begin{prop}
Let  $\cM \subeq \cB(\cH)$ be a factor
and let $\alpha_t\in\Aut(\cM)$ is implemented by a Borchers--Arveson
one-parameter group $(U_t)_{t \in \R}$ of $\cM$ with non-negative
spectrum,  satisfying the spectral gap condition \eqref{eq:gap}.
Then there exists at most countably many projections
$(P_j)_{j \in J}$ and pairwise different $\lambda_j \geq 0$ with
$U_t = \sum_{j \in J} e^{-it\lambda_j} P_j$. For $j \not=k$, we further have
$|\lambda_j - \lambda_k| \geq \eps$.
\end{prop}

\begin{prf} Let $0 \leq a < b$ such that $2(b-a)\leq \eps$ and
$M \in P[a,b]\cM P[a,b]$. Then
$\Spec_\alpha(M) \subeq [a,b] - [a,b]  \subeq [-\eps,\eps]$ implies that
\[ P[a,b]\cM P[a,b] \subeq \cM_0. \]
For disjoint compact subsets
$S_1, S_2 \subeq [a,b]$, this further leads to
\[ P(S_1)\cM P(S_2) \subeq \cM^\alpha(S_1 - S_2) = \{0\}\]
because $S_1 - S_2 \subeq [-\eps,\eps]$ does not contain $0$.
In view of \cite[Prop.~1.10.7]{Sa71}, the central supports of
$P(S_1)$ and $P(S_2)$ are disjoint. As $\cM$ is a factor, we obtain
$P(S_1) = 0$ or $P(S_2) = 0$. This implies that $U$ has at most a single spectral
value in the interval $[a,b]$, and from that we  derive that
$\Spec(U)$ is discrete,
so that
$U_t = \sum_{j \in J} e^{-it\lambda_j} P_j$ as asserted.
Then the differences $\lambda_j - \lambda_k$ are contained in $\Spec_\alpha(\cM)$,
which implies that $|\lambda_j -\lambda_k|\geq\eps$ for $j \not=k$.
\end{prf}

\begin{ex} In general, if $\cM$ is not a factor, the assumption that $0$ is isolated in $\Spec_\alpha(\cM)$
does not imply that $\Spec(U)$ is discrete.
In $\cM := \ell^\infty(\N,\cB(\ell^2))$, we consider the \minimal  one-parameter group
given by $U_t=(U_t^{(1)},U_t^{(2)},\ldots)$ with $U_t^{(n)}\in\cB(\ell^2)$ defined by
\[ U_t^{(n)} := P_{1} + \sum_{j = 0}^\infty e^{it(j+1 + f(n))} P_{j+2},\]
where $P_j, j \in \N,$ is the orthogonal projection onto $\C e_j$ and
$f \: \N \to \Q_+$ is surjective. Then
\[ \Spec(U) = \{0\} \cup [1,\infty). \]
and, for $\alpha = \Ad(U)$ the block diagonal structure leads to
$\Spec_\alpha(M) = \oline{\bigcup_n \Spec_{\alpha^{(n)}(\cM)}},$
which in turn leads to
\[ \Spec_\alpha(\cM) = (-\infty, -1] \cup \{0\} \cup [1,\infty).\]
\end{ex}

A specific instance where  $0$ is isolated in $\Spec_\alpha(\cM)$ is the periodic case.
We continue  analysis of  the periodic case, started above
in Example~\ref{Zmatrix}.
Let $(\cM,\T, \alpha)$ be a $W^*$-dynamical system and
$(U_t)_{t \in \R}$ be a weakly continuous unitary one-parameter group in $\cM$
with positive spectrum such that $\alpha_{e^{it}} = \Ad(U_t)$ for $t \in \R$ and
$U$ is  minimal (cf.~Definition~\ref{def:minimal}).

The $2\pi$-periodicity of $U$ implies the existence
of projections $(P_n)_{n \in \N_0}$ in $\cM$ with
\[ U_t = \sum_{n = 0}^\infty e^{-int} P_n \quad \mbox{ for } \quad t \in \R.\]
In this case $P_0 = P[0,\eps)$ for $0 \leq \eps \leq 1$, so that
Lemma~\ref{lem:minimal} implies that the central support of $P_0$ is $\1$. With
Lemma~\ref{CP1cyclic} this leads to
\[ \cM = \oline{\cM P_0 \cM}^w.\]
Put $\chi_n(t) := e^{-int}$ and
\[ \cM_n := \{ M \in \cM\; \mid\; (\forall t \in \R)\ \alpha_{e^{it}}(M) = e^{-int} M\}.\]
Then the subspaces $P_n\cM P_m$ are $\alpha$-eigenspaces with respect to the character
$\chi_{n-m}$ and the direct vector space sum
 $\sum_{k = -n}^\infty P_{k+n}\cM P_k$ is weakly dense in $\cM_n$ for $n \in \Z$.
In particular, the fixed point algebra
$\cM_0$ is the weak closure of $\sum_{k = 0}^\infty P_{k}\cM P_k$, where
the subalgebras  $P_k \cM P_k$ of $\cM$ are  two-sided ideals of $\cM_0$  (as $[\cM_0,P_n]=0$).

Above in Example~\ref{Zmatrix} we noted that the fixed point projection
$p_0 \: \cM \to \cM_0$ by
\[
p_0(M) := \sum_{k=0}^\infty P_kMP_k= \int_\T \alpha_z(M)\, dz\qquad
\hbox{for}\;M\in\cM
\]
is completely positive.  Hence we can use the Stinespring dilation to build
representations on $\cM$ from representations on $\cM_0$.

Consider the case of a $\T$-action, i.e.\  a $W^*$-dynamical system.
For $A_n \in \cM_n$ we have
\[ \alpha_f(A_n)
= \int_\T f(t)\alpha_t(A_n)\, dt
= \int_\T f(t)e^{- int}\, dt \cdot A_n
= \hat f(n) A_n.\]
Therefore
\[ \pi(\alpha_f(A_n)) \Omega = \hat f(n) \pi(A_n) \Omega \]
vanishes if $\Omega$ is a ground state vector and $\supp(\hat f) \subeq -\N$.

\section{KMS states and modular groups.}
\mlabel{KMSsect}

A major area where covariant representations
of singular actions are studied is
that of KMS states and their representations.
This is  a fundamental
part of the study of thermal quantum systems, and the literature
in this area is vast. This section is only a scratch on the surface, and we will concentrate
on some of  the main structural issues. The standard references
include \cite{BR96}, \cite{SZ79} and for the
case of $W^*$-actions, a useful review of results is in \cite{DJP03}.
For a particularly interesting application in QFT, see \cite{CR94}.

\subsection{Modular group of a weight on a von Neumann algebra.}
\label{ModGps}

First, we need to define the modular group (proofs and
constructions are in  \cite[Ch.~10]{SZ79} and
\cite[Sect.~III.4]{Ta03}).
 Let $\cM\subset\cB(\cH)$ be a von
Neumann algebra, and let $\varphi$ be a
faithful, normal semifinite weight on $\cM$. Recalling the GNS construction for
it, consider the left ideal
\[
\cN_\varphi:=\{A\in\cM\,\mid\,\varphi(A^*A)<\infty\}.
\]
By faithfulness of $\varphi$ the sesquilinear form
$\la A,B\ra:={\varphi(A^*B)},$ $A,B\in \cN_\varphi$ is positive definite,
hence we may complete $\cN_\varphi$ to obtain the Hilbert space $\cH_\varphi$.
Let $\xi:\cN_\varphi\to\cH_\varphi$ denote the faithful linear imbedding.
Then the GNS representation $\pi_\varphi:\cM\to \cB(\cH_\varphi)$ is given by
\[
\pi_\varphi(A)\xi(B):=\xi(AB)\quad\hbox{for}\quad A\in \cM,\;B\in\cN_\varphi
\]
and it is faithful. There may be no cyclic vector in $\cH_\varphi$,
unless $\varphi$ is bounded.
By Theorem~\ref{Def-SF}, this GNS representation is unitarily equivalent
to the standard form realization  of $\cM$.
On the subspace $\cD_\varphi:=\xi(\cN_\varphi\cap\cN_\varphi^*)\subset \cH_\varphi$ there
is a closable conjugate linear operator $S_0$ defined by
\[
S_0\xi(A):=\xi(A^*)\quad\hbox{for}\quad A\in \cN_\varphi\cap\cN_\varphi^*.
\]
Denote its closure by $S_\varphi$. Then the {\it modular operator of}
$\varphi$ is the invertible positive  operator (in general, unbounded)
given by
\[
\Delta_\varphi:=S_\varphi^*S_\varphi.
\]
The {\it modular conjugation of} $\varphi$ is the operator $J_\varphi:=\Delta^{1/2}_\varphi
S_\varphi$.
Then $(\Delta^{it}_\varphi)_{t \in \R}$
defines a strong operator continuous one parameter unitary
group, and as $\Delta^{it}_\varphi\cM\Delta^{-it}_\varphi = \cM$, this defines the
{\it modular automorphism group} $(\sigma^\varphi_t)_{t \in \R}$ in $\Aut\cM$ by
\[
\sigma^\varphi_t(A):=\Delta^{it}_\varphi A\Delta^{-it}_\varphi \quad \mbox{ for }
\quad A\in\cM,\; t\in\R,\]
which is obviously a $W^*$-dynamical system. It is covariant by construction, and the generator
of its implementers is $L_\varphi:=-\ln\Delta$ (called the {\it standard Liouvillean}),
 i.e.
\[\Delta^{it}_\varphi=\exp(-itL_\varphi).\]
The relation between different modular groups
on the same von Neumann algebra is given by:
\begin{thm} {\rm(Connes)} Let $\cM$ be a von Neumann algebra, and let $\varphi$ be a
faithful, normal semifinite weight on $\cM$.
\begin{itemize}
\item[\rm(i)] If $\psi$ is another faithful, normal semifinite weight on $\cM$,
then there is a unique \strongH continuous path of unitaries $(u_t)_{t\in\R}\subset\cM$
such that
\[
\sigma^\psi_t(A)=u_t\sigma^\varphi_t(A)u_t^*\qquad\hbox{and}\qquad
u_{t+s}=u_t\sigma_t^\varphi(u_s).
\]
\nin{\rm We write $(D\psi :D\varphi)_t := u_t$.}
\item[\rm(ii)] Conversely,  if  a \strongH continuous path of unitaries $(u_t)\subset\cM$
satisfies $u_{t+s}=u_t\sigma_t^\varphi(u_s)$ for all $t,s\in\R$, then there is a unique
faithful, normal semifinite weight $\psi$ on $\cM$ with
 $(u_t)=(D\psi :D\varphi)_t$ for all $t$.
\end{itemize}
\end{thm}

This is proved in \cite[Thm.~III.4.7.5]{Bla06}, and  in
\cite[Thms.~VIII.3.3, VIII.3.8]{Ta03}.
The modular group also affects the adjoint action of unitary one-parameter groups with positive spectrum on $\cM$
(cf. \cite{F98} for a direct proof and \cite[Th. 2.1]{ArZs05} for a general version):

\begin{thm}{\rm(Borchers' Theorem on modular inclusions;\cite{Bo92})}
\label{BorTMI}
Let $\cM\subset\cB(\cH)$ be a von Neumann algebra, and let $\Omega\in\cH$ be a cyclic
and separating vector with associated  vector state
$\omega(M):=\la\Omega,M\Omega\ra$.
Let $(U_s)_{s \in \R}$ be a   unitary one-parameter group with positive spectrum on $\cH$ such that
\[U_s\Omega=\Omega\quad \mbox{ for } \quad s \in \R
\qquad\mbox{ and } \qquad
U_s\cM U_s^*\subseteq\cM\quad \mbox{ for } \quad s \geq 0. \]
Then
\begin{itemize}
\item[\rm(i)] $\sigma^\omega_t(U_s) = \Delta^{it}_\omega U_s\Delta^{-it}_\omega
=U_{e^{-2\pi t}s}$ for $s,t \in \R$, and
\item[\rm(ii)] $J_\omega U_sJ_\omega=U_s^*$ for $s\in\R$.
\end{itemize}
\end{thm}

It is quite remarkable that there exist homomorphisms
$\alpha \: \R \to \Inn(\cM)\equiv$ the inner automorphisms of $\cM$,
which define $W^*$-dynamical systems which do not
lift to $\cU(\cM)$, i.e.\  the corresponding central extension
$\hat \R := \alpha^*\cU(\cM)$ of $\R$ by $\cU(Z(\cM))$ is non-trivial
(\cite[\S 15.16]{Str81}). Here is the main result behind these examples:

\begin{thm}\label{D2}
A $W^*$-algebra $\cM$ is semifinite if and only if
the modular automorphism group
of one of its faithful normal semifinite weights is implemented by a
unitary one-parameter group in $\cU(\cM)$.
Then the modular automorphism groups
of all faithful normal semifinite weights are implemented by a unitary
one-parameter  group in $\cU(\cM)$.
\end{thm}

\begin{prf}
See \cite[Th. 7.4]{PT73}, which goes back to \cite[Ch. 14]{Ta70}.
\end{prf}

It follows by Theorem~\ref{D2} that for any faithful normal semifinite weight of a
factor $\cM$ of type~III, its corresponding
modular automorphism group cannot be implemented by a unitary one-parameter
group of $\cM$.
Consequently, the factor $\cM$ of type III given by \cite[Cor.~1.5.8(c)]{Co73}
has the remarkable property that, for every faithful normal semifinite weight,
its modular automorphism group consists of inner automorphisms
and yet it is not implemented by any one-parameter unitary group in~$\cM$.
As explained in \cite[p.~21]{Ta83}, this property can be shared only by
(possibly countably decomposable) $W^*$-algebras with nonseparable predual.

On the positive side, there are nice results of Kallman and Moore
building on measurable sections and Polish group structures.
This requires $\cM_*$ to be separable.
More concretely, in \cite{Ka71} one finds that, for $G = \R$ and
$\cM_*$ separable, all inner $W^*$-dynamical systems
can be implemented by one-parameter groups $U \: \R \to \cU(\cM)$.
Note that the separability of $\cM_*$ implies that
the standard representation of $\cM$ is separable because
the cone $\cC \cong \cM_{*,+}$ is separable.

\subsection{KMS condition for a weight with respect to a  $C^*$-action $(\cA, \R, \alpha)$.}

A weight $\varphi$  and its modular group $\sigma^\varphi$ satisfy the
modular condition:
\begin{defn}\mlabel{defKMS}
Given a  $C^*$-action $(\cA, \R, \alpha)$, possibly singular, then
 a lower semicontinuous weight $\varphi$
on $\cA$ is said to satisfy the {\it KMS condition for $\alpha$ at
$\beta\not=0$} if
\begin{itemize}
\item[\rm(i)] $\varphi=\varphi\circ\alpha_t$ for all $t\in\R$,
\item[\rm(ii)] for every pair $A,B\in\cN_\varphi\cap\cN_\varphi^*$, there exists
a bounded
continuous function $F$ on the closed horizontal strip $S_\beta\subset\C$
where
\begin{eqnarray*}
S_\beta&:=&\{z\in\C\,\mid\, 0\leq \pm{\rm Im}(z)\leq\pm\beta\}\qquad\hbox{if}\quad \pm\beta>0
\qquad\hbox{(matched signs).}
\end{eqnarray*}
Moreover, $F$ is analytic on the interior of $S_\beta$ and satisfies for all $t\in\R$:
\begin{equation}\label{KMScond}
F(t)=\varphi(\alpha_t(A)B),\qquad
F(t+i\beta)=\varphi(B\alpha_t(A)).
\end{equation}
\end{itemize}
For the case $\beta=1$ we call the KMS condition the {\it modular condition}.
By rescaling  $\alpha$, we see that
$\varphi$ satisfies the KMS condition for $\alpha$ at
$\beta\not=0$ if and only if it satisfies the modular condition for $\alpha_{\beta t}$.
If $\varphi$ is a state, it will be called a {\it KMS state for $\alpha$ at $\beta$} or just a KMS state for short.
\end{defn}
\begin{rem}
(a)
In physical models with KMS states, $\beta$ is identified with the (negative) inverse temperature.
In the case that $\varphi$ is a state
(which is the case if $\cA$ is unital and $\1\in\cN_\varphi$),
the invariance condition (i) is redundant,
as invariance then follows from (ii). To see this,
note that condition~(\ref{KMScond})
implies that for every $A\in\cN_\varphi\cap\cN_\varphi^*$, there exists
a bounded continuous function $F$ on the closed horizontal strip $S_\beta\subset\C$ which is analytic on the interior, and such that
\[
F(t)=\varphi(\alpha_t(A))=
F(t+i\beta)\,.
\]
This is obtained by either substituting $\1$ for $B$ into \eqref{KMScond}, or by substituting an approximate
identity for $B$ into \eqref{KMScond}, and taking the limit (which is uniform in $t$).
This means that we can define a new function $\tilde{F}$ on the entire complex plane by tiling $\C$ with
vertical translates of the strip
$S_\beta$, carrying along the values of $F$ on $S_\beta$.
Then $\tilde{F}$ is continuous, bounded and analytic everywhere except on the horizontal lines where
the strips join. By Morera's Theorem, it is  in fact analytic also on these joining lines,
i.e.\ it is entire, and as it is bounded, by Liouville's
theorem it is constant. Thus $F(t)=\varphi(\alpha_t(A))$ is constant,
i.e.\ (i) holds (see \cite[Prop.~5.3.3]{BR96} for more details).

(b)
For $C^*$-dynamical systems,  Pusz and Woronowicz  showed that both
ground states and KMS states are  ``passive'' states (\cite[Thm.~1.2]{PW78}), i.e.
\[
\omega\big(-iU^*\delta(U)\big)\geq 0\quad\hbox{for all}\quad
U\in \cU_0(\cA)\cap \cD(\delta)
\]
where $\delta$ is the generator of $\alpha$ with domain $\cD(\delta)\subseteq\cA$
and $\cU_0(\cA)$ denotes the identity component of the group of unitaries in $\cA$
w.r.t. the norm topology
(\cite[Thm.~2.1]{PW78}).
Conversely, if a passive state is weakly clustering, then
it is either KMS or a ground state \cite[Thm.~1.3]{PW78}.
\end{rem}

The modular condition in fact uniquely characterizes the modular group of a weight by:
\begin{thm}\label{ModProp}
Let $\cM$ be a von Neumann algebra, and let $\varphi$ be a
faithful normal semifinite weight on $\cM$. Then
the modular automorphism group $(\sigma^\varphi_t)_{t \in \R}$ in $\Aut\cM$ satisfies
the modular condition  for $\varphi$. Conversely,
for any  $W^*$-dynamical system $(\cM, \R, \alpha)$
which satisfies the modular condition for $\varphi$,
the modular group $\sigma^\varphi$ coincides with $\alpha$.
\end{thm}

This is proven in \cite[Thm.~VIII.1.2]{Ta03} and
\cite[Thm. p.~289]{SZ79}. Thus, every faithful, normal semifinite weight on $\cM$
is a KMS weight for a unique one-parameter automorphism group.

\begin{thm}\label{CKMS}
Given a  $C^*$-action $(\cA, \R, \alpha)$,  let $\omega$ be a
faithful state on $\cA$ which satisfies the modular condition
for $\alpha$.
Then the normal extension $\tilde\omega$ of $\omega$ to
$\cM:=\pi_\omega(\cA)''$ is faithful, and satisfies
\[
\pi_\omega\circ\alpha_t=\sigma_t^{\tilde\omega}\circ\pi_\omega\quad\mbox{ for } \quad t\in\R.
\]
\end{thm}

This is proven in \cite[Prop.~VIII.1.5]{Ta03}.
In fact, the requirement that $\omega$ is faithful is too strong, one only
needs that  the vector state ${(\Omega_\omega,\cdot\,\Omega_\omega)}$ is faithful on
$\cM:=\pi_\omega(\cA)''$ (cf.  \cite[Thm.~5.3.10]{BR96}), and this can happen even when
$\pi_\omega$ is not faithful.
A state on $\cA$ which satisfies the KMS condition
for $\alpha$ can therefore be characterized by this condition, i.e.
that its GNS representation $\pi_\omega$ intertwines $\alpha$ with a
rescaled copy of its modular group.

\begin{prop}\label{KMSsingtoW}
Given a  $C^*$-action $(\cA, \R, \alpha)$, possibly singular, let $\omega$ be a
KMS state for $\alpha$ at $\beta$. Then the following hold:
\begin{itemize}
\item[\rm(i)] $(\pi_\omega,U^\omega)$ is covariant.
\item[\rm(ii)] the normal extension $\tilde\omega$ of $\omega$ to $\cM:=\pi_\omega(\cA)''$
by $\tilde\omega(M):={\la\Omega_\omega,M\Omega_\omega\ra}$
is faithful.
\item[\rm(iii)] If $\tilde\alpha_t:={\rm Ad}U^\omega_t$, then
$(\cM, \R, \tilde\alpha)$ is a $W^*$-dynamical system
 for which $\tilde\omega$ is a KMS state for $\tilde\alpha$ at~$\beta$.
 \item[\rm(iv)] $\cM'\cap\cM\subseteq\cM^{\tilde\alpha}$, the set of invariant
elements of $\cM$ with respect to ~$\tilde\alpha$ (modular automorphisms
act trivially on the center).
\item[\rm(v)] Let $\cN\subseteq\cM$ be a commutative von Neumann subalgebra
such that $\tilde\alpha_t(\cN)\subseteq\cN$ for all $t\in\R$.
Then $\cN\subseteq\cM^{\tilde\alpha}$.
\item[\rm(vi)] $\cM^{\tilde\alpha}=\{A\in\cM\,\mid (\forall
B \in\cM) \;\omega([A,B])=0\}$.
\end{itemize}
\end{prop}

\begin{prf} (i) By the KMS condition, for every $A, B\in\cA$, the function
$t\mapsto\omega(\alpha_t(A)B)$ is continuous. By invariance of $\omega$ this implies for
the GNS unitaries  that $(U^\omega_t)_{t \in \R}$ is strong operator continuous, and so $(\pi_\omega,U^\omega)$ is covariant.

(iii) By assumption $\tilde\omega$ satisfies the KMS condition with respect to ~$\tilde\alpha$
on the strong operator dense subalgebra $\pi_\omega(\cA)\subset\cM$.
By substituting an approximate identity for $B$ in the KMS condition (\ref{KMScond}),
taking the limit and using Liouville's theorem, we conclude that $\omega$
is  $\alpha$-invariant on $\pi_\omega(\cA)$, hence on all of $\cM$.
By  Lemma~\ref{KMSdom} below,
it then follows that  $\tilde\omega$ satisfies the KMS condition with respect to ~$\tilde\alpha$
on all of $\cM$.

(ii) By \cite[Theorem~5.3.10]{BR96} it follows from (iii) that $\tilde\omega$
is faithful on $\cM$.

(iv) Let $C\in\cM'\cap\cM$ and $A,\,B\in \cM$. Then, by (iii), we have for some
continuous bounded function $F$ on the strip $S_\beta$ that it is
holomorphic on the interior, and on the boundary
\[
F(t)=\tilde\omega(\alpha_t(AB)C)=\tilde\omega(C\alpha_t(AB))=
F(t+i\beta).
\]
Proceeding as above,  we define a new function $\tilde{F}$ on the entire complex plane by tiling $\C$ with
vertical translates of the strip
$S_\beta$, carrying along the values of $F$ on $S_\beta$.
Then $\tilde{F}$ is continuous, bounded and analytic everywhere except on the horizontal lines where
the strips join. By Morera's Theorem, it is  in fact analytic also on these joining lines,
i.e.\ it is entire, and as it is bounded, by Liouville's
theorem it is constant. Thus,  $F$ is constant, and equal to
\[
F(t)=\tilde\omega(\alpha_t(A)C\alpha_t(B))=\tilde\omega(A\alpha_{-t}(C)B)
=\la A^*\Omega_\omega, \alpha_{-t}(C)B\Omega_\omega\ra\quad\mbox{ for } \quad t\in\R.\]
As $\Omega_\omega$ is cyclic and $F$ is constant, we get that $C\in\cM^{\tilde\alpha}$.

(v) As the restriction $\omega_0$ of $\tilde\omega$ to $\cN$ is still a KMS-state
with respect to the restriction $\alpha^{(0)}$ of $\tilde\alpha$ to $\cN$, it follows that $\alpha^{(0)}$
coincides with the modular automorphism with respect to  $\omega_0$.
Thus by (iv), as $\cN$ is commutative, we have that
$\pi_{\omega_0}(\cN)\subseteq\pi_{\omega_0}(\cN)^{\alpha^{(0)}}$, i.e.
$\pi_{\omega_0}(\tilde\alpha_t(N)-N)=0$ for all $N\in \cN$ and $t\in\R$.
By (ii), $\tilde\omega$ is faithful, hence its restriction $\omega_0$ is faithful,
and so $\pi_{\omega_0}$ is faithful. Thus $\tilde\alpha_t(N)=N$ for all $N,\;t$,
i.e. $\cN\subseteq\cM^{\tilde\alpha}$.

(vi) is proven in  \cite[Prop.~5.3.28]{BR96}.
\end{prf}

\begin{lem}
\label{KMSdom}
Let $(\cM, \R, \alpha)$ be a $W^*$-dynamical system, and let $\omega$ be
a normal $\alpha$-invariant state
satisfying the KMS condition \eqref{KMScond} for all $A,\,B$ in some
$\sigma(\cM,\cM_*)$-dense $\alpha$-invariant unital *-subalgebra $\cD$ of $\cM$. Then $\omega$
satisfies \eqref{KMScond} on all of $\cM$, hence is a
KMS state for $\alpha$ at $\beta$.
\end{lem}
\begin{prf} (Adapted from that of \cite[Prop.~5.3.7]{BR96})
 Let $A,\,B\in \cM$ be arbitrary, and let $(A_\nu)_{\nu\in\Gamma}$
and $(B_\nu)_{\nu\in\Gamma'}$ be nets in $\cD$ which $W^*$-converge to $A$ and $B$
respectively. We can choose the same directed set $\Gamma=\Gamma'$ for both nets,
and by Kaplansky's density theorem (cf. \cite[Thm.~2.4.16]{BR02}) we may choose
$\|A_\nu\|\leq\|A\|$, $\|B_\nu\|\leq\|B\|$ for all $\nu\in\Gamma$. By assumption,
for each pair $A_\nu,\,B_\nu\in\cD$, there exists a bounded
continuous function $F_\nu$ on the closed horizontal strip $S_\beta\subset\C$
which is holomorphic on the interior of $S_\beta$, and satisfies
\[
F_\nu(t)=\omega(\alpha_t(A_\nu)B_\nu),\qquad
F_\nu(t+i\beta)=\omega(B_\nu\alpha_t(A_\nu))\quad \mbox{ for } \quad t\in\R.
\]
Let $\nu>\mu\in\Gamma$.
Then, by the so-called three-line theorem \cite[Prop.~5.3.5]{BR96}, the positive function
$z\mapsto{\big|F_\nu(z)-F_\mu(z)\big|}$ takes its maximum on the boundary
of $S_\beta$, and hence for any $z\in S_\beta$ we obtain  on $\cM$:
\[ \big|F_\nu(z)-F_\mu(z)\big|\leq \max\Big\{\sup_{t\in\R}
\big|\omega(\alpha_t(A_\nu)B_\nu-\alpha_t(A_\mu)B_\mu)\big|,\;
\sup_{t\in\R}\big|\omega(B_\nu\alpha_t(A_\nu)-B_\mu\alpha_t(A_\mu))\big|\Big\}.\]
Now using
the $\alpha$-invariance of $\omega$ we have
\begin{align*}
&\big|\omega(\alpha_t(A_\nu)B_\nu -\alpha_t(A_\mu)B_\mu)\big|\\
=& \big|\omega\big(\alpha_t(A_\nu-A)B_\nu-\alpha_t(A_\mu-A)B_\mu
+  \alpha_t(A)(B_\nu-B)-\alpha_t(A)(B_\mu-B)\big)\big| \\
\leq &\|B\|\big(\|\pi_\omega(A_\nu^*-A^*)\Omega_\omega\|
+\|\pi_\omega(A_\mu^*-A^*)\Omega_\omega\|\big)
+\|A\|\big(\|\pi_\omega(B_\nu-B)\Omega_\omega\|
+\|\pi_\omega(B_\mu-B)\Omega_\omega\|\big).
\end{align*}
This expression converges uniformly with respect to ~$t$
 to zero as $\nu$ and $\mu\nearrow\infty$.
Likewise the other term ${\big|\omega(B_\nu\alpha_t(A_\nu)-B_\mu\alpha_t(A_\mu))\big|}$
converges uniformly with respect to ~$t$
 to zero as $\nu$ and $\mu\nearrow\infty$, hence ${\big|F_\nu(z)-F_\mu(z)\big|}$
 converges uniformly with respect to ~$z$
 to zero as both $\nu$ and $\mu\nearrow\infty$, hence $(F_\nu)_{\nu\in\Gamma}$ is a Cauchy net
 which converges uniformly, hence the limit function $F(z)$ is continuous and bounded on
 $S_\beta\subset\C$ and analytic on its interior. As
\begin{eqnarray*}
F(t)&=& \lim_\nu\omega(\alpha_t(A_\nu)B_\nu)=\omega(\alpha_t(A)B)\quad\hbox{and}\\[1mm]
F(t+i\beta)&=&\lim_\nu F_\nu(t+i\beta) = \lim_\nu\omega(B_\nu\alpha_t(A_\nu))
=\omega(B\alpha_t(A)),
\end{eqnarray*}
it follows that  $\omega$ satisfies (\ref{KMScond}) for all $A,\,B\in \cM$.
\end{prf}

\begin{rem} Recall the context of Proposition~\ref{KMSsingtoW}.
\begin{itemize}
\item[\rm(a)] By Proposition~\ref{KMSsingtoW}(ii), $\cM=\pi_\omega(\cA)''$ is in standard form.
\item[\rm(b)]
By Theorem~\ref{CKMS}, $\tilde\alpha$ coincides with a
rescaled copy of the modular group of $\tilde\omega$, hence there are
strong restrictions on the existence of a KMS state for a given
 $C^*$-action.
 \item[\rm(c)] By Proposition~\ref{KMSsingtoW}(iv),
the modular automorphism group of a KMS state acts trivially
on the center of the corresponding von Neumann algebra. This
means that it adapts to the central disintegration of this
algebra into factors. Therefore the main point in understanding modular
automorphism groups concerns factors.
 \item[\rm(d)] By the fact that the modular automorphism group of a KMS state acts trivially
on the center of the corresponding von Neumann algebra, it is easy to give an example
of a $W^*$-dynamical system which has no normal faithful KMS states. Just take any one with
an automorphism group which is not trivial on the center.
\item[\rm(e)] By Proposition~\ref{KMSsingtoW}(iv),
if $\cM$ is commutative, it can only have KMS states for the
 trivial action. Compare this with the analogous property for ground states
 (cf.~Remark~\ref{RemBorAr}(b)).
 Moreover, by Proposition~\ref{KMSsingtoW}(v), the group $\tilde\alpha$ cannot have normal
 eigenvectors, unless they are invariant.
 \item[\rm(f)]
 The spectrum of the implementing group $U^\omega$ has been examined, and under some conditions
one can even prove that ${\rm Sp}(U^\omega)$ is independent of $\omega$ and
$\beta$ (cf.~\cite[Thm.~A]{tBW76}).
However, as $\cM$ is in standard form and  $U^\omega_t$ are the
standard form implementers given by Proposition~\ref{prop:2.9.1}
 (using \cite[Prop.~IX.1.17]{Ta03}),
 and the spectrum of $U^\omega$ equals the Arveson spectrum of $\tilde\alpha$
 by Proposition~\ref{SpectraNormalF},  the reason for this is clear.
 \end{itemize}
\end{rem}
 We list a few  equivalent conditions, where the extension
of the $\R$-action to $\pi_\omega(\cA)''$ is assumed;
criteria for this are given in Corollary~\ref{cor:bor1}(iii).

\begin{thm}\label{KMSequiv}
Given a  $C^*$-action $(\cA, \R, \alpha)$, possibly singular, let $\omega$ be a
state on $\cA$ such that the induced action of $\R$ on $\pi_\omega(\cA)$ extends to an
action  $\wt\alpha:\R\to\Aut(\pi_\omega(\cA)'')$, and defines a $W^*$-dynamical system.
Denote the normal extension of $\omega$ to $\cM:=\pi_\omega(\cA)''$ by $\tilde\omega$.\\[1mm]
Then the following are equivalent for $\beta>0$:
\begin{itemize}
\item[\rm(i)] $\omega$ is a
KMS state on $\cA$ for $\alpha$ at $\beta$
\item[\rm(ii)]
$
\tilde\omega\big(A\tilde\alpha_{i\beta}(B)\big)=\tilde\omega\big(BA\big)\;\;
$
for all $A,\,B$ in some $W^*$-dense $\tilde\alpha$-invariant *-subalgebra
of the entire elements  in $\cM$ of $\tilde\alpha$.
 \item[\rm(iii)]
  $\tilde\omega$ is $\tilde\alpha$-invariant, and satisfies the spectral condition:
  \begin{equation}
    \label{eq:6.10}
  \tilde\omega(A^*A)\leq e^{\beta\lambda}\tilde\omega(AA^*)\qquad\hbox{for all}\quad
   A\in\cM^{\tilde\alpha}(-\infty,\lambda)
  \quad\hbox{and}\quad  \lambda\in\R,
  \end{equation}
 where $\cM^{\tilde\alpha}(-\infty,\lambda)$ denotes the Arveson spectral subspaces.
 \item[\rm(iv)] For all $A,B\in\cA$ and $f$ with $\hat{f} \in C^\infty_c(\R)$, we have:
\[
 \int_{\R}f(t)\, \omega\big(A\alpha_t(B)\big)\,dt=\int_{\R}f(t+i\beta)\,\omega\big(\alpha_t(B)A\big)\, dt\,.
 \]
\end{itemize}
\end{thm}

\begin{prf}
(i)$\Leftrightarrow$(ii): (i) gives via Proposition~\ref{KMSsingtoW} the
$W^*$-dynamical system
$\wt\alpha:\R\to\Aut(\pi_\omega(\cA)'')$, satisfying the KMS condition for $\tilde\omega$. By
 \cite[Prop.~5.3.7, Def.~5.3.1]{BR96},
 this is equivalent to the condition given in (ii).\\[1mm]
 (ii)$\Leftrightarrow$(iii): The restriction of the  $W^*$-dynamical system
$\wt\alpha:\R\to\Aut(\cM)$ to its $W^*$-dense continuous subalgebra $\cM_c$
is a $C^*$-dynamical system. Assume (ii). Then
the restriction of $\tilde\omega$ to  $\cM_c$
is still KMS, hence by \cite[Thm.~1.1]{dC82}, this is
equivalent for the    $\tilde\alpha$-invariant $\tilde\omega$ to satisfy
\begin{equation}
\label{Canniere}
  \tilde\omega(A^*A)\leq e^{\beta\lambda}\tilde\omega(AA^*)\qquad\hbox{for all}\quad
   A\in(\cM_c)^\alpha(-\infty,\lambda)
  \quad\hbox{and}\quad  \lambda\in\R\,.
 \end{equation}
  Recall from \cite[Lemma~3.2.39(4)]{BR02} that
  \[
 \cM^{\tilde\alpha}(-\infty,\lambda) = \cM_0^{\tilde\alpha}(-\infty,\lambda)
 =\overline{{\rm Span}\big\{\tilde\alpha_f(M)\,\big|\,M\in\cM,\;
 f\in L^1(\R),\;{\rm supp}\hat{f}\subset(-\infty,\lambda)\big\}}^\sigma,   \]
  where the closure is $W^*$-closure, but for the  $C^*$-dynamical system on  $\cM_c$,
  the corresponding expression has a norm closure. However, as the maps
  $M\mapsto \tilde\alpha_f(M)$ are {$\sigma(\cM,\cM_*)\mathord{-}\sigma(\cM,\cM_*)\hbox{-continuous}$}
  (cf.~\cite[Prop.~3.1.4]{BR02})
  it follows that
  $(\cM_c)^{\tilde\alpha}(-\infty,\lambda)$ is $W^*$-dense in  $\cM^{\tilde\alpha}(-\infty,\lambda)$.
  As $\tilde\omega$ is normal, by substituting for $A$ in condition (\ref{Canniere}) a net
  in  $(\cM_c)^{\tilde\alpha}(-\infty,\lambda)$ which  $\sigma(\cM,\cM_*)\hbox{-converges}$
  to some $M\in  \cM^{\tilde\alpha}(-\infty,\lambda)$, we obtain
\eqref{eq:6.10} for $A = M$.

  For the converse, assume (iii).
Then the condition restricts to the  $C^*$-dynamical system
on the $W^*$-dense continuous subalgebra $\cM_c$, using
$(\cM_c)^{\tilde\alpha}(-\infty,\lambda)\subseteq\cM^{\tilde\alpha}(-\infty,\lambda)$.
Thus, by \cite[Thm.~1.1]{dC82}, $\tilde\omega$ is KMS on $\cM_c$, hence it satisfies (ii) for the norm-dense
subalgebra of $\cM_c$ consisting of the entire elements of $\tilde\alpha$ (cf.  \cite[Prop.~5.3.7, Def.~5.3.1]{BR96}).
As this subalgebra is $W^*$-dense in $\cM$, (ii) is satisfied.\\[2mm]
  (ii)$\Leftrightarrow$(iv):
  First write condition (iv) as
  \[
  \omega\big(A\alpha_f(B)\big)=\omega\big(\alpha_{f_\beta}(B)A\big)\quad\hbox{for}\quad
  f_\beta(t):=f(t+i\beta)
  \]
  where  $A,B\in\cA$ and $f$ with $\hat{f}\in\cD$. As the maps
  $M\mapsto \tilde\alpha_f(M)$ are {$\sigma(\cM,\cM_*)\mathord{-}\sigma(\cM,\cM_*)\hbox{-continuous}$}
  (cf.~\cite[Prop.~3.1.4]{BR02}), we can extend this condition to all $\cM$. Then
  the equivalence of (iv) with (ii) on the $C^*$-dynamical subsystem of $\tilde\alpha$
  restricted to $\cM_c$ is given in \cite[Prop.~5.3.12]{BR96}.
  As the dense $\tilde\alpha$-invariant *-subalgebra
of the entire elements  in $\cM_c$ of $\tilde\alpha$ are $W^*$-dense in $\cM$, the equivalence
with (ii) follows.
\end{prf}
Only condition (iii) needs $\beta>0$.
Note that condition (ii), whilst commonly used for $C^*$-dynamical systems, is not that useful for
singular actions, as the subalgebra of  analytic elements on which it holds, may  have zero
intersection with $\pi_\omega(\cA)$ by Example~\ref{Ac0}.
There is a range of other equivalent conditions for the KMS condition, e.g.
in terms of correlation functions (cf.~\cite{FvB77}), Green's functions
(cf.~\cite{GJO94}),
spectral passivity (cf.~\cite{dC82}), and in terms of stability with respect to ~local perturbations of the dynamics
(cf.~\cite{HKTP74}).

Regarding the question of the existence of KMS states for a given $C^*$-action, there are very few general
results, and most analyses are done in particular contexts. In the $C^*$-dynamical case, existence of
KMS states is proven for
approximately inner dynamics if there is a trace state (cf.~\cite{PoSa75}),
time evolutions of Haag--Kastler quantum field theories, satisfying a
nuclearity condition (cf.~\cite{BJ89}), for the Cuntz algebra (cf.~\cite{OP78}),
for the CAR-algebra (cf.~\cite{RST69})
and many others. For a singular action on the Weyl algebra there is an
existence condition in \cite{RST70}.

For a general condition for existence of KMS states,
the only one we know of is by Woronowicz (cf.~\cite{W85}).

\begin{thm}{\rm(Woronowicz)}
Let $(\cA, \R, \alpha)$ be a unital $C^*$-dynamical system.
Then there is a KMS state $\omega$ on $\cA$
 for $\alpha$ at $\beta=1$
if and only if ${\cal L}\not=\cA^{\rm op}\otimes\cA$
(maximal tensor product), where
$\cA^{\rm op}$ is the opposite algebra of $\cA$ and ${\cal L}$ is the smallest closed left ideal
in $\cA^{\rm op}\otimes\cA$ containing the set
\[
\{{A}\otimes\1-\1\otimes\alpha_{i/2}(A^*)\,\mid\,
A\in\cA\;\hbox{an entire element}\}.
\]
\end{thm}
This is \cite[Thm.~3]{W85}.
The set $\fS_\beta$ of  KMS states for $\alpha$ at $\beta$
has an interesting structure.
\begin{thm}
Let $(\cM, \R, \alpha)$ be a  $W^*$-dynamical system. Then
\begin{itemize}
\item[\rm(i)] ${\got S}_\beta\subset\cM_*$ is convex and weakly closed, but need not be compact nor
have extreme points.
\item[\rm(ii)] $\omega\in{\got S}_\beta$ is extremal in ${\got S}_\beta$ if and only if it is a
factor state.
\item[\rm(iii)] Two extremal points of ${\got S}_\beta$ are either equal or disjoint
(i.e.\ have disjoint GNS representations).
\end{itemize}
\end{thm}
See the paragraph below \cite[Prop.~5.3.30]{BR96}. Note that the proofs of (ii) and
(iii)
carry over directly from the corresponding proofs in \cite[Prop.~5.3.30]{BR96}.
If $(\cA, \R, \alpha)$ is a  $C^*$-dynamical system, then far stronger properties
listed in \cite[Prop.~5.3.30]{BR96} hold.

There is a great deal more structure for KMS states, e.g.\ much is known
about the behavior of KMS states
with respect to ~perturbation of $\alpha$ (cf. \cite[Ch.~5.4]{BR96}, and \cite{DJP03}).
We leave this large topic for the monographs.

\section{Ergodic states for $C^*$-actions.}
\mlabel{subsec:3.5}

\begin{defn}
\label{ErgSta}
(a) Let $(\cM,G,\alpha)$ be a $W^*$-dynamical system. We say that
it is {\it ergodic} if $\cM^G = \C \1$
(cf. \cite[Def.~X.3.13]{Ta03}).

(b) For a $C^*$-action $(\cA,G,\alpha)$ a $G$-invariant state
$\omega$ is called {\it ergodic} if it is an extreme point
of the convex set $\fS(\cA)^G$ of all $G$-invariant states.
The state $\omega$ is called {\it weakly ergodic}
if $\cH_\omega^G =\C \Omega_\omega$ holds in the corresponding
covariant GNS representation $(\pi_\omega, U_\omega, \Omega_\omega)$.
\end{defn}

\begin{rem} \mlabel{rem:b.2}
(a) For a $C^*$-action $(\cA,G,\alpha)$, if  $\fS(\cA)^G\not=\emptyset$,
then extreme points, i.e. ergodic states,  exist: First, if $\cA$ is unital,
then the state space  $\fS(\cA)$ is weak-$*$-compact, and it is easy to see
that the subspace of invariant states  $\fS(\cA)^G\subset\fS(\cA)$ is weak-$*$-closed,
hence is also weak-$*$-compact and convex, and nonempty.
 It follows from the Krein--Milman theorem that
$\fS(\cA)^G$ has extreme points, and $\fS(\cA)^G$ is equal to the closed convex set they generate.
So ergodic states exist.

 If $\cA$ is nonunital,
then augment $\cA$ with the identity to obtain $\tilde\cA$.
It contains the maximal ideal $\cA$ and $\tilde\cA/\cA\cong \C$.
Extend the action $\alpha$
to $\tilde\cA$ by setting $\tilde\alpha_g(\1)=\1$ for all $g\in G$.
We identify the set $\fS(\cA)$ of states of $\cA$ with those
states $\omega$ of $\tilde\cA$ for which $\omega\restriction\cA$
is a state of~$\cA$, so that $\omega$ is uniquely determined
by this restriction. Each state $\omega$
on $\tilde\cA$ has a unique decomposition
\[
\omega=\lambda\omega_0 +(1-\lambda)\varphi
\quad \mbox{ with }  \lambda = 1 -
\|\omega\restriction\cA\| \in[0,1], \]
where $\omega_0$ is the unique state satisfying $\omega(\cA)=0$
and $\varphi \in \fS(\cA)$.
Then $\omega$ is invariant
if and only if  $\varphi\in\fS(\cA)^G$, so that
\[ \fS(\tilde A)^G = \conv\big(\{\omega_0\} \cup \fS(\cA)^G)\big) \]
and therefore we have
\[ \Ext(\fS(\tilde A)^G)
= \{\omega_0\} \dot\cup \Ext\big(\fS(\cA)^G\big), \]
where $\Ext(C)$ denotes the set of extreme points of the convex set~$C$.
Here the inclusion $\subeq$ is immediate and for the converse
we use that $\fS(\cA)$ is a face of the convex set $\fS(\tilde\cA)$
which follows from the convexity of the functional
$\omega \mapsto \|\omega\restriction\cA\|$.
This describes the ergodic states of $(\tilde\cA,G,\tilde\alpha)$
in terms of those of $(\cA,G,\alpha)$.

(b)  If $(\cM,G,\alpha)$ is a $W^*$-dynamical system and
a normal state $\omega$ is an extreme point in $ \fS_n(\cM)^G$, then it also is
an extreme point in the larger set $\fS(\cM)^G$ of all $G$-invariant
states of the $C^*$-algebra~$\cM$. This is due to the fact that
$\fS_n(\cM) \subeq \fS(\cM)$ is a face, which in turn follows from
the continuity characterization in \cite[Thm.~2.4.21]{BR02}.
Conversely, if it is an extreme point of $\fS(\cM)^G$ it is an extreme point in
$\fS_n(\cM)^G$. Hence a normal state is ergodic if and only if it is extreme in
the set of invariant normal states $\fS_n(\cM)^G$.

(c) If $\omega$ is an ergodic state of  a $C^*$-action $(\cA,G,\alpha)$,
then the associated  $W^*$-dynamical system $(\pi_\omega(\cA)'',G,\tilde\alpha)$
need not be ergodic, though the converse is true.
For instance, if $G=\{\1\}$ or, more generally,
if $G$ is arbitrary and its action on $\cA$ is trivial,
then the ergodic states of $(\cA,G,\alpha)$ are exactly the pure states of~$\cA$, and for every pure state $\omega$ of $\cA$ one has
$(\pi_\omega(\cA)'')^G=\pi_\omega(\cA)''=B(\cH_\omega)$.
Hence the $W^*$-dynamical system $(\pi_\omega(\cA)'',G,\tilde\alpha)$ is not ergodic unless $\dim\cH_\omega=1$.
Examples of this type can also be constructed for nontrivial  group actions, cf. Example~\ref{ex:ergdis} below.
This discrepancy between ergodicity of the state $\omega$
and ergodicity of the $W^*$-dynamical system $(\pi_\omega(\cA)'',G,\tilde\alpha)$ is discussed in Theorem~\ref{thm:b.3} below.

(d) Ergodic states for singular actions need not have covariant GNS representations, unlike ground states and KMS states,
so are less useful.
To get a covariant GNS representation, one needs also a condition in Proposition~\ref{fixstatecov}. It seems for
singular actions this must be added to obtain useful ergodic states. We now
give an example of an ergodic state
where the GNS-representation is not covariant.
\end{rem}

\begin{ex}
We continue the context of Example~\ref{ExoticEx}.
Let $G$ be an abelian exotic topological group.
Take the left regular representation on $\ell^2(G)$, i.e.\ $(V_g\psi)(h):= \psi(g^{-1}h)$ for $\psi\in \ell^2(G)$,
$g,\,h\in G$.
Let $\cA=\cK(\ell^2(G))$ which is a simple ideal of  $\cB(\ell^2(G))$.
Define $\alpha:G\to \Aut(\cA)$ by $\alpha_g(A):=V_gAV_g^*$.
We showed above that the $C^*$-action $(\cA,G,\alpha)$ has no covariant representations,
so it suffices to show that it has ergodic states. As $G$ is abelian, it is amenable (with respect to  any topology), hence
$(\cA,G,\alpha)$ has an invariant state, i.e.~$\fS(\cA)^G\not=\emptyset$.
By (a) above, it has ergodic states.
\end{ex}
\begin{thm} \mlabel{thm:b.3} Let $(\cA,G,\alpha)$ be a $C^*$-action,
$\omega \in \fS(\cA)^G$ and $(\pi_\omega, U_\omega, \cH_\omega, \Omega_\omega)$ be the
corresponding covariant GNS representation.
Consider the following properties:
\begin{itemize}
\item[\rm(a)] $(\pi_\omega(\cA)'')^G = \C \1$, i.e., the action of
$G$ on $\pi_\omega(\cA)''$ is ergodic.
\item[\rm(b)] $\cH_\omega^G = \C \Omega_\omega$, i.e.,
$\omega$ is weakly ergodic.
\item[\rm(c)] $\omega$ is $G$-ergodic, i.e., an extreme point of $\fS(\cA)^G$.
\item[\rm(d)] $\pi_\omega(\cA)\cup U_\omega(G)$ acts irreducibly on $\cH_\omega$.
\item[\rm(e)] $\pi_\omega(\cA)''$ is of type III or $\Omega_\omega$
is a trace vector for $\pi_\omega(\cA)'$.
\end{itemize}
Then the implications {\rm(a)} $\Rarrow$ {\rm(b)} $\Rarrow$ {\rm(c)}
$\Leftrightarrow$ {\rm(d)} and  {\rm(b)} $\Rarrow$ {\rm(e)} hold.
Moreover, {\rm(a)} implies that $\Omega_\omega$ is separating for
$\pi_\omega(\cA)''$.
On the other hand,  if  $\Omega_\omega$ is a separating vector for
$\pi_\omega(\cA)''$, then the four conditions {\rm (a)--(d)} are equivalent.
\end{thm}
The relations between {\rm(a)} to {\rm(d)} are in \cite[Thm.~4.3.20]{BR02},
whereas the implication {\rm(b)} $\Rarrow$ {\rm(e)} is in \cite[Thm.~1]{Lo79},
which is a Theorem by Hugenholtz and St\o{}rmer
(cf.  \cite{Hu67, St67}).

\begin{lem} \mlabel{lem:commute}
Let $\cM \subeq \cB(\cH)$ be a von Neumann algebra and
$G \subeq \cU(\cH)$ be a subgroup normalizing~$\cM$.
Suppose further that $\Omega \in \cH^G$ is a $G$-invariant
cyclic separating unit vector for $\cM$. Then $G$ commutes with
the corresponding modular objects $J$ and $\Delta$.
\end{lem}

\begin{prf} Denote the action of $G$ on $\cM$ by
$\alpha_g(M) := g M g^*$. As $G$ fixed $\Omega$, we have
$g(M\Omega) = \alpha_g(M)\Omega,$
and this implies that the unbounded antilinear involution
defined by $S(M\Omega) := M^*\Omega$ for $M \in \cM$
commutes with $G$. Now $J$ and $\Delta$ are uniquely determined
by the polar decomposition $S = J \Delta^{1/2}$, hence also commute
with~$G$.
\end{prf}

The following theorem is a refinement of the preceding
one for von Neumann algebras with a cyclic vector $\Omega$. It clarifies
in particular to which extent (c) implies (a), resp., (b). Note
 that $\Omega$ is separating if and only if $p = s(\omega) = \1$.

\begin{thm} Let $\cM \subeq \cB(\cH)$ be a von Neumann algebra,
$G \subeq \cU(\cH)$ be a subgroup normalizing $\cM$,
$\Omega \in \cH^G$ be an $\cM$-cyclic vector and $\omega \in \fS_n(\cM)^G$
be the corresponding state. We write $p = s(\omega) \in \cM$ for its
carrier projection. Then the following are equivalent:
  \begin{itemize}
  \item[\rm(i)] $\cM\cup U_\omega(G)$ acts irreducibly on $\cH_\omega$.
  \item[\rm(ii)] $\omega$ is a $G$-ergodic state of $\cM$.
\item[\rm(iii)] $(\cM')^G = \C \1$, i.e., the $G$-action on $\cM'$ is ergodic.
  \item[\rm(iv)] $\omega$ is a $G$-ergodic state of $\cM_p=p\cM p$.
  \item[\rm(v)] $(\cM_p)^G = \C\1$, i.e., the $G$-action on $\cM_p$ is ergodic.
  \item[\rm(vi)] $\cH_p^G = \C \Omega$, i.e., the state $\omega\res_{\cM_p}$
is weakly ergodic.
  \end{itemize}
\end{thm}

\begin{prf} (i) $\Leftrightarrow$ (ii) follows from
the equivalence of (c) and (d) in Theorem~\ref{thm:b.3} and
Remark~\ref{rem:b.2}.

(i) $\Leftrightarrow$ (iii) follows from $(\cM \cup G)' = (\cM')^G$.

(iii) $\Leftrightarrow$ (iv): As $G$ fixes $\omega$, it commutes with $p$,
hence leaves the subspace $\cH_p := p\cH$ invariant.
The cyclic representation of $\cM_p$ on the subspace $\cH_p$
has the commutant $(\cM_p)'=(\cM')_p$ (cf. Lemma~\ref{lem:proj}).
Since $\Omega$ is cyclic for $\cM$, it is separating for $\cM'$,
and thus $(\cM')_p \cong \cM'$.
Therefore the equivalence
of (iii) and (iv) follows by applying the equivalence of (i) and (iii)
to $\cM_p$ instead of~$\cM$.

(iii) $\Leftrightarrow$ (v): As the representation of
$\cM_p$ on $\cH_p$ is standard by Lemma~\ref{lem:3.1},
the corresponding conjugation $J$ yields an antilinear
$G$-equivariant  bijection $\cM_p \to \cM'$.
Here the $G$-equivariance follows from the fact that, on $\cH_p$,
the $G$-action commutes with $J$ by Lemma~\ref{lem:commute}.
Hence (iii) and (v) are equivalent.

(v) $\Leftrightarrow$ (vi) follows from Theorem~\ref{thm:b.3}
because $\Omega$ is a separating cyclic vector for~$\cM_p$.
\end{prf}

\begin{rem} We have seen above that a weakly ergodic state is
in particular ergodic. So it is natural to look for sufficient conditions
for the converse to be true.
Suppose that $\cA$ is a separable $C^*$-algebra,
$G$ locally compact separable and
$(\cA,G,\alpha)$ a $C^*$-dynamical system. Then
$\cA$ is $G$-abelian (i.e.\ $\fS(\cA)^G$ is a simplex)
if and only if every  invariant
ergodic state $\omega \in \fS(\cA)$ is weakly ergodic (\cite[Thm.~2]{DNN75}).
\end{rem}

\begin{prop} \mlabel{prop:d.2}
Let $(\cM,\cH,J,\cC)$ be a von Neumann algebra in standard form,
identify $\Aut(\cM)$ with $\cU(\cH)_\cM$ and consider
a subgroup $G \subeq \cU(\cH)_\cM$.
The following are equivalent
  \begin{itemize}
  \item[\rm(i)] $\cM\cup G$ acts irreducibly on $\cH$.
  \item[\rm(ii)] $(\cM')^G = \C \1$.
  \item[\rm(iii)] $\cM'\cup G$ acts irreducibly on $\cH$.
  \item[\rm(iv)] $\cM^G = \C \1$.
  \end{itemize}
\end{prop}

\begin{prf} Conjugating with $J$ implies the equivalence of (i)/(iii),
(ii)/(iv).
The equivalence between (i) and (ii) and of (iii) and (iv)
follows from Schur's Lemma.
\end{prf}

\begin{rem} Suppose that $(\cM,G, \alpha)$ is a $W^*$-dynamical
system where $\cM$ is commutative and
$\omega \in \fS_n(\cM)$ is a faithful separating normal state.
Then $\cM$ is countably decomposable, hence isomorphic to
$L^\infty(X,\fS,\mu)$ for a finite measure space.
Then $\cM_* \cong L^1(X,\fS,\mu)$ and the standard representations
can be realized on $\cH := L^2(X,\fS,\mu)$. The group
$G$ acts on this space by
\[ U_{g}f = \delta(g)^{1/2} (g_*f),\]
where $\delta(g) \in L^1(X,\fS,\mu)$ is the Radon--Nikodym derivative
defined by $g_*\mu = \delta(g)\mu$.
Note that the implementability of $G$ on the measurable
space $(X,\fS)$ may be problematic if $G$ is not
locally compact second countable, but in any case the unitary
representation on $\cH$ exists and so does the action of
$G$ on the Boolean $\sigma$-algebra $\fS_\mu = \fS/\sim$, where
$E \sim F$ with $\mu(E \Delta F) = 0$. This Boolean $\sigma$-algebra is the space
of projections in $\cM$.

That $\mu$ is ergodic means that $(\cM')^G = \cM^G = \C\1$.
Now $\cH^G \not=\{0\}$ holds only if $[\mu]$ contains a $G$-invariant
finite measure. In fact, $f \in \cH^G$ implies that $|f|^2 \mu$ is $G$-invariant.
For the translation action of $\R$ on itself we have $\cH^G = \{0\}$.
\end{rem}

\appendix

\section{Auxiliary results} \mlabel{app:3}

\begin{lem} \mlabel{lem:3.1b} Let $G$ be a connected topological group acting
on a nonempty set $X$. We consider the corresponding unitary representation
$(\pi, \ell^2(X))$. Then
\begin{itemize}
\item[\rm(i)] every $G$-continuous vector $\xi \in \ell^2(X)$ is fixed, and
\item[\rm(ii)] $\ell^2(X)^G$ is generated by the characteristic functions
of the finite $G$-orbits in $X$.
\end{itemize}
\end{lem}

\begin{prf} (i) Let $\xi \in \ell^2(X)$ be non-zero $G$-continuous vector
and $c \in \C^\times$ be such that $\xi_x = c$ for some $x \in X$.
Then $F_c := \{ x \in X \: \xi_x = c\}$ is a finite subset of $X$.
We write $P_c \: \ell^2(X) \to \ell^2(F_c)$ for the corresponding orthogonal projection.
Let $\eps > 0$ be such that $|\xi_y-c| > \eps$ for $y \not\in F_c$.
If $g \in G$ satisfies $\|P_c (g.\xi - \xi) \| < \eps$,
then, for every $x \in F_c$, we have $|\xi_x - \xi_{g^{-1}.x}| < \eps$,
hence $g^{-1}.x \in F_c$. Now  the finiteness of $F_c$ implies
that $g.F_c = F_c$ and hence $P_c (g.\xi - \xi) = 0$.
We conclude that
\[ U := \{ g \in G \: \|P_c (g.\xi - \xi)  \| < \eps \} \]
is an open closed identity neighborhood of $G$. Since $G$ is connected,
it follows that $G = U$. This shows that all the subsets $F_c$ are $G$-invariant,
and this in turn entails that $\xi$ is fixed under $G$.

(ii) is trivial.
\end{prf}

\begin{lem} \mlabel{lem:3.2b} If $\omega \in \cA^*$ is a tracial state
of a $C^*$-algebra $\cA,$ then
\[ \ker \pi_\omega = \{ A \in \cA \,\mid\, \omega A = 0\}.\]
\end{lem}

\begin{prf} This follows from the fact that the cyclic element
$\Omega \in \cH_\omega$ is also separating: If $\omega A = 0$ and $B \in \cA$, then
\cite[Rem.~3.2.66]{BR02} yields
\[ \la \pi_\omega(AB)\Omega,  \pi_\omega(AB)\Omega \ra
= \omega(B^*A^*AB) = \omega(BB^*A^*A) = 0.\qedhere\]
\end{prf}

\begin{lem} \mlabel{lem:UVW} Let $(U_t)_{t \in \R}$ and $(V_t)_{t \in \R}$
be two commuting continuous unitary one-parameter groups on $\cH$ with
non-negative spectrum, and put $W_t := U_t V_t$.
If $A$ and $B$ are the infinitesimal generators of $U$ and $V$, respectively, then
$A + B$ is closed and the infinitesimal generator of $W$.
\end{lem}

\begin{prf} Decomposing $\cH$ into cyclic subspaces with respect to the representation of
$\R^2$, defined by $(t,s) \mapsto U_t V_s$, we may without loss of generality\ assume that
$\cH = L^2(\R^2, \mu)$ for a finite measure $\mu$ and that
\[ (U_t F)(x,y) = e^{-itx} F(x,y) \quad \mbox{ and } \quad
 (V_s F)(x,y) = e^{-isy} F(x,y).\]
Our assumption now implies that $\supp(\mu) \subeq [0,\infty)^2$.
We further have $(AF)(x,y) = x F(x,y)$ and $(BF)(x,y) = y F(x,y)$. We define
$(CF)(x,y) := (x+y)F(x,y)$ on its maximal domain
\[ \cD(C) := \Big\{ F \in L^2(\R^2, \mu) \,\mid\, \int_{\R^2} (x+y)^2 |F(x,y)|^2\,  d\mu(x,y)
< \infty\Big\}\]
and note that this is the infinitesimal generator of $W$.
Then $\cD(C) = \cD(A) \cap \cD(B)$ follows from the positivity of the functions
$x$ and $y$ $\mu$-almost everywhere.
\end{prf}

\begin{lem} \mlabel{lem:fin}
Let $\cA$ be a unital $C^*$-algebra for which the spectrum of every
hermitian element is finite. Then $\cA$ is finite dimensional.
\end{lem}

\begin{prf} Let $\cC \subeq \cA$ be  maximal abelian. Then $\cC$
inherits the finite spectrum property from $\cA$, and this implies
that $\cC \cong C(X)$, where $X$ is a compact Hausdorff space on which
every continuous function has finitely many values. This implies that
$X$ is finite.

If $|X| = n$, then $\cC$ has a basis $(p_1, \ldots, p_n)$ consisting
of minimal mutually orthogonal projections. Now
\[ \1 = p_1 + \cdots + p_n \quad \mbox{ and } \quad p_i p_j
= \delta_{ij} p_i.\]
This leads to the decomposition
$\cA = \sum_{i,j=1}^n p_i \cA p_j$.
Put $\cA_{ij} := p_i \cA p_j$.
The minimality of each $p_i$ implies that
$\cA_{ii}  = \C p_i$ is one-dimensional.
Now let $i \not=j$ and $0 \not= z \in \cA_{ij}$.
Then $0 \not= zz^* \in \cA_{ii} = \C p_i$. Hence
\[ zw^* := h(z,w) p_i \]
defines a positive definite hermitian form $h$ on $\cA_{ij}$.
If $w \in \cA_{ij}$ is orthogonal to $z$, then
$zw^* = 0$ leads to $zw^*w = 0$. As $w^*w \in \cA_{jj} = \C p_j$
is non-zero if $w\not=0$, it follows that $w^*w =0$.
Therefore $\dim \cA_{ij} = 1$ and thus $\dim \cA \leq n^2$.
\end{prf}

With the preceding lemma one easily verifies the following
(see the proof of \cite[Thm.~1]{CM80}):

\begin{prop}\label{finite_erg}
	Let $\cA$ be a unital $C^*$-algebra and let
$\Gamma \subeq \Aut(\cA)$ be a subgroup which is compact in the
norm topology.
If $\Gamma$ acts ergodically
on $\cA$, i.e., $\cA^\Gamma = \C \1$,  then $\cA$ is finite dimensional.
\end{prop}

\begin{prf} We consider the conditional expectation
\[ f \: \cA \to \C, \quad
f(A)\1 = \int_\Gamma \alpha_\gamma(A)\, d\gamma, \]
where  $d\gamma$ refers to the normalized Haar measure $\mu_\Gamma$
on $\Gamma$, using the assumption that $\cA^\Gamma = \C \1$.

For $\eps \in (0,1)$ we pick an open $\1$-neighborhood
$U \subeq \Gamma$ such that
$\|\alpha_\gamma - \id_\cA\| < \eps$ for $\gamma \in U$.
For $0 \leq A \in \cA$ we then have
\[ f(A) \1
\geq \int_U \alpha_\gamma(A)\, d\gamma
= \int_U (\alpha_\gamma(A)-A)\, d\gamma  + \mu_\Gamma(U)A
\geq 0. \]
As
$\big\|\int_U (\alpha_\gamma(A)-A)\, d\gamma\big\|
\leq \eps \mu_\Gamma(U)  \|A\|,$ this
 leads to
\begin{equation}
  \label{eq:esti2}
f(A)=\|f(A) \1\| \geq \mu_\Gamma(U) \|A\|-\mu_\Gamma(U)\eps \|A\|= c \|A\|
\end{equation}
where  $c := \mu_\Gamma(U) (1 - \eps)$.
If $p_1, \ldots, p_n \in \cA$ satisfy
$0 \leq p_i \leq \1$, $\|p_i\| = 1$, and $\sum_{j = 1}^n p_j = \1$, then
$1 = f(\1) = \sum_{i = 1}^n f(p_i) \geq c n,$
and hence  $n \leq c^{-1}$.
Thus, if $\cC \cong C(X)$ is a commutative subalgebra
of $\cA$, then all partitions of unity of $X$ are finite, and hence
$X$ is finite.  Now the proof of Lemma~\ref{lem:fin}  shows that
$\cA$ is finite dimensional with $\dim \cA \leq c^{-2}$.
\end{prf}

\begin{example}
  \mlabel{ex:ergdis}
Examples of an ergodic state $\omega$  of  a $C^*$-action $(\cA,G,\alpha)$,
where the associated  $W^*$-dynamical system $(\pi_\omega(\cA)'',G,\tilde\alpha)$
need not be ergodic,  for nontrivial  group actions.

Let $(\cA,G,\alpha)$ be a $C^*$-dynamical system where $G$ is a compact group,
and consider the faithful conditional expectation
\[ E\colon\cA\to \cA^G, \qquad
E(A)=\int_G\alpha_\gamma(A)d\gamma,\]
 obtained by averaging with respect to the probability Haar measure on $G$.
Then it is easily checked that $\fS(\cA)^G=\{\omega\in\fS(\cA)\mid \omega\circ E=\omega\}$ and the map
$\fS(\cA^G)\to\fS(\cA)^G$, $\omega_0\mapsto \omega_0\circ E $,
is an affine isomorphism.
Hence the ergodic states of~$\cA$ are exactly the states $\omega=\omega_0\circ E$ where $\omega_0=\omega\vert_{\cA^G}\in\fS(\cA^G)$ is a pure state of $\cA^G$.

For any $\omega=\omega_0\circ E\in\fS(\cA)^G$ with  $\omega_0\in\fS(\cA)$,
the inclusion map $\cA^G\hookrightarrow\cA$ leads to an isometric embedding of Hilbert spaces $\cH_{\omega_0}\hookrightarrow\cH_\omega$
and the corresponding orthogonal projection $P\colon\cH_\omega\to\cH_{\omega_0}$ is the extension by continuity of the conditional expectation $E\colon\cA\to\cA^G$.
Moreover, for every $A\in\cA^G$ one has
$E(AB)=AB$ for all $B\in\cA^G$, hence $P\pi_\omega(A)\vert_{\cH_0}=\pi_{\omega_0}(A)$.
This shows that one has the well-defined surjective linear map
$\pi_\omega(\cA)\to \pi_{\omega_0}(\cA)$,
$T\mapsto PT\vert_{\cH_0}$,
which implies $\dim\pi_{\omega_0}(\cA)\leq\dim\pi_\omega(\cA)$.

If, moreover, the group $G$ is finite and  $\omega\in\fS(\cA)$ is a state whose corresponding $W^*$-dynamical system $(\pi_\omega(\cA)'',G,\tilde\alpha)$ is ergodic,
then $\dim\pi_\omega(\cA)''<\infty$ by Proposition~\ref{finite_erg}, hence
$\dim\pi_{\omega_0}(\cA)<\infty$ by the preceding paragraph.
But at least for the permutation group $G=S_n$, there are many dynamical systems $(\cA,G,\alpha)$ and pure states $\omega_0\in\fS(\cA^G)$ with $\dim\pi_{\omega_0}(\cA)=\infty$,
with $S_n$ acting by permutations on $\cA=\cB^{\otimes n}$
for various $C^*$-algebras $\cB$.
See for instance \cite[Ex. 2.3]{BN16}.
\end{example}

\section{Commutative von Neumann algebras}

Let $(X,\fS,\mu)$ be a $\sigma$-finite measure space.
Then we may identify $L^\infty(X,\fS,\mu)$ with the algebra $\cM$
of multiplication operators on $L^2(X,\fS,\mu)$ and
any function $f \in L^2(X,\fS,\mu)$ for which $f^{-1}(0)$ is a
zero-set is a cyclic separating vector, from which one easily derives
that $\cM = \cM'$ is maximal abelian in $B(\cH)$; in particular $\cM$
is a commutative von Neumann algebra.

The following theorem
provides an effective tool to determine when a $*$-invariant
subset $S \subeq \cM$ generates $\cM$ as a von Neumann algebra,
i.e., $S'' = \cM$. This is achieved by a description of all
von Neumann subalgebras of the von Neumann subalgebra
$\cM = L^\infty(X,\fS,\mu) \subeq \cB(L^2(X,\fS,\mu))$.

\begin{thm} \mlabel{thm:vonNeumann-sigmaalg}
{\rm(The $L^\infty$-Subalgebra Theorem)}
Let $(X,\fS,\mu)$ be a $\sigma$-finite measure space
and $\cA \subeq L^\infty(X,\fS, \mu) \subeq \cB(L^2(X,\fS,\mu))$
be a von Neumann algebra. Then
$$ \fA := \{ E \in \fS \,\mid\, \chi_E \in \cA \} $$
is a $\sigma$-subalgebra of $\fS$ and
\[ \cA \cong L^\infty(X,\fA, \mu\res_{\fA}).\]
Conversely, for every $\sigma$-subalgebra $\fA \subeq \fS$,
$L^\infty(X,\fA, \mu\res_{\fA})$ is a von Neumann subalgebra of
\break $L^\infty(X,\fS, \mu)$.
\end{thm}

\begin{prf} {\bf Step 1:} First we show that $\fA$ is a $\sigma$-algebra.
Clearly $0 \in \cA$ implies
$\eset \in \fA$, and since $\1 \in \cA'' = \cA$, we also have
$\chi_{E^c} = \1 - \chi_E \in \cA$ for each $E \in \fA$.
From $\chi_E \cdot \chi_F = \chi_{E \cap F}$ we derive
that $\fA$ is closed under finite intersections.
Now let $(E_n)_{n \in \N}$ be a sequence  of elements in $\fA$.
It remains to show that $E :=  \bigcap_{n \in \N} E_n \in \fA$. Let
$F_n := E_1 \cap \cdots \cap E_n$. Then $F_n \in \fA$
implies $\chi_{F_n} \in \cA$. Moreover,
$\chi_{F_n}\to \chi_F$ holds pointwise, so that $\chi_{F_n} \to \chi_F$
in the weak operator topology, so that $\chi_F \in \cA$
and thus $F \in \fA$. This proves that $\fA$ is a $\sigma$-algebra.

{\bf Step 2:} That $\cA \supeq L^\infty(X,\fA, \mu\res_\fA)$ follows directly from the fact
that $\cA$ contains all finite linear combinations $\sum_j c_j \chi_{E_j}$, $E_j \in \fA$,
the norm-closedness of $\cA$ and the fact that every element
$f \in L^\infty(X,\fA,\mu\res_{\fA})$ is a norm-limit of a sequence of step
functions~$f_n$.

{\bf Step 3:} Finally we show that $\cA \subeq L^\infty(X,\fA, \mu\res_\fA)$, i.e., that
all elements of $\cA$ are $\fA$-measurable (if possibly modified on sets of measure
zero).

Note that $\cA$ is closed under bounded pointwise limits.
Let $(p_n)_{n\in\N}$ be the sequence of polynomials
converging on $[0,1]$ uniformly  to the square
root function. For $0 \not= f\in \cA$, we consider the functions
$p_n\big(\frac{|f|^2}{\|f\|_\infty^2}\big)$,
which also belong to $\cA$. Since they converge pointwise
to $\frac{|f|}{\|f\|_\infty}$, we see that
$|f| \in \cA$. For real-valued elements $f,g \in \cA$, this further implies that
\[   \max(f,g)=\frac{1}{2}(f+g+|f-g|) \in \cA.\]
For any $c \in \R$, it now follows that
$\max(f,c) \in \cA$. The sequence
$e^{-n (\max(f,c)-c)} \in \cA$ is bounded and converges
pointwise to the characteristic function $\chi_{\{f \leq c\}}$ of the set
\[ \{f \leq c\} := \{ x \in X \,\mid\,  f(x) \leq c\}. \]
We thus obtain that $\chi_{\{f \leq c\}} \in \cA$. We conclude that the set
$\{f \leq c\}$ is contained in the $\mu$-completion $\fA_\mu$ of $\fA$,
and this finally shows that
$f \in L^\infty(X,\fA_\mu,\mu)= L^\infty(X,\fA,\mu)$.

{\bf Step 4:} To show the converse,
let $\fA \subeq \fS$ be a $\sigma$-subalgebra, and consider
the closed subspace $\cH_\fA := L^2(X,\fA, \mu) \subeq L^2(X,\fS, \mu)$
generated by the characteristic functions $\chi_E$, $E \in \fA$,
$\mu(E) < \infty$.
Then a projection operator defined by a characteristic function
$\chi_E \in L^\infty(X,\fS, \mu)$ preserves $\cH_\cA$ if and only if
$E \in \fA_\mu$. Therefore
$L^\infty(X,\fA, \mu\res_{\fA}) = L^\infty(X,\fA_\mu, \mu\res_{\fA})$
is the von Neumann subalgebra of
$L^\infty(X,\fS, \mu)$ consisting of operators preserving~$\cH_A$.
\end{prf}

\begin{cor} \mlabel{cor:genvonneumann}
If $(X,\fS,\mu)$ is a $\sigma$-finite measure space
and
$\cF \subeq L^\infty(X,\fS,\mu)$
is a subset with the property that
$\fS$ is the smallest $\sigma$-algebra for which all elements of $\cF$  are
measurable, then $\cF'' = L^\infty(X,\fS,\mu)$, i.e.,
$\cF$ generates $L^\infty(X,\fS,\mu)$ as a von Neumann algebra.
\end{cor}

\begin{prf} We have seen in Theorem~\ref{thm:vonNeumann-sigmaalg}  that
$\cF'' = L^\infty(X,\fA, \mu\res_{\fA})$ holds for a
$\sigma$-subalgebra $\fA \subeq  \fS$. Then
all elements of $\cF$ are measurable with respect to the $\mu$-completion
$\fA_\mu$ of $\fA$, so that $\fS \subeq \fA_\mu$. This implies that
\[ \cF''
= L^\infty(X,\fA, \mu\res_{\fA})
= L^\infty(X,\fA_\mu, \mu\res_{\fA})
\supeq L^\infty(X,\fS, \mu).\qedhere\]
\end{prf}

\section{A corrigendum to \cite{Ne14}}

In this short section we provide a corrigendum
for a few wrong statements in \cite{Ne14} which have
no consequences in that paper.

We consider a $C^*$-action $(\cA,G,\alpha)$.
In the introduction of \cite{Ne14} and in \cite[p.~314]{Ne14}
we say that in \cite{Bo83}
a state $\omega \in \fS(\cA)$ occurs in a covariant representation if and only if
$\omega \in (\cA^*)_c$. This is not correct in general and rectified by
Theorem~\ref{thm:2.7}, but it is ok for $C^*$-dynamical systems
(Corollary~\ref{cor:2.21}). We need, in addition, that
$\cA\omega\cA \subeq(\cA^*)_c$.

Note also that \cite[Cor.~6.3(ii)]{Ne14} is correct because there it is assumed that
the action of $T$ on $G$ is continuous.

\section{A corrigendum to \cite{GrN14}}

In \cite[Thm.~II.3]{Bo69} Borchers states conditions
which imply that  $\fS(\cA)_c$ is a folium, but there he assumes the
regular case. This leads to a false statement in \cite[Prop.~8.9(ii)]{GrN14},
where it is claimed that $(\cA^*)_c
= \pi_{co}(\cA)''_*$. In general this is false by Example~\ref{AcNotFol}.

\section{Index of terms and notation}

\begin{minipage}[t]{0.5\textwidth}
\fontsize{8pt}{10pt}\selectfont
${\cal A}_c$, Remark~\ref{cont}\\
 $\cA^*,$ Sect.~\ref{NotTerm}\\
 $(\cA^*)_c$,   Def.~\ref{def:1.1}(iii)\\
 Arveson spectrum of $\alpha$, ${\rm Spec}(\alpha)$, Def.~\ref{def:arveson}(1)\\
  Arveson spectral subspace of $\alpha$, $\cM^\alpha(S)$, Def.~\ref{def:arveson}(2)\\
 $\cB(X,Y)$, Sect.~\ref{NotTerm}\\
$ \cB_2(\cH)$, Hilbert--Schmidt operators, Ex.~\ref{ex:3.4}\\
  Borchers--Arveson group, Def.~\ref{def:minimal}(b)\\
 Borchers--Arveson projections $Q_t$, Def.~\ref{def:minimal}(b)\\
Borchers--Arveson Theorem~\ref{BA-thm}\\
Borchers--Halpern Theorem~\ref{thm:bor}\\
Borchers' Theorem on modular inclusions~\ref{BorTMI}\\
$C^*$-action,  $(\cA, G, \alpha)$,
Def.~\ref{def:1.1}(i)\\
$C^*$-dynamical system, Def.~\ref{def:1.1}(ii)\\
$C$-spectral condition, Def.~\ref{DefSpecCon}\\
carrier projection of a state,  $s(\omega)$, Def.~\ref{def:carrier}\\
carrier projection of a vector, $s(\Omega)$, Def.~\ref{standP}\\
central support  of a folium, $z(F),$ Remark~\ref{rem:2.x}(3)\\
central support of a projection, $z(P)$, Lemma~\ref{CP1cyclic}(i)\\
central support of a state,  $z(\omega)$, Def.~\ref{def:carrier}\\
central support of a representation, $z(\pi),$ Remark~\ref{rem:2.x}(3)\\
countably decomposable, above Prop.~\ref{prop:2.9} \\
covariant representation $(\pi,U)$,  Def.~\ref{def:2.2}(a)\\
covariant state,  Def.~\ref{def:2.2}(c)\\
Doplicher  ideal, Rem.~\ref{DopDef}\\
dual action $\alpha^*:G\to\cB(\cA^*)$, Def.~\ref{def:1.1}(iii)\\
ergodic  $W^*$-dynamical system,  Def.~\ref{ErgSta}(a)\\
ergodic state, Def.~\ref{ErgSta}(b)\\
exotic group, Example~\ref{ExoticEx}\\
folium, $F(\pi)$, Def.~\ref{FolDef}, Remark~\ref{rem:2.15}(b)\\
$\Fol(\omega)$, $\Fol_G(\omega)$, Theorem~\ref{thm:2.7}\\
$G_d$, Sect.~\ref{NotTerm}\\
generating projection, below Lemma~\ref{CP1cyclic}\\
ground state, Def.~\ref{defgroundst0}\\
ground state vectors,  Def.~\ref{def:groundvec}\\
induced von Neumann algebra $(\cM')\s P.$, Def.~\ref{DefRed}(i)\\
inner covariant representation,   Sect.~\ref{InnCovRep}\\
inner minimal positive one-parameter group, Def.~\ref{def:minimal}(a)\\
KMS condition,  Def.~\ref{defKMS}\\
KMS state,  Def.~\ref{defKMS}\\
Kallman's Theorem~\ref{thm:Kallman}\\
Longo's Lemma~\ref{lem:longo}\\
$\cM_*$, Sect.~\ref{NotTerm}\\
$\cM_{co}$,  Def.~\ref{def:1.1d}\\
$\cM_0^\alpha(S)$, Def.~\ref{def:arveson}(2)\\
$\cM$-cyclic,  Def.~\ref{DefRed}(iii)\\
$\cM$-generating,  Def.~\ref{DefRed}(iii)\\
minimal folium $\Fol(E)$, Remark~\ref{rem:2.15}(e)\\
minimal one-parameter group, Def.~\ref{def:minimal}(a)\\
modular condition,  Def.~\ref{defKMS}
\end{minipage}
\begin{minipage}[t]{0.5\textwidth}
\fontsize{8pt}{10pt}\selectfont
modular operator of $\varphi$, $\Delta_\varphi$,  Sect.~\ref{ModGps}\\
modular conjugation of $\varphi$,  $J_\varphi$,  Sect.~\ref{ModGps}\\
modular automorphism group, $(\sigma^\varphi_t)_{t \in \R}$,  Sect.~\ref{ModGps}\\
Murray--von Neumann equivalence, below Remark~\ref{rem:2.x}\\
$\cN_\varphi$  Sect.~\ref{ModGps}\\
$\omega_S$, Lemma~\ref{lem:2.3}\\
opposite algebra,  $\cM^{\rm op}$, Remark~\ref{OppAlg}(1)\\
$P$-standard representation, 
Def.~\ref{PstandRep}\\
$p$-topology  of $ \Aut(\cM)$, Def.~\ref{PUtops}\\
positive spectrum for a representation, Def.~\ref{DefSpecCon}\\
quasi-covariant representation, Def.~\ref{def:2.2}(b)\\
quasi-equivalent representations, Remark~\ref{rem:2.15}(c)\\
quasi-invariant state, Def.~\ref{def:2.2}(c)\\
\strong continuous, Def.~\ref{def:1.1}(ii)\\
reduced von Neumann algebra $\cM_P$, Def.~\ref{DefRed}(i)\\       
reduced action $(\cM_P,G,\beta^P)$,   Sect.~\ref{InnCovRep}\\
\usual case, Def.~\ref{def:1.1}(ii)\\
resolvent algebra $\rsl$, Exmp.~\ref{ccrexmp2} \\
$S_\beta$, horizontal strip,  Def.~\ref{defKMS}\\
 $\fS(\cA),$ $\fS_n(\cM)$ Sect.~\ref{NotTerm}\\
 $\fS_\alpha(\cA)$,  Def.~\ref{def:2.2}(c)\\
 $\fS_\alpha^0(\cA)$, set of ground states,  Sect.~\ref{SectGrndSt}\\
  $\fS(\cA)_c$,   Def.~\ref{def:1.1}(iii)\\
 $\fS_{co}$,  Def.~\ref{def:1.1d}\\
$\fS_{n,\alpha}(\cM)$,  Remark~\ref{remFolia}(c)\\
self-dual cone, Def.~\ref{def:a.1}\\
semifinite weight, Remark~\ref{OppAlg}(2)\\
separating projection, below Lemma~\ref{CP1cyclic}\\
singular action, Def.~\ref{def:1.1}(ii)\\
spectrum of an $A$ with respect to  $\alpha$, $\Spec_\alpha(A)$, Def.~\ref{def:arveson}(1)\\
standard form, $(\cM, \cH, J,\cC)$, Def.~\ref{def:a.1}\\
standard (form) representation, Def.~\ref{def:a.1} \\
standard projection, Def.~\ref{standP}\\
Stinespring dilation (minimal), $(\pi_\phi, \cH_\phi,V_\phi)$,  Def.~\ref{StineDef}\\
strong topology, point-norm topology, Sect.~\ref{NotTerm}\\
subrep. with positive spectrum
 $\pi^{(+)}$, above Prop.~\ref{posrep}\\
support of a state,  $s(\omega)$, Def.~\ref{def:carrier}\\
support of a representation, $s(\pi)$, Def.~\ref{suppDef}\\
$\cU(\cH)_\cM$, Remark~\ref{rem:aut-w*}(a)\\
$u$-topology of $ \Aut(\cM)$, Def.~\ref{PUtops}\\
uniform topology, Sect.~\ref{NotTerm}\\
universal covariant rep.,
$(\pi_{co},U_{co},\cH_{co})$ Def.~\ref{def:1.1d}\\
 $W^*$-dynamical system, Def.~\ref{def:1.1.1}\\
weakly ergodic state, Def.~\ref{ErgSta}(b)\\
 Weyl algebra $\ccr \cH,\sigma.$  Exmp.~\ref{ccrexmp1}
\end{minipage}


\bigbreak

\noindent\textbf{Acknowledgment}.
The work of the first-named author was supported by a grant of the Romanian National Authority for Scientific Research and
Innovation, CNCS--UEFISCDI, project number PN-II-RU-TE-2014-4-0370.

We extend our sincere thanks to the referee for his careful reading
of the manuscript. His suggestions were most helpful in improving the
exposition of this paper.

\end{document}